\documentclass[12pt]{amsart}
\usepackage{amssymb}
\usepackage[noadjust]{cite}
\usepackage{booktabs}
\usepackage{url}
\usepackage{hyphenat}
\usepackage{mathtools}
\usepackage{cancel} 
\usepackage{ifpdf}
\usepackage[T1]{fontenc}
\usepackage[utf8]{inputenc}
\usepackage{mymathabx}
\usepackage{mdwtab}
\usepackage{enumitem}
\ifpdf
  \usepackage[pdftex]{graphicx}
  \usepackage[pdftex,margin=1in]{geometry}
  \usepackage[bookmarks=true, bookmarksopen=true,%
    bookmarksdepth=3,bookmarksopenlevel=2,%
    colorlinks=true,%
    linkcolor=blue,%
    citecolor=blue,%
    filecolor=blue,%
    menucolor=blue,%
    urlcolor=blue]{hyperref}
  \hypersetup{pdftitle={Bordered Floer homology and the branched
      double cover I}}
  \hypersetup{pdfauthor={Robert Lipshitz, Peter S. Ozsváth, and Dylan
      P. Thurston}}
\else
  \usepackage[dvips]{graphicx}
  \usepackage[dvips,margin=1in]{geometry}
  \usepackage[draft]{hyperref}
\fi
\usepackage{color}
\usepackage{tikz}

\usetikzlibrary{matrix,arrows}
\tikzstyle{every picture}=[> = latex']
\tikzset{cdlabel/.style={above,sloped,
    execute at begin node=$\scriptstyle,execute at end node=$}}
\tikzset{algarrow/.style={->, thick}}
\tikzset{blgarrow/.style={->, thick}}
\tikzset{clgarrow/.style={->, thick}}
\tikzset{tensoralgarrow/.style={double, double equal sign distance, -implies}}
\tikzset{tensorblgarrow/.style={->, thin, double}}
\tikzset{tensorclgarrow/.style={->, thin, double}}
\tikzset{modarrow/.style={->, dashed}}
\tikzset{Amodar/.style={->, dashed}}
\tikzset{Dmodar/.style={->, dashed}}

\newcommand{\RR}{\mathbb R}
\newcommand{\CC}{\mathbb C}
\newcommand{\ZZ}{\mathbb Z}

\newcommand{\FF}{\mathbb F}

\newcommand{\co}{\colon}

\newcommand{\OneHalf}{{\textstyle\frac{1}{2}}}

\newcommand{\bdy}{\partial}

\newcommand{\lbracket}{[}
\newcommand{\rbracket}{]}

\newcommand{\Hyph}{\text{-}}

\DeclareMathOperator{\ind}{ind}

\DeclareMathOperator{\gr}{gr}

\DeclareMathOperator{\Cone}{Cone}

\theoremstyle{plain}
\newtheorem{theorem}{Theorem}
\numberwithin{equation}{section}

\newtheorem{proposition}[equation]{Proposition}
\newtheorem{lemma}[equation]{Lemma}
\newtheorem{corollary}[equation]{Corollary}

\newtheorem{definition}[equation]{Definition}

\newtheorem{example}[equation]{Example}

\theoremstyle{definition}

\newtheorem{hypothesis}[equation]{Hypothesis}

\theoremstyle{remark}
\newtheorem{remark}[equation]{Remark}

\newcommand{\HF}{\mathit{HF}}
\newcommand{\HFa}{\widehat {\HF}}

\newcommand{\CFa}{\widehat {\mathit{CF}}}

\newcommand{\x}{\mathbf x}
\newcommand{\y}{\mathbf y}

\newcommand\HH{\mathit{HH}}

\newcommand\Hochschild\HH

\newcommand{\Ainf}{\mathcal A_\infty}
\newcommand{\Alg}{\mathcal{A}}

\newcommand\Blg{\mathcal{B}}
\newcommand\BlgDeg{\mathcal{B}_\degen}
\newcommand\Clg{\mathcal{C}}
\newcommand\Dlg{\mathcal{D}}
\newcommand\Mlg{\mathcal{M}}

\newcommand{\alphas}{{\boldsymbol{\alpha}}}
\newcommand{\betas}{{\boldsymbol{\beta}}}
\newcommand{\gammas}{{\boldsymbol{\gamma}}}

\newcommand{\DD}{\textit{DD}}
\newcommand{\DA}{\textit{DA}}

\newcommand{\AAm}{\textit{AA}} 
\newcommand{\CFD}{\mathit{CFD}}
\newcommand{\CFDD}{\mathit{CFDD}}
\newcommand{\CFDDAa}{\widehat{\mathit{CFDDA}}}
\newcommand{\CFA}{\mathit{CFA}}

\newcommand{\CFDA}{\mathit{CFDA}}
\newcommand{\CFDAa}{\widehat{\CFDA}}
\newcommand{\CFAA}{\mathit{CFAA}}
\newcommand{\CFAAa}{\widehat{\CFAA}}
\newcommand{\CFDa}{\widehat{\CFD}}

\newcommand{\CFDDa}{\widehat{\CFDD}}

\newcommand{\CFAa}{\widehat{\CFA}}

\newcommand{\cZ}{\mathcal{Z}}
\newcommand{\PtdMatchCirc}{\cZ}
\newcommand{\PMC}{\PtdMatchCirc}

\newcommand{\CircPts}{{\mathbf{a}}}

\newcommand{\dg}{\textit{dg} }

\newcommand\Id{\mathbb{I}}
\newcommand\Ground{\mathbf k}
\newcommand\Groundl{\mathbf l}

\newcommand\DT{\boxtimes}
\newcommand\Gen{\mathfrak{S}}

\newcommand\Tensor{\mathcal T}

\newcommand{\Field}{{\FF_2}}
\DeclareMathOperator{\nbd}{nbd}

\newcommand{\Heegaard}{\mathcal{H}}
\newcommand{\HD}{\Heegaard}

\renewcommand{\th}{^\text{th}}

\newcommand{\ModCat}{\mathsf{Mod}}

\DeclareMathOperator{\Mor}{Mor}

\newcommand{\op}{\mathrm{op}}

\newcommand{\Zint}{\mathcal{J}}
\newcommand{\Zpint}{\mathcal{J}'}
\newcommand{\degen}{\diamond}

\makeatletter
\newcommand\honestalg[3]{\bigl\lbracket
\begin{smallmatrix} #1\@ifempty{#3}{}{&#3} \\ #2 \end{smallmatrix}
\bigr\rbracket}

\makeatother
\newcommand{\lab}[1]{$\scriptstyle #1$}

\newcommand{\lsub}[2]{{}_{#1}#2}
\newcommand{\lsup}[2]{{}^{#1}\mskip-.6\thinmuskip#2}

\newcommand\Dehn{D}

\newread\testin

\graphicspath{{draws/}{mpdraws/}{}}
\makeatletter
\def\input@path{{}{draws/}}
\makeatother

\def\mathcenter#1{%
  \vcenter{\hbox{$#1$}}%
}

\def\mfigb#1{
        \mathcenter{\includegraphics[trim=-1 -1 -1 -1]{#1}}
}

\makeatletter
\newcommand\mi@kern[1]{%
  \settowidth\@tempdima{$\mi@obj^{#1}$}
  \kern-\@tempdima
  #1
  \settowidth\@tempdima{$\mi@obj$}
  \kern\@tempdima
}

\newtoks\mi@toksp
\newtoks\mi@toksb
\DeclareRobustCommand{\manyindices}[5]{
  \def\mi@obj{#5}
  \mi@toksp\expandafter{\mi@kern{#2}}
  \mi@toksb\expandafter{\mi@kern{#1}}
  \@mathmeasure4\textstyle{#5_{#1}^{#2}}
  \@mathmeasure6\textstyle{#5_{#3}^{#4}}
  \dimen0-\wd6 \advance\dimen0\wd4
  \@mathmeasure8\textstyle{\hphantom{{}_{#1}^{#2}}#5^{\the\mi@toksp#4}_{\the\mi@toksb#3}}
  \hbox to \dimen0{}{\kern-\dimen0\box8}
}
\makeatother 

\usepackage{ifpdf}
\ifpdf
  \let\textalt\texorpdfstring
\else
  \newcommand{\textalt}[2]{#1}
\fi

\newcommand\rKh{\widetilde{\mathit{Kh}}}

\newcommand\SetS{\mathbf{s}}
\newcommand\SetT{\mathbf{t}}

\begin{document}
\title[Bordered Floer homology and the  branched double cover I]
{Bordered Floer homology and the spectral sequence of a
branched double cover I}

\author[Lipshitz]{Robert Lipshitz}
\thanks{RL was supported by NSF grant DMS-0905796, MSRI, and a Sloan Research
  Fellowship.}
\address{Department of Mathematics, Columbia University\\
  New York, NY 10027}
\email{lipshitz@math.columbia.edu}

\author[Ozsv\'ath]{Peter~S.~Ozsv\'ath}
\thanks{PSO was supported by NSF grant DMS-0804121, MSRI, a Clay
  Senior Scholar Fellowship and a Guggenheim Fellowship.}
\address {Department of Mathematics, Princeton University\\ 
Princeton, NJ 08544}
\email {petero@math.princeton.edu}

\author[Thurston]{Dylan~P.~Thurston}
\thanks {DPT was supported by NSF
  grant DMS-1008049, MSRI, and a Sloan Research Fellowship.}
\address{Department of Mathematics,
         Indiana University\\
         Bloomington, IN 47405}
\email{dpthurst@indiana.edu}

\begin{abstract}
  Given a link in the three-sphere, Z.~Szab\'o
  and the second author constructed a
  spectral sequence starting at the Khovanov homology of the link and
  converging to the Heegaard Floer homology of its branched
  double-cover.  The aim of this paper and its sequel is to explicitly calculate this
  spectral sequence, using  bordered Floer homology.
  There are two primary ingredients in this computation: an explicit calculation
  of filtered bimodules associated to Dehn twists and a
  pairing theorem for polygons.  In this paper we give the first
  ingredient, and so obtain a combinatorial spectral sequence from Khovanov
  homology to Heegaard Floer homology; in the sequel we show that this
  spectral sequence agrees with the previously known one.
\end{abstract} 

\maketitle 

\tableofcontents
\section{Introduction}

Let $L$ be a link in the three-sphere. Two recently-defined
invariants one can associate to $L$ are the reduced Khovanov homology
$\rKh(r(L))$ of the mirror of $L$
and the Heegaard Floer homology $\HFa(\Sigma(L))$ of the
double-cover of $S^3$ branched over $L$. (We take both homology groups with coefficients in
$\Field$.) These two link invariants, $\rKh(r(L))$ and $\HFa(\Sigma(L))$,
have much in common; in particular, they agree when $L$ is an unlink, and they both
satisfy the same skein exact sequence. In~\cite{BrDCov}, this
observation was parlayed into the following:

\begin{theorem} [Ozsv{\'a}th-Szab{\'o}]
  \label{thm:MultiDiagramSpectralSequence}
  For any link $L$ in the three-sphere, there is a spectral sequence
  whose $E_2$ page is $\rKh(r(L))$ and whose $E_\infty$ page
  is $\HFa(\Sigma(L))$.
\end{theorem}

Baldwin~\cite{Baldwin11:ss} has shown that the spectral
sequence appearing in Theorem~\ref{thm:MultiDiagramSpectralSequence}
is itself a link invariant.

The aim of the present paper and its sequel~\cite{LOT:DCov2} is to
calculate the spectral sequence from
Theorem~\ref{thm:MultiDiagramSpectralSequence}, when $L$ is a link
expressed as the plat closure of a braid, using techniques from
bordered Heegaard Floer homology. In the present paper, we achieve the
more modest goal of explicitly describing a spectral sequence from
reduced Khovanov homology to $\HFa(\Sigma(L))$.
As we will see, the existence of a spectral sequence follows from formal properties
of bordered Floer homology; the main contribution of the present paper is to 
compute the spectral sequence.
In the sequel, we identify the spectral sequence from this paper with
the one from Theorem~\ref{thm:MultiDiagramSpectralSequence}.

Our arguments here are a combination of basic topology, formal properties of the theory,
and combinatorics. There is analysis in the background, in knowing
that the invariants are well-defined, and that they satisfy
pairing theorems (i.e., that gluing manifolds corresponds to tensoring
their bordered invariants), but this analysis was all done
in~\cite{LOT1} and~\cite{LOT2}. The main work in
identifying the spectral sequences in the sequel~\cite{LOT:DCov2} is a new
analytic result: a pairing theorem for polygon maps. This
pairing theorem, in turn, follows from a degeneration argument similar to the one used
in~\cite{LOT1} to prove the pairing theorem for modules, but with an
extra step.

\subsection{Description of the spectral sequence}

\begin{figure}
  \centering
  \includegraphics{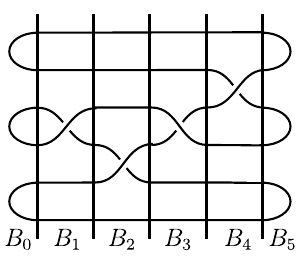}
  \caption{\textbf{Cutting up a braid projection.} The branched double covers of the top and bottom slices are handlebodies; the branched double covers of the other slices are mapping cylinders of Dehn twists.}
  \label{fig:CutUpProjection}
\end{figure}

To describe our spectral sequence, present the link
as the plat closure
of a braid and slice it as
follows.   Think of the link as lying on its side in the
three-ball. If it has $n$ crossings and $2k$ strands, cut the
three-ball along $n+1$ vertical slices. To the left of the
$i=1$ slice, there should be no crossings, just the $k$ extrema in the link
projection. Between the slice $i$ and slice $i+1$ there should be $2k$
strands, exactly two of which should be permuted.
Finally, to the right of the  $i=n+1$ slice, there should be no
crossings, just the $k$ other critical points of the link
projection. See Figure~\ref{fig:CutUpProjection}
for an illustration. (Note that we have departed from standard conventions from
knot theory, where the link is typically thought of as having $k$ maxima and minima,
rather than the $k$ leftmost and $k$ rightmost critical points we have here.)
We
can now think of the three-ball as cut into $n+2$ layers
$B_0,\dots,B_{n+1}$, where the layer $B_i$ lies between the $i\th$
and $(i+1)^{\text{st}}$ slices. The
branched double-covers of $B_0$ and $B_{n+1}$ are handlebodies, while
the branched double-covers of all other $B_i$ (for $i=1,\dots, n$) are
product manifolds $[0,1]\times \Sigma$, where $\Sigma$ is a surface of
genus $k-1$. 

With a little more care, one can pick parametrizations of the boundaries,
making the $B_i$ \emph{bordered} in the sense of~\cite{LOT1}.
We can think of $\Sigma(B_i)$ for $i=1,\dots,n$ as the mapping
cylinder of a diffeomorphism which is a positive or negative Dehn twist along a
(homologically essential) curve $\gamma$ in $\Sigma$. Call this
bordered $3$-manifold $\Dehn_\gamma$ or $\Dehn_\gamma^{-1}$, depending
on the sign of the Dehn twist.

Thinking of $\Sigma$ as a surface equipped with a convenient
parameterization, we can now associate bordered invariants (in the
sense of~\cite{LOT1} and \cite{LOT2}) to each of the
pieces. Specifically, we form $\CFAa(\Sigma(B_0))$,
$\CFDAa(\Sigma(B_i))$ for $i=1,\dots n$, and $\CFDa(\Sigma(B_{n+1}))$.
According to a pairing theorem for bordered Floer homology,
the Heegaard Floer complex for $\Sigma(L)$ can be obtained by forming
the tensor products of these bordered pieces, i.e.,
there is a homotopy equivalence
\begin{equation}
  \label{eq:PairingTheorem}
  \CFa(\Sigma(L))=\CFAa(B_0)\DT \CFDAa(B_1)\DT\cdots\DT \CFDAa(B_n)\DT \CFDa(B_{n+1}).
\end{equation}

We now turn to studying the bimodules associated to the Dehn twists
appearing in Equation~\eqref{eq:PairingTheorem}
(i.e., the modules $\CFDAa(B_i)$ for $i=1,\dots,n$). First we have the following:

\begin{definition}
  Let $\gamma\subset \Sigma$ be an embedded curve.
  Let $Y_{0(\gamma)}$ denote the three-manifold obtained as $0$-framed surgery on 
  $[0,1]\times \Sigma$ along
  $\gamma$, thought of as supported in $\{\OneHalf\}\times \Sigma$
  and equipped with its surface framing.
\end{definition}

\begin{figure}
  \centering
  \includegraphics{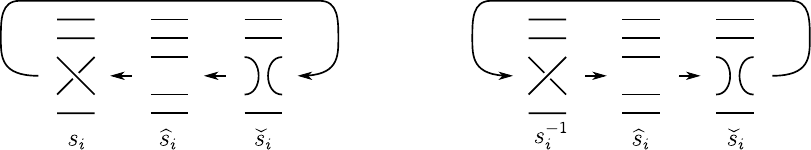}
  \caption{\textbf{Braid-like and anti-braid-like resolutions of braid
      generators.} Left: a positive braid generator $s_i$, its
    braid-like resolution, and its anti-braid-like resolution. Right:
    a negative braid generator $s_i^{-1}$, its braid-like
    resolution, and its anti-braid-like resolution. In both cases, the
    arrows indicate the directions in the skein sequence.}
  \label{fig:resolve-braid-gen}
\end{figure}

Given a braid generator $s_i$, by the \emph{braid-like
  resolution} $\widehat{s}_i$ we simply
mean the identity braid, and by the \emph{anti-braid-like resolution}
$\widecheck{s}_i$ we mean a cap followed by a
cup in the same column as $s_i$; see
Figure~\ref{fig:resolve-braid-gen}.
In these terms, the (unoriented) skein exact
triangles for Khovanov homology and $\HFa(\Sigma(L))$ have the forms
\[
\begin{tikzpicture}
  \node at (0,0) (tc) {$s_i$};
  \node at (-1, -1) (bl) {$\widehat{s}_i$};
  \node at (1, -1) (br) {$\widecheck{s}_i$};
  \draw[->] (tc) to (br);
  \draw[->] (br) to (bl);
  \draw[->] (bl) to (tc);
\end{tikzpicture}
\qquad\qquad
\begin{tikzpicture}
  \node at (0,0) (tc) {$s_i^{-1}$};
  \node at (-1, -1) (bl) {$\widehat{s}_i$};
  \node at (1, -1) (br) {$\widecheck{s}_i$};
  \draw[->] (br) to (tc);
  \draw[->] (bl) to (br);
  \draw[->] (tc) to (bl);
\end{tikzpicture}.
\]

The branched double cover of $s_i$ is $\Dehn_\gamma$ (for an
appropriate non-separating curve $\gamma$); the branched double cover of
$\widehat{s}_i$ is the identity cobordism $\Id$; and the
branched double cover of $\widecheck{s}_i$ is $Y_{0(\gamma)}$
(for the same curve~$\gamma$).  So, it is natural to
expect that there would be an exact triangle of bordered Floer
bimodules of these manifolds. The next theorem says that this is the
case.
\begin{theorem}
  \label{thm:MappingCone}
  There are distinguished bimodule morphisms $F^-\co
  \CFDDa(\Id)\to\CFDDa(Y_{0(\gamma)})$ and 
  $F^+\co
  \CFDDa(Y_{0(\gamma)})\to\CFDDa(\Id)$
  which are uniquely characterized
  up to homotopy by the following properties:
  \begin{itemize}
  \item the mapping cones of $F^+$ and $F^-$ are identified
    with $\CFDDa(\Dehn_\gamma)$ and $\CFDDa(\Dehn_\gamma^{-1})$
    respectively, and
  \item the maps $F^\pm$ respect the gradings in a suitable sense (the
    condition of being \emph{triangle-like} of
    Definition~\ref{def:triangle-like}).
  \end{itemize}
\end{theorem}
Note that the order of the two terms is determined by the type of the
crossing, as illustrated in Figure~\ref{fig:resolve-braid-gen}.

The theorem can be divided into two parts. The first is an existence
statement (restated and proved as Theorem~\ref{thm:Dehn-is-MC} below), which is
essentially a local version of the usual exact triangle. The second is the uniqueness statement
which is constructive: the bimodule morphisms $F^\pm$ are calculated
explicitly.
More specifically,
let $B_i^0$ and $B_i^1$ denote the constituent bimodules in the
mapping cone description of $\CFDDa(\Sigma(B_i))$. One of these, the
identity type \DD\ bimodule, is calculated in~\cite{LOT4} (and recalled
in Proposition~\ref{prop:RecallIdentity} below). The other is a
composite of two elementary cobordisms. As such, its information is
formally contained in~\cite{LOT4}; however, we prefer to give a more
explicit version of the answer in
Proposition~\ref{prop:CalculateVertex} below. 
The bimodule morphism
between them is calculated, using only the properties from the statement of
Theorem~\ref{thm:MappingCone}, in
Propositions~\ref{prop:CalculateNegMorphism}
and~\ref{prop:CalculatePosMorphism} (depending on the crossing type).

Tensoring both sides of Theorem~\ref{thm:MappingCone} with the type
\AAm\ module associated to the identity cobordism, $\CFAAa(\Id)$, and
using the fact that this is an equivalence of categories,
gives the analogous statement for $\CFDAa$:
\begin{corollary}
  \label{cor:MappingCone-DA}
  There are bimodule morphisms
  \begin{align*}
  F^-&\co
  \CFDAa(\Id)\to\CFDAa(Y_{0(\gamma)})\\
  F^+&\co
  \CFDAa(Y_{0(\gamma)})\to\CFDAa(\Id)
  \end{align*}
  which are uniquely characterized
  up to homotopy by the property that the mapping cones of $F^+$ and $F^-$
  are identified with $\CFDAa(\Dehn_\gamma)$ and $\CFDAa(\Dehn_\gamma^{-1})$ respectively.  
\end{corollary}
In~\cite{LOT4} we computed the bimodule $\CFAAa(\Id)$. So, computing
the modules and maps in Corollary~\ref{cor:MappingCone-DA} reduces to
computing the modules and maps in Theorem~\ref{thm:MappingCone}.

Theorem~\ref{thm:MappingCone} equips all of the
factors $\CFDAa(B_i)$ for $i=1,\dots,n$ with a two-step filtration.
Thus, the tensor product appearing in
Equation~\eqref{eq:PairingTheorem} is naturally equipped with a
filtration by the hypercube $\{0,1\}^n$.  (The reader should not be
lulled into a sense of complacency, here.  For an ordinary tensor
product, the tensor product of $n$ mapping cones is also a complex
which is filtered by the hypercube; however, the differentials
appearing in this complex never drop filtration level by more than
one. This is not the case with the relevant $\Ainf$ tensor product
appearing in Equation~\eqref{eq:PairingTheorem}, where the
differentials can drop arbitrarily far.)

For $L$ an $n$-crossing link (with the crossings of $L$ enumerated,
say) and $v=(v_1,\dots,v_n)\in\{0,1\}^n$, there is a
corresponding resolution $v(L)$ of $L$, gotten by taking the
$v_i\in\{0,1\}$ resolution at the $i\th$ crossing of $L$, for each
$i$. The reduced Khovanov homology of $r(L)$ is the homology of
an $n$-dimensional hypercube where at each vertex $v$ we place the
group $\rKh(v(r(L)))$. It is easy to see that $\rKh(v(r(L)))\cong
\HFa(\Sigma(v(L))$. (The resolution $v(r(L))$ is an unlink, whose branched
double-cover is a connected sum of $(S^1\times S^2)$'s, a manifold
whose Heegaard Floer homology is easy to calculate.)
These considerations lead to the following:

\begin{theorem}
  \label{thm:BorderedSpectralSequence}
  Equation~\eqref{eq:PairingTheorem} and
  Theorem~\ref{thm:MappingCone} equip a complex for
  $\CFa(\Sigma(L))$ with a filtration by $\{0,1\}^n$.  In particular,
  for each $v\in \{0,1\}^n$, we obtain a chain complex
  $\CFa(\Sigma(v(L)))$, whose homology, $\HFa(\Sigma(v(L)))$, is
  identified with the Khovanov homology of $v(L)$, and whose edge maps
  from $\HFa(\Sigma(v(L)))$ to $\HFa(\Sigma(v'(L)))$ (where $v'$ and
  $v$ differ in one place, with $v_i=0$ and $v_i'=1$) are identified
  with the differentials in Khovanov homology. More succinctly, the
  associated spectral sequence has $E_2$ term identified with
  (reduced) Khovanov homology $\rKh(r(L))$, and $E_\infty$ term identified with
  $\HFa(\Sigma(L))$.
\end{theorem}
Most of Theorem~\ref{thm:BorderedSpectralSequence} follows from the preceding discussion, together with 
the pairing theorem. In particular, this information is sufficient to identify the $E_1$
page of the spectral sequence with Khovanov's (reduced) chain complex,
as $\Field$-vector spaces, and to identify the $E_\infty$ page with $\HFa(\Sigma(L))$.
Identifying the differential on the $E_1$ page with Khovanov's differential uses
a little more geometric input, namely, a version of the
pairing theorem for triangles, which we defer to~\cite{LOT:DCov2}.

Theorem~\ref{thm:BorderedSpectralSequence} gives a spectral sequence of the same form as
Theorem~\ref{thm:MultiDiagramSpectralSequence}, and such that all the
differentials can be explicitly determined.  In the
sequel~\cite{LOT:DCov2}, we prove the following:

\begin{theorem}
  \label{thm:IdentifySpectralSequences}
  The spectral sequence coming from
  Theorem~\ref{thm:BorderedSpectralSequence} is identified with the
  spectral sequence for multi-diagrams from
  Theorem~\ref{thm:MultiDiagramSpectralSequence}. 
\end{theorem}

\subsection{Further remarks}

The present paper is devoted to computing a spectral sequence from
Khovanov homology to the Heegaard Floer homology of a branched double
cover. It is worth pointing out that this spectral sequence has a
number of generalizations to other contexts.  For example,
Roberts~\cite{Roberts13} constructed an analogous spectral sequence from skein
homology of Asaeda, Przytycki, and Sikora~\cite{APS}, converging to knot Floer homology
in a suitable branched cover.  Building on this, Grigsby and Wehrli established an analogous
spectral sequence starting at various colored versions of Khovanov
homology, converging to knot homology of the branch locus in various
branched covers of $L$, leading to a proof that these colored Khovanov
homology groups detect the
unknot~\cite{GrigsbyWehrli:detects}. Bloom~\cite{Bloom:ss} proved an
analogue of Theorem~\ref{thm:MultiDiagramSpectralSequence} using
Seiberg-Witten theory in place of Heegaard Floer homology. More
recently, Kronheimer and Mrowka~\cite{KronheimerMrowka11:detect} have
shown an analogue with $\ZZ$ coefficients, converging to a version of
instanton link homology, showing that Khovanov homology detects the
unknot. 

During the preparation of this paper, Zolt{\'a}n Szab{\'o} announced a
combinatorially-defined spectral sequence starting at Khovanov's 
homology (modulo $2$) and converging to a knot invariant~\cite{Szabo}. It would be interesting
to compare his spectral sequence with the one from 
Theorem~\ref{thm:BorderedSpectralSequence}. 

In this paper we have relied extensively on the machinery of bordered
Floer homology, which gives a technique for computing $\HFa$ for
three-manifolds. Another powerful technique for computing this
invariant is the technique of {\em nice diagrams},
see~\cite{SarkarWang07:ComputingHFhat}. At present, it is not known
how to use nice diagrams to compute the spectral
sequence from $\rKh(r(K))$ to $\HFa(\Sigma(K))$. See
also~\cite{MOS06:CombinatorialDescrip,MOT:grid,ManolescuOzsvath:surgery}
for another approach to the spectral sequence.

In~\cite{LOT4} we computed the mapping class group action in bordered
Floer theory, giving an algorithm for computing $\HFa(Y)$ for any
$3$-manifold~$Y$. The present paper independently gives as a corollary
a part of that computation, for the case where we glue two
handlebodies by a map~$\phi$ that is the double-cover of an element of
the braid group, i.e., $\phi$ is hyperelliptic.

\subsection{Organization}
This paper is organized as
follows. Section~\ref{sec:FilteredComplexes} collects some algebraic
background on the various types of modules, and discusses filtered
versions of them. Section~\ref{sec:SkeinSequence} establishes the skein
sequence for $\CFDDa$ (Theorem~\ref{thm:MappingCone}) and $\CFDAa$
(Corollary~\ref{cor:MappingCone-DA}). Using this, and the pairing
theorem for triangles
from~\cite{LOT:DCov2}, Theorem~\ref{thm:BorderedSpectralSequence} is
also established there. Next, we turn to the work of computing the
spectral sequence. In Section~\ref{sec:diagrams} we compute the type
$D$ invariants of the relevant handlebodies. In
Section~\ref{sec:vertices} we compute the vertices of the cube of
resolutions, and in Section~\ref{sec:edges} we compute the maps
corresponding to edges. These sections compose the heart of the
paper. In Section~\ref{sec:cube} we assemble the pieces to explain how
to compute the spectral sequence, and illustrate how to compute the
spectral sequence with a simple example.

\subsection{Acknowledgements}
The authors thank MSRI for hosting them while most of this research
was undertaken at the program on Homology Theories of Knots and Links
in Spring 2010, and the Columbia mathematics department, where the
rest of the research was completed. We also thank the referee for corrections and helpful suggestions.


\section{Filtered complexes and their products}
\label{sec:FilteredComplexes}

\subsection{Filtered modules}
\label{sec:filtered-modules}

Fix an $\Ainf$-algebra $\Alg$ over a commutative ground ring $\Ground$ of
characteristic $2$. (For the applications in this paper, $\Alg$ will
be the algebra $\Alg(\PMC)$ associated to a pointed matched circle,
and $\Ground$ is the sub-algebra of idempotents, which is isomorphic
to a direct sum $\bigoplus_{i=1}^N\Field$.) A morphism $F\co M_\Alg\to
N_\Alg$ of $\Ainf$-modules over $\Alg$ is just a $\Ground$-module map
$F\co M\otimes \Tensor^*(\Alg[1])\to N$, where 
$$\Tensor^*(V)=\bigoplus_{i=0}^{\infty} \overbrace{V\otimes \dots \otimes V}^i$$
denotes the tensor algebra of a $\Ground$-bimodule $V$. (Here and elsewhere, $\otimes$ denotes the
tensor product over $\Ground$, and $[1]$ denotes a degree shift by $1$.\footnote{In the rest of this paper, we will be working with ungraded modules and chain complexes, but in this background section we will keep track of gradings.})

There is a natural differential on the
set of such maps, making $\Mor(M,N)$ into a chain complex; the grading on this chain complex is defined so that degree $0$ (i.e., grading-preserving) maps $M\otimes \Tensor^*(\Alg[1])\to N$ have grading $0$ in $\Mor(M,N)$. Composition of morphisms is a chain map, so the
category $\ModCat_\Alg$ of right $\Ainf$-modules over $\Alg$ is a
\dg category. The cycles in $\Mor(M,N)$ are the ($\Ainf$-)
\emph{homomorphisms} from $M$ to $N$. See, for instance,~\cite[Section
2]{LOT2} for more details.

\begin{definition}
  Fix an $\Ainf$ algebra $\Alg$, and let $S$ be a finite partially ordered
  set.  An {\em $S$-filtered $\Ainf$-module over $\Alg$} is a 
  collection $\{M^s\}_{s\in S}$ of $\Ainf$-modules over $\Alg$,
  equipped with preferred degree-zero morphisms $F^{s<t}\in \Mor(M^s,M^t[1])$ (for
  each pair $s,t\in S$ with $s<t$)
  satisfying the compatibility equation for each $s<t$:
  $$dF^{s<t} = \sum_{\{u\mid s<u<t\}} F^{u<t}\circ F^{s<u},$$
  where here $d$ denotes the differential on the morphism
  space $\Mor(M^s,M^t)$. Recall that we call an $\Ainf$-module
  $(M^2,\{m_{1+n}\})$ (respectively $\Ainf$-morphism $\{F_{1+n}\}$)
  bounded if all but finitely many of the $m_{1+n}$
  (respectively $F_{1+n}$) vanish. 
  We say that a filtered $\Ainf$-module $M$ is \emph{bounded} if each component
  $M^s$ and each morphism $F^{s<t}$ is bounded.
\end{definition}
  
In particular, the compatibility condition implies
that if $s<t$ are consecutive (i.e. there is no $u$ such that $s<u<t$),
then $F^{s<t}$ is a homomorphism of $\Ainf$-modules.

Consider $M=\bigoplus_{s\in S} M^s$. The higher products $m_n$ on the
$M^s$ and the maps $F^{s<t}$ assemble to give a map $M\otimes
\Tensor^*(\Alg[1])\to M[1]$ given, for $x_s\in M^s$ and
$a_1,\dots,a_n\in\Alg$, by
\[
x_s\otimes a_1\otimes\cdots\otimes a_n\mapsto m^{s}(x_s\otimes
a_1\otimes\dots\otimes a_n)+\sum_{s<t}F^{s<t}(x_s\otimes
a_1\otimes \dots\otimes a_n).
\]
where $m^s\co M^s\otimes \Tensor^*(\Alg[1])\to M^s[1]$ denotes the $\Ainf$-multiplication on $M^s$.
The compatibility conditions on the $F^{s<t}$ can be interpreted as the
condition that this induced map gives $M$ the structure of an 
$\Ainf$-module over $\Alg$.

A \emph{morphism} of $S$-filtered $\Ainf$-modules 
$$\{M^s,F^{s<t}\}\to \{N^s,G^{s<t}\}$$ is a collection of
maps $H^{s\leq t}\in\Mor(M^s,N^t)$ for each pair $s,t\in S$ with
$s\leq t$. The space of morphisms of $S$-filtered $\Ainf$-modules
inherits a differential in an obvious way, making
the category of $S$-filtered $\Ainf$-modules (or type $D$ structures)
into a \dg category. In particular, it makes sense to talk about
$S$-filtered homomorphisms (the cycles), homotopy classes of
$S$-filtered maps (the quotient by the boundaries) and $S$-filtered
homotopy equivalences.

For convenience, assume now that $\Alg$ is a \dg algebra.
Let $\lsup{\Alg}P$ be a (left) type $D$ structure over $\Alg$,
as in~\cite[Chapter 2]{LOT1} or~\cite[Section 2]{LOT2}. That is, $P$
is a $\Ground$-module (which we also denote $P$) together with a
map $\delta^1\co P\to \Alg[1]\otimes P$ satisfying a suitable
compatibility condition. If $P$ and $Q$ are two type $D$ structures,
the morphism space $\Mor(P,Q)$ between them consists of the set of
$\Ground$-module maps from $P$ to $\Alg\otimes Q$, and is
equipped with the differential 
$$df^1 = (m_1\otimes\Id_{Q})\circ f^1 + (m_2\otimes \Id_{Q})\circ (\Id\otimes f^1)
\circ \delta^1_P+  (m_2\otimes \Id_{Q})\circ (\Id_{\Alg}\otimes \delta^1_Q) \circ f^1$$
and composition law, if $g^1\in \Mor(P,Q)$, $f^1\in \Mor(Q,R)$:
$$(f\circ g)^1 = (m_2\otimes \Id_{R})\circ (\Id_\Alg\otimes f)\circ g.$$
(Note that, as in~\cite{LOT2}, we include the superscript $1$
in the notation for morphisms of type $D$ structures.)
These make the collection of left type $D$ structures over $\Alg$ into
a \dg category. 

\begin{definition}
  Let $\Alg$ be a \dg algebra.
  An \emph{$S$-filtered type $D$ structure} is a collection
  $\{P_s\}_{s\in S}$ of type $D$ structures over $\Alg$,
  equipped with preferred morphisms
  $f_{s<t}\co P_s\to \Alg[1]\otimes P_t$ for each pair
  $s,t\in S$ with $s<t$,
  satisfying the following compatibility equation for each $s<t$:
  $$df_{s<t}=\sum_{\{u\mid s<u<t\}} f_{u<t}\circ f_{s<u}.$$
  We say it is \emph{bounded} if all components $P_s$ are bounded type
  $D$ structures.
\end{definition}

Consider $P=\bigoplus_{s\in S} P_s$. Adding up the various $\delta^1_{P_s} \co
P_s\to \Alg[1]\otimes P_s$ and $f^1_{s<t}\co P_s\to \Alg[1]\otimes P_t$, we obtain a map
$$D^1\co P\to \Alg[1]\otimes P.$$
The compatibility condition on the $f^1_{s<t}$ can be interpreted as the
condition that $(P,D^1)$ is a type
$D$ structure.
  
As with $S$-filtered $\Ainf$-modules, $S$-filtered type $D$ structures can be
made into a \dg category in a natural way.

\begin{remark}
  If we were working over an $\Ainf$-algebra $\Alg$, the category of
  type $D$ structures would be an $\Ainf$-category, not a \dg
  category. In this case, the compatibility condition for an $S$-filtered
  type $D$ structure is:
  $$\sum_{s=s_0<\dots<s_n=t}
  m_n(f_{s_{n-1}<s_n},\dots,f_{s_0<s_1})=0,$$ where $m_n$ denotes the
  $n\th$ higher composition on the category of type $D$ structures.
\end{remark}

Fix an $S$-filtered $\Ainf$-module $M=\{M^s\}_{s\in S}$ and an
$S'$-filtered type $D$ structure $P=\{P_{s'}\}_{s'\in S'}$, with at
least one
of $M$ and~$P$ bounded.  Then $S \times S'$ is also a partially
ordered set (with $(s,s') \le (t,t')$ if $s \le t$ and $s' \le t'$),
and we
can consider their box tensor product $M\DT P$ from~\cite[Section~\ref*{LOT:sec:DT}]{LOT1}, equipped with its induced differential, as a
chain complex filtered by $S \times S'$.
Specifically, we have a group 
$$M\DT P = \bigoplus_{s\times s'\in
  S\times S'} M^s\otimes P_{s'}.$$ This is endowed with a differential
$\partial$ which is a sum of maps
\[\partial_{s\leq t}^{s'\leq t'}\co M^s\otimes P_{s'}\to M^t\otimes P_{t'}[1]\]
given by the following expressions (in four slightly different cases):
\begin{figure}
  \centering
  \hspace{-1em}
  \begin{tikzpicture}
    \node at (-1, -1.5) (difftype) {$\partial_{s=t}^{s'=t'}=$};
    \node at (0,1) (tl) {};
    \node at (1,1) (tr) {};
    \node at (1,0) (delta) {$\delta_{s'}$};
    \node at (0,-1) (m) {$m^s$};
    \node at (0,-2) (bl) {};
    \node at (1,-2) (br) {};
    \node at (0,-2.5) (alabel) {(a)};
    \draw[modarrow] (tr) to (delta);
    \draw[modarrow] (delta) to (br);
    \draw[modarrow] (tl) to (m);
    \draw[modarrow] (m) to (bl);
    \draw[tensoralgarrow] (delta) to (m);
    \node at (-1,-6.5) (bdifftype) {$\partial_{s<t}^{s'=t'}=$};
    \node at (0,-5) (btl) {};
    \node at (1,-5) (btr) {};
    \node at (1,-6) (bdelta) {$\delta_{s'}$};
    \node at (0,-7) (bF) {$F^{s<t}$};
    \node at (0,-8) (bbl) {};
    \node at (1,-8) (bbr) {};
    \node at (0,-8.5) (blabel) {(b)};
    \draw[modarrow] (btr) to (bdelta);
    \draw[modarrow] (bdelta) to (bbr);
    \draw[modarrow] (btl) to (bF);
    \draw[modarrow] (bF) to (bbl);
    \draw[tensoralgarrow] (bdelta) to (bF);
  \end{tikzpicture}
  \quad
  \begin{tikzpicture}
    \node at (-3.5,-4) (difftype) {$\partial_{s=t}^{s'<t'}=\hspace{-2em}\displaystyle{\sum_{s'=s_0'<\dots<s'_n=t'}}$};
    \node at (-2,0) (tlblank) {};
    \node at (.25,0) (trblank) {};
    \node at (.25,-1) (delta1) {$\delta_{s_0'}$};
    \node at (.25,-2) (f1) {$f^1_{s_0'<s_1'}$};
    \node at (.25,-3) (delta2) {$\delta_{s_{1}'}$};
    \node at (.25,-4) (rdots) {$\vdots$};
    \node at (.25,-5) (delta3) {$\delta_{s_{n-1}'}$};
    \node at (.25,-6) (f2) {$f^1_{s_{n-1}'<s_n'}$};
    \node at (.25,-7) (delta4) {$\delta_{s_n'}$};
    \node at (-2,-8) (m) {$m^s$};
    \node at (.25,-9) (brblank) {};
    \node at (-2,-9) (blblank) {};
    \node at (-1,-9.5) (clabel) {(c)};
    \draw[modarrow] (tlblank) to (m);
    \draw[modarrow] (m) to (blblank);
    \draw[modarrow] (trblank) to (delta1);
    \draw[modarrow] (delta1) to (f1);
    \draw[modarrow] (f1) to (delta2);
    \draw[modarrow] (delta2) to (rdots);
    \draw[modarrow] (rdots) to (delta3);
    \draw[modarrow] (delta3) to (f2);
    \draw[modarrow] (f2) to (delta4);
    \draw[modarrow] (delta4) to (brblank);
    \draw[tensoralgarrow, bend right=10] (delta1) to (m);
    \draw[tensoralgarrow, bend right=8] (delta2) to (m);
    \draw[tensoralgarrow, bend right=5] (delta3) to (m);
    \draw[tensoralgarrow] (delta4) to (m);
    \draw[algarrow, bend right=9] (f1) to (m);
    \draw[algarrow, bend right=2] (f2) to (m);
  \end{tikzpicture}
  \quad
  \begin{tikzpicture}
    \node at (-3.5,-4) (difftype) {$\partial_{s<t}^{s'<t'}=\hspace{-2em}\displaystyle{\sum_{s'=s_0'<\dots<s'_n=t'}}$};
    \node at (-2,0) (tlblank) {};
    \node at (.25,0) (trblank) {};
    \node at (.25,-1) (delta1) {$\delta_{s_0'}$};
    \node at (.25,-2) (f1) {$f^1_{s_0'<s_1'}$};
    \node at (.25,-3) (delta2) {$\delta_{s_{1}'}$};
    \node at (.25,-4) (rdots) {$\vdots$};
    \node at (.25,-5) (delta3) {$\delta_{s_{n-1}'}$};
    \node at (.25,-6) (f2) {$f^1_{s_{n-1}'<s_n'}$};
    \node at (.25,-7) (delta4) {$\delta_{s_n'}$};
    \node at (-2,-8) (m) {$F^{s<t}$};
    \node at (.25,-9) (brblank) {};
    \node at (-2,-9) (blblank) {};
    \node at (-1.25,-9.5) (dlabel) {(d)};
    \draw[modarrow] (tlblank) to (m);
    \draw[modarrow] (m) to (blblank);
    \draw[modarrow] (trblank) to (delta1);
    \draw[modarrow] (delta1) to (f1);
    \draw[modarrow] (f1) to (delta2);
    \draw[modarrow] (delta2) to (rdots);
    \draw[modarrow] (rdots) to (delta3);
    \draw[modarrow] (delta3) to (f2);
    \draw[modarrow] (f2) to (delta4);
    \draw[modarrow] (delta4) to (brblank);
    \draw[tensoralgarrow, bend right=10] (delta1) to (m);
    \draw[tensoralgarrow, bend right=8] (delta2) to (m);
    \draw[tensoralgarrow, bend right=5] (delta3) to (m);
    \draw[tensoralgarrow] (delta4) to (m);
    \draw[algarrow, bend right=9] (f1) to (m);
    \draw[algarrow, bend right=2] (f2) to (m);
  \end{tikzpicture}
  \caption{\textbf{Differential on the filtered $\DT$ product.} (a)
    Terms with $s=t$ and $s'=t'$. (b) Terms with $s<t$ and $s'=t'$. (c) Terms with $s=t$ and
    $s'<t'$. (d) Terms with $s<t$ and $s'<t'$.}
  \label{fig:filt-diff}
\end{figure}
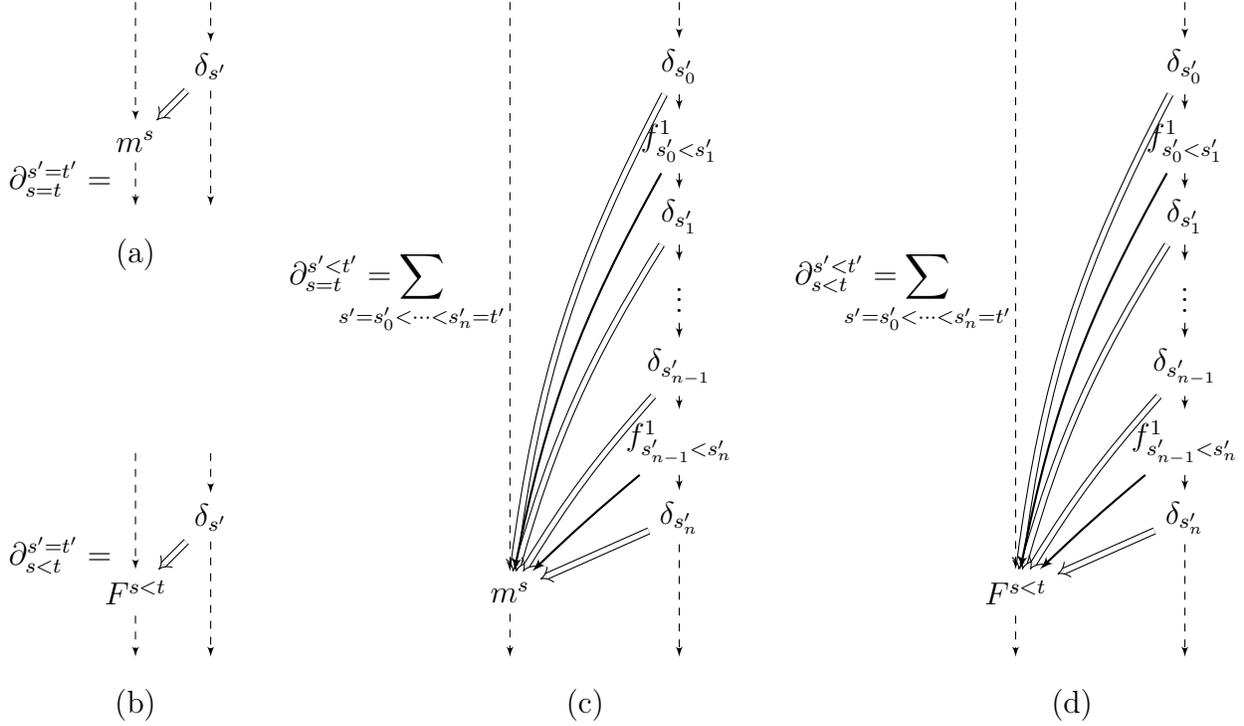
\begin{enumerate}
\item When $s=t$, and $s'=t'$, $\partial_{s=t}^{s'=t'}$ is
simply the usual differential on $M^s\DT P_{s'}$ is given in part (a)
of Figure~\ref{fig:filt-diff}.
Here, the diagram is read as follows. The node $\delta_{s'}$ denotes the map
$P_{s'}\to \Tensor^*(\Alg)\otimes P$ obtained by iterating the map
$\delta^1_{s'}$ an arbitrary number of times. The arrow from
$\delta_{s'}$ to $m^s$ means that we feed this element of
$\Tensor^*(\Alg)$ into the operation $m_s$ on $M^s$. The result is a
map from $M^s\otimes_\Field P_{s'}$ to itself. (A doubled arrow will
always denote an element of a tensor algebra, a solid arrow will
denote a single element of an algebra, and a dashed arrow will denote
an element of a module or type $D$ structure.)
\item When $s<t$ and $s'=t'$, $\partial_{s<t}^{s'=t'}$ is given part (a)
  of Figure~\ref{fig:filt-diff}.
\item When $s=t$ and $s'<t'$, $\partial_{s=t}^{s'<t'}$ is given as a
  sum over all increasing sequences $s'=s'_0<\dots <s'_n=t'$ by part
  (c) of Figure~\ref{fig:filt-diff}.
\item When $s<t$ and $s'<t'$, $\partial_{s<t}^{s'<t'}$ is the sum over
  all increasing sequences $s'=s'_0<\dots <s'_n=t'$ of part (d) of
  Figure~\ref{fig:filt-diff}.  Here, $\delta_{s_i}$ denotes the map
  from $P_{s_i}\to \Tensor^*(\Alg[1])\otimes P_{s_i}[1]$ gotten by
  iterating $\delta^1_{P_{s_i}}$ an arbitrary number of times.
\end{enumerate}

\begin{lemma}\label{lem:DT-descends}
  Let $M$ be an $S$-filtered $\Ainf$-module over a \dg algebra $\Alg$,
  and $P$ be an $S'$-filtered type $D$ structure over $\Alg$.
  Suppose moreover that $M$ is bounded or $P$ is bounded. 
  Then
  the above map $\partial$ endows $M\otimes P$ with the structure of
  an $(S\times S')$-filtered chain complex. Moreover, suppose that $N$ is an $S$-filtered
  $\Ainf$ module which is $S$-filtered homotopy equivalent to $M$,
  and $Q$ is an $S'$-filtered type $D$ structure which is ($S'$-filtered) homotopy
  equivalent to $P$. Suppose that one of the following boundedness conditions
  is satisfied:
  \begin{enumerate}
    \item $M$, $N$ and the homotopy equivalence between them (including its homotopies)
      are bounded.
    \item $P$ and $Q$ are bounded.
  \end{enumerate}
  Then
  the complex $N\DT Q$ is $(S\times S')$-filtered homotopy equivalent to 
  $M\DT P$.
\end{lemma}

\begin{proof}
  We can think of $M$ as an $\Ainf$ module with action given by
  $\sum_s m^s + \sum_{s<t} F^{s<t}$,
  and $P$ as a type $D$ structure with differential
  $$D^1=\sum_{s'} \delta^1_{s'} + \sum_{s'<t'} f^1_{s'<t'}.$$
  The $\DT$ product of these two modules is a chain complex whose
  components are the maps $\partial_{s\leq t}^{s'\leq t'}$ given above.
  Thus, the lemma is essentially a restatement of the fact that $\DT$
  gives a chain complex. The homotopy equivalence property is 
  also a standard property of $\DT$. (See
  \cite[Section~\ref*{LOT:sec:DT}]{LOT1}.)
\end{proof}

\subsection{Generalizations to \textalt{$\DA$}{DA} bimodules}
Fix \dg algebras $\Alg$ and $\Blg$ over ground rings $\Ground$ and
$\Groundl$, respectively.
Recall from \cite[Section~\ref*{LOT2:sec:bimod-var-types}]{LOT2} that an $\Alg\Hyph\Blg$ type \DA\
bimodule is a $(\Ground,\Groundl)$-bimodule $M$ equipped with a map
$$\delta^1\co M\otimes \Tensor^*(\Blg[1])\to \Alg[1]\otimes M,$$
satisfying a structural equation, generalizing the structural
equations for type $D$ structures and $\Ainf$-modules. (Here, tensor
products to the left of $M$ are over $\Ground$, and tensor products to
the right of $M$ are over $\Groundl$.)

Let $\delta^1_0$ denote the restriction of $\delta^1$ to $M=M\otimes
\Blg^{\otimes 0}$ and $\delta^n_0\co M\to \Alg^{\otimes n}\otimes M$
the $n\th$ iterate of $\delta^1_0$.
Under a suitable boundedness condition---for instance, the assumption
that $\delta^n_0$ vanishes for $n$ sufficiently large---by iterating $\delta^1$ we obtain a map
$$\delta\co M\otimes \Tensor^*(\Blg[1])\to \Tensor^*(\Alg[1])\otimes M.$$
The structural equation for a type \DA\ structure is equivalent to the condition that
the map $\delta$ is a chain map.

If $M$ and $N$ are type \DA\ bimodules, a morphism 
$f^1\in\Mor(M,N)$ is a map
$$f^1\co M\otimes \Tensor^*(\Blg[1])\to \Alg\otimes N.$$

Composition is defined by 
\[
(g\circ f)^1(m\otimes
b_1\otimes\cdots\otimes b_n)=\sum_{i=0}^n \mu_\Alg\left((\Id_\Alg\otimes g^1)\left((f^1(m\otimes b_1\otimes\cdots\otimes
b_i)\otimes b_{i+1}\otimes\cdots\otimes b_n)\right)\right).
\]
Here, $\mu_\Alg$ denotes the multiplication on $\Alg$. Graphically,
this is:
\[
\mathcenter{(g\circ f)^1(x\otimes a_1\otimes\cdots\otimes a_n)=}
\mathcenter{
\begin{tikzpicture}
  \node at (0,0) (tc) {$x$};
  \node at (1,0) (tr1) {$a_1$};
  \node at (1.5,0) (tr2) {$\cdots$};
  \node at (2,0) (tr3) {$a_i$};
  \node at (3,0) (tr3b) {$a_{i+1}$};
  \node at (3.5,0) (tr4) {$\cdots$};
  \node at (4,0) (tr5) {$a_n$};
  \node at (0,-2) (deltaN) {$f^1$};
  \node at (0,-3) (f) {$g^1$};
  \node at (-1,-4) (mu) {$\mu_\Alg$};
  \node at (-2,-5) (bl) {};
  \node at (0,-5) (bc) {};
  \draw[modarrow] (tc) to (deltaN);
  \draw[modarrow] (deltaN) to (f);
  \draw[modarrow] (f) to (bc);
  \draw[blgarrow] (tr1) to (deltaN);
  \draw[blgarrow] (tr3) to (deltaN);
  \draw[blgarrow] (tr3b) to (f);
  \draw[blgarrow] (tr5) to (f);
  \draw[algarrow] (f) to (mu);
  \draw[algarrow] (deltaN) to (mu);
  \draw[algarrow] (mu) to (bl);
\end{tikzpicture}}.
\]
As for $\Ainf$-modules and type $D$ structures, the morphism space
$\Mor(M,N)$ between type \DA\ bimodules has a differential, given by:
\begin{align*}
d(f^1)&=\!\!
\mathcenter{\begin{tikzpicture}
  \node at (0,0) (tc) {$x$};
  \node at (1,0) (tr1) {$a_1$};
  \node at (1.5,0) (tr2) {$\cdots$};
  \node at (2,0) (tr3) {$a_n$};
  \node at (0,-1) (f) {$f^1$};
  \node at (-1,-2) (diff) {$d_\Alg$};
  \node at (-2,-3) (bl) {};
  \node at (0,-3) (bc) {};
  \draw[modarrow] (tc) to (f);
  \draw[modarrow] (f) to (bc);
  \draw[blgarrow] (tr1) to (f);
  \draw[blgarrow] (tr3) to (f);
  \draw[algarrow] (f) to (diff);
  \draw[algarrow] (diff) to (bl);
\end{tikzpicture}}
\!\!{+}\!\!
\mathcenter{\begin{tikzpicture}
  \node at (0,0) (tc) {$x$};
  \node at (1,0) (tr1) {$a_1$};
  \node at (1.5,0) (tr2) {$\cdots$};
  \node at (2,0) (tr3) {$a_i$};
  \node at (2.5,0) (tr4) {$\cdots$};
  \node at (3,0) (tr5) {$a_n$};
  \node at (0,-2) (f) {$f^1$};
  \node at (1,-1) (diff) {$d_\Alg$};
  \node at (-1,-3) (bl) {};
  \node at (0,-3) (bc) {};
  \draw[modarrow] (tc) to (f);
  \draw[modarrow] (f) to (bc);
  \draw[blgarrow] (tr1) to (f);
  \draw[blgarrow] (tr3) to (diff);
  \draw[blgarrow] (tr5) to (f);
  \draw[algarrow] (diff) to (f);
  \draw[algarrow] (f) to (bl);
\end{tikzpicture}}
\!\!{+}\!\!
\mathcenter{\begin{tikzpicture}
  \node at (0,0) (tc) {$x$};
  \node at (1,0) (tr1) {$a_1$};
  \node at (1.5,0) (tr2) {$\cdots$};
  \node at (2,0) (tr3) {$a_i$};
  \node at (3,0) (tr3b) {$a_{i+1}$};
  \node at (3.5,0) (tr4) {$\cdots$};
  \node at (4,0) (tr5) {$a_n$};
  \node at (0,-2) (f) {$f^1$};
  \node at (1,-1) (mu) {$\mu_\Alg$};
  \node at (-1,-3) (bl) {};
  \node at (0,-3) (bc) {};
  \draw[modarrow] (tc) to (f);
  \draw[modarrow] (f) to (bc);
  \draw[blgarrow] (tr1) to (f);
  \draw[blgarrow] (tr3) to (mu);
  \draw[blgarrow] (tr3b) to (mu);
  \draw[blgarrow] (tr5) to (f);
  \draw[algarrow] (mu) to (f);
  \draw[algarrow] (f) to (bl);
\end{tikzpicture}}\\
&\quad+\!\!\!
\mathcenter{\begin{tikzpicture}
  \node at (0,0) (tc) {$x$};
  \node at (1,0) (tr1) {$a_1$};
  \node at (1.5,0) (tr2) {$\cdots$};
  \node at (2,0) (tr3) {$a_i$};
  \node at (3,0) (tr3b) {$a_{i+1}$};
  \node at (3.5,0) (tr4) {$\cdots$};
  \node at (4,0) (tr5) {$a_n$};
  \node at (0,-2) (f) {$f^1$};
  \node at (0,-3) (deltaN) {$\delta^1_N$};
  \node at (-1,-4) (mu) {$\mu_\Alg$};
  \node at (-2,-5) (bl) {};
  \node at (0,-5) (bc) {};
  \draw[modarrow] (tc) to (f);
  \draw[modarrow] (f) to (deltaN);
  \draw[modarrow] (deltaN) to (bc);
  \draw[blgarrow] (tr1) to (f);
  \draw[blgarrow] (tr3) to (f);
  \draw[blgarrow] (tr3b) to (deltaN);
  \draw[blgarrow] (tr5) to (deltaN);
  \draw[algarrow] (f) to (mu);
  \draw[algarrow] (deltaN) to (mu);
  \draw[algarrow] (mu) to (bl);
\end{tikzpicture}}
\!\!\!{+}\!\!\!
\mathcenter{\begin{tikzpicture}
  \node at (0,0) (tc) {$x$};
  \node at (1,0) (tr1) {$a_1$};
  \node at (1.5,0) (tr2) {$\cdots$};
  \node at (2,0) (tr3) {$a_i$};
  \node at (3,0) (tr3b) {$a_{i+1}$};
  \node at (3.5,0) (tr4) {$\cdots$};
  \node at (4,0) (tr5) {$a_n$};
  \node at (0,-2) (deltaN) {$\delta^1_M$};
  \node at (0,-3) (f) {$f^1$};
  \node at (-1,-4) (mu) {$\mu_\Alg$};
  \node at (-2,-5) (bl) {};
  \node at (0,-5) (bc) {};
  \draw[modarrow] (tc) to (deltaN);
  \draw[modarrow] (deltaN) to (f);
  \draw[modarrow] (f) to (bc);
  \draw[blgarrow] (tr1) to (deltaN);
  \draw[blgarrow] (tr3) to (deltaN);
  \draw[blgarrow] (tr3b) to (f);
  \draw[blgarrow] (tr5) to (f);
  \draw[algarrow] (f) to (mu);
  \draw[algarrow] (deltaN) to (mu);
  \draw[algarrow] (mu) to (bl);
\end{tikzpicture}}
\end{align*}
where we sum over the index $i$ when it occurs.
This makes the type \DA\ bimodules over $\Alg$ and $\Blg$ into a \dg
category.
(See~\cite{LOT2} for the analogue for type \DA\ structures of
$\Ainf$-algebras, which in some ways looks simpler, though the result
is an $\Ainf$-category.)

\begin{definition}
  An {\em $S$-filtered  $\Alg\Hyph\Blg$  type \DA\ bimodule}
  is a collection of $\Alg\Hyph\Blg$ type \DA\ bimodules $M^s$ indexed
  by $s\in S$, and a collection of morphisms $F^{s<t}\co
  M^s\to M^t[1]$ indexed by $s, t\in S$ with $s<t$, satisfying
  the composition law that if $s<t<u$, then $\sum_t F^{t<u}\circ
  F^{s<t}=d(F^{s<u})$.
\end{definition}

Given an $S$-filtered $\Alg\Hyph\Blg$ type \DA\ bimodule
$\{M^s,F^{s<t}\}$, we define its \emph{total (type \DA) bimodule}
$(M,\delta^1)$ by setting $M=\bigoplus_{s\in S}M^s$ and 
\[
\delta^1(x_s\otimes a_1\otimes\cdots\otimes a_j)=\delta^1_{s}(x_s\otimes a_1\otimes\cdots\otimes a_j)+\sum_{s<t}F^{s<t}(x_s\otimes a_1\otimes\cdots\otimes a_j).
\]

\begin{example}\label{eg:mapping-cone}
  Let $B_1$ and $B_2$ be two $\Alg\Hyph\Blg$ type \DA\ bimodules, and
  let $f\co M^0\to M^1$ be a bimodule morphism. Then the
  \emph{mapping cone of $f$}, denoted $\Cone(f)$, is the $\Alg\Hyph\Blg$ type \DA\ bimodule
  whose underlying vector space is $M^0[1]\oplus M^1$, and whose
  differential is given by the matrix
  $$
  \left(\begin{array}{ll}
      \partial_0 & 0 \\
      f & \partial_1
    \end{array}\right),
  $$ 
  where here $\partial_0$ and $\partial_1$ are the
  differentials on $M$ and $N$ respectively. The mapping cone is the
  total bimodule of the $\{0,1\}$-filtered type \DA\ bimodule
  $\{M^0[1],M^1,f^{0<1}=f\}$.
\end{example}

The following is immediate from the definitions (compare Lemma~\ref{lem:DT-descends}):
\begin{lemma}
  \label{lem:FilteredProduct}
  Suppose that $M$ is an $S$-filtered type \DA\ bimodule over
  $\Alg\Hyph\Blg$, and $N$ is a $T$-filtered type \DA\ bimodule over
  $\Blg\Hyph\Clg$. Assume that either $M$ is right-bounded or $N$ is
  left-bounded, so the tensor product $M\DT N$ is defined. Then $M\DT N$ is naturally an
  $(S\times T)$-filtered $\Alg\Hyph\Clg$ type \DA\ bimodule.
\end{lemma}

As for $S$-filtered $\Ainf$-modules or $S$-filtered type $D$
structures, the category of $S$-filtered type \DA\ bimodules forms a
\dg category in an obvious way. In particular, it makes sense to talk
about homotopy equivalences between $S$-filtered type \DA\
bimodules. We have the following analogue of the rest of Lemma~\ref{lem:DT-descends}:
\begin{lemma}
  \label{lem:DT-descends-bimod}
  Suppose that $M$ and $M'$ are homotopy equivalent $S$-filtered type \DA\ bimodules over
  $\Alg\Hyph\Blg$, and $N$ and $N'$ are homotopy equivalent $T$-filtered type \DA\ bimodule over
  $\Blg\Hyph\Clg$. Suppose, moreover, that either
  \begin{enumerate}
  \item $M$, $M'$ and the homotopy equivalence between them (including
    the homotopies) are bounded or
  \item $N$, $N'$ and the homotopy equivalence between them (including
    the homotopies) are bounded.
  \end{enumerate}
  Then $M\DT N$ and $M'\DT N'$ are homotopy equivalent $(S\times
  T)$-filtered type \DA\ bimodules over $\Alg\Hyph\Clg$.
\end{lemma}
As was the case for Lemma~\ref{lem:DT-descends}, the proof of 
Lemma~\ref{lem:DT-descends-bimod} is an easy exercise in unwinding the
definitions.

\subsection{Other kinds of bimodules}

In~\cite{LOT2}, we described bimodules of various types.
Lemmas~\ref{lem:DT-descends-bimod} has natural generalizations to
other kinds of bimodules. Rather than attempting to catalogue all of these here,
we will describe the one special case which is of interest to us.

Let $\Alg$ and $\Blg$ be \dg algebras.  By a \emph{type \AAm\ bimodule
  over $\Alg$ and $\Blg$} we simply mean an $\Ainf$-bimodule; see
also~\cite[Section 2]{LOT2}. By a \emph{type \DD\ bimodule over $\Alg$
  and $\Blg$} we mean a type $D$ structure over $\Alg\otimes\Blg^\op$,
where $\Blg^\op$ denotes the opposite algebra to $\Blg$. Given a type
\DD\ bimodule $\lsup{\Alg}M^\Blg$ and a type \AAm\ bimodule
$\lsub{\Blg}P_\Clg$, such that either $M$ or $P$ is (appropriately)
bounded, their $\DT$ tensor product
$\lsup{\Alg}M^\Blg\DT\lsub{\Blg}P_\Clg$ is a type \DA\
bimodule. Moreover, given another type \DD\ bimodule
$\lsup{\Alg}N^\Blg$ and a morphism $f\in\Mor(M,N)$, so that either $P$
or $M$ and $N$ are bounded, there is an induced morphism
\[
(f\DT\Id)\co \lsup{\Alg}M^\Blg\DT\lsub{\Blg}P_\Clg\to \lsup{\Alg}N^\Blg\DT\lsub{\Blg}P_\Clg.
\]
\begin{lemma}
  \label{lem:DDMappingCones}
  Suppose that $$f\co P \to Q$$ is a homomorphism of type \DD\
  bimodules over $\Blg$-$\Clg$ and $M$ is a bounded type \AAm\
  bimodule over $\Alg$-$\Blg$. Then there is a canonical isomorphism
  $$\Cone(\Id_M\DT f)\cong M\DT \Cone(f).$$
\end{lemma}
\begin{proof}
  This follows immediately from the definitions.
\end{proof}


\section{The truncated algebra}
\label{sec:Trunc}

Throughout this paper, we will be using the bordered Floer homology
package developed in~\cite{LOT1, LOT2}.  We will say little about the
particular constructions, referring the interested reader to those
sources for details. Almost all of our constructions here will involve
the simpler type $D$ structures, rather than $\Ainf$-modules.

More importantly, we find it convenient to work with the smaller model
of the bordered algebra, introduced in
\cite[Proposition~\ref*{LOT2:prop:SmallerModel}]{LOT2}. This is the
quotient of the usual bordered algebra introduced in~\cite{LOT1},
divided out by the differential ideal generated by strands elements
which have multiplicity at least $2$ somewhere. This algebra is
smaller than the bordered algebra, and hence it is more practical for
calculations. Using this smaller algebra also
makes case analyses, such as the one in
Proposition~\ref{prop:factor-M-minus} below, simpler. 

As shown in
\cite[Proposition~\ref*{LOT2:prop:SmallerModel}]{LOT2}, the two
algebras are quasi-isomorphic, and hence their module categories are
identified. Thus, constructions such as the one in
Section~\ref{sec:SkeinSequence} can be verified for either algebra. 

In~\cite{LOT2}, the algebras are distinguished in their notation:
$\Alg(\PMC)$ is the untruncated version, while $\Alg'(\PMC)$ is the
quotient. Since we will never need the untruncated version
here, we will drop the prime from the notation, and write $\Alg(\PMC)$
for the {\em truncated} algebra.

 \section{Skein sequence}
\label{sec:SkeinSequence}

The aim of the present section is to realize the bimodule of a Dehn
twist as a mapping cone of bimodules, corresponding to the identity
cobordism, and a filling of the surgery curve in the identity
cobordism. 

The results of this section can be interpreted in terms of Type~$D$
modules over the original bordered algebras of~\cite{LOT1}, or
Type~$D$ modules over the truncated algebra discussed in
Section~\ref{sec:Trunc}.  The two points of view are interchangeable
as the two algebras are quasi-isomorphic; see
\cite[Proposition~\ref*{LOT2:prop:SmallerModel}]{LOT2}. While in this
paper we always use $\Alg(\PMC)$ to denote the truncated algebra, the
discussion in the present section works for either the original
or the truncated algebra.

Let $\gamma\subset F(\PMC)$ be a curve. Then,
$\gamma\times \{\OneHalf\}\subset [0,1]\times F(\PMC)$ inherits a
surface framing, and we can form the three-manifold $Y_{0(\gamma)}$,
which is $0$-framed surgery along this curve (with its inherited
surface framing).
 
If we perform $-1$ surgery on $\gamma$, the result is naturally the
mapping cylinder of a positive Dehn twist $\Dehn_\gamma$ along
$\gamma$, while $+1$ surgery is the mapping cylinder of a negative
Dehn twist along $\gamma$ (written $\Dehn_\gamma^{-1})$. In fact, we have the
following:

\begin{theorem}\label{thm:Dehn-is-MC}
  There is a map $F^+_\gamma\co
  \CFDDa(Y_{0(\gamma)})\to\CFDDa(\Id)$ so that
  $\CFDDa(\Dehn_\gamma)$ is homotopy equivalent to the mapping cone
  $\Cone(F^+_\gamma)$ of $F^+_\gamma$.  Similarly, there is a map
  $F^-_\gamma\co \CFDDa(\Id)\to \CFDDa(Y_{0(\gamma)})$ so that
  $\CFDDa(\Dehn_\gamma^{-1})$ is homotopy equivalent to the mapping
  cone $\Cone(F^-_\gamma)$ of~$F^-_\gamma$.

  Analogous statements hold for type \DA\ and \AAm\ bimodules.
\end{theorem}

\begin{proof}
  We prove the statement for $F^+_\gamma$; the statement for
  $F^-_\gamma$ is similar.  The three $3$-manifolds $Y_\Id=[0,1]\times
  F(\PMC)$, $Y_{\Dehn(\gamma)}$ and $Y_{0(\gamma)}$ are all gotten by
  Dehn filling on $Y=([0,1]\times F(\PMC))\setminus
  \nbd(\{\OneHalf\}\times\gamma)$.  Let $T$ be $\bdy
  \nbd(\{\OneHalf\}\times\gamma)$, with orientation agreeing with $\bdy
  Y$.  The three filling slopes on $T$ are $\infty$, $-1$ and
  $0$, respectively. Change the framing on $T$ so that these three
  slopes become $-1$, $0$ and $\infty$, respectively.

  Choose a framed arc in $Y$ connecting $T$ to one of the other two boundary
  components of $Y$. Together with the framing of $T$, this makes $Y$ into a
  strongly bordered $3$-manifold with three boundary components. Let
  $\HD$ be a bordered Heegaard diagram for $Y$.

  Let $\HD_\infty$, $\HD_{-1}$ and $\HD_0$ denote the standard bordered Heegaard
  diagrams for the $\infty$-, $(-1)$- and $0$-framed solid tori,
  respectively:
  \[
  \HD_\infty:\mfigb{torus-30}\qquad
  \HD_{-1}:\mfigb{torus-20} \qquad
  \HD_0:\mfigb{torus-10}.
  \]
  It is straightforward to verify that there is a map $\psi\co
  \CFDa(\HD_\infty)\to \CFDa(\HD_{-1})$ so that $\CFDa(\HD_0)\simeq
  \Cone(\psi)$; see~\cite[Section 11.2]{LOT1}.

  The union $\HD\cup \HD_\infty$ (respectively $\HD\cup \HD_{-1}$,
  $\HD\cup\HD_0$) is a bordered Heegaard diagram for $Y_{0(\gamma)}$
  (respectively $Y_\Id$, $Y_{\Dehn(\gamma)}$). Thus,
  \begin{align*}
    \CFDDa(Y_{\Dehn(\gamma)})&\simeq \CFDDa(\HD\cup\HD_0)\\
    &\simeq \CFDDAa(\HD)\DT\CFDa(\HD_0)\\
    &\simeq \CFDDAa(\HD)\DT\Cone(\psi\co
  \CFDa(\HD_\infty)\to \CFDa(\HD_{-1}))\\
    &\simeq\Cone((\Id\DT\psi)\co (\CFDDAa(\HD)\DT\CFDa(\HD_\infty))\to
    (\CFDDAa(\HD)\DT\CFDa(\HD_{-1})))\\
    &\simeq \Cone(\psi'\co \CFDDa(Y_\Id)\to \CFDDa(Y_{0(\gamma)})),
  \end{align*}
  from the pairing theorem, the fact that $\CFDa(\HD_0)$ is a mapping
  cone, Lemma~\ref{lem:DDMappingCones} and the pairing theorem, again.

  The analogous statements for $\CFDAa$ and $\CFAAa$ follow by
  tensoring with $\CFAAa(\Id_{\PMC})$.
\end{proof}

Suppose that $Y_L$ and $Y_R$ are bordered $3$-manifolds with $\bdy Y_L=F(\PMC)=-\bdy Y_R$. Consider the three closed $3$-manifolds
$Y_L\cup_\bdy Y_R$, $Y_L\cup_\bdy \Dehn_\gamma\cup_\bdy Y_R$ and $Y_L\cup_\bdy Y_{0(\gamma)}\cup_\bdy Y_R$. It follows from Theorem~\ref{thm:Dehn-is-MC} that there are exact triangles
\[
\begin{tikzpicture}
  \node at (0,0) (tc) {$\HFa(Y_L\cup_\bdy Y_R)$};
  \node at (2.3, 2) (br) {$\HFa(Y_L\cup_\bdy \Dehn_\gamma\cup_\bdy Y_R)$};
  \node at (-2.3, 2) (bl) {$\HFa(Y_L\cup_\bdy Y_{0(\gamma)}\cup_\bdy Y_R)$};
  \draw[->] (tc) to (br);
  \draw[->] (br) to (bl);
  \draw[->] (bl) to node[above, left]{\lab{\Id\DT F^+_\gamma\DT \Id}} (tc);
\end{tikzpicture}\hspace{-.55in}
\begin{tikzpicture}
  \node at (0,2) (tc) {$\HFa(Y_L\cup_\bdy Y_R)$};
  \node at (-2.4, 0) (br) {$\HFa(Y_L\cup_\bdy \Dehn_\gamma^{-1}\cup_\bdy Y_R)$};
  \node at (2.4, 0) (bl) {$\HFa(Y_L\cup_\bdy Y_{0(\gamma)}\cup_\bdy Y_R).$};
  \draw[->] (br) to (tc);
  \draw[->] (bl) to (br);
  \draw[->] (tc) to node[below, right]{\lab{\Id\DT F^-_\gamma\DT \Id}} (bl);
\end{tikzpicture}
\]

Previously, by counting holomorphic triangles, Ozsv\'ath-Szab\'o constructed another exact triangle of the same form~\cite{OS04:HolDiskProperties}. We will show in the sequel that these triangles agree; see also~\cite{LOTCobordisms}:
\begin{proposition}\label{prop:triangles-agree}
  The exact triangle on $\HFa$ defined in~\cite[Chapter 11]{LOT1}
  agrees with the exact triangle on $\HFa$ defined
  in~\cite{OS04:HolDiskProperties}. 

  More precisely, consider bordered $3$-manifolds $Y_L$, $Y^i_C$
  ($i=1,2,3$) and $Y_R$, where $Y_C^i$ has two boundary
  components $\bdy_LY_C^i$ and $\bdy_RY_C^i$, and $\bdy
  Y_L=-\bdy_LY_C^i$, $\bdy_RY_C^i=-\bdy Y_R$ (for $i=1,2,3$).  Assume
  that $Y_C^i$ are gotten by $\infty$-, $(-1)$- and $0$-framed Dehn
  filling on the same $3$-manifold with three boundary components.
  Fix bordered Heegaard diagrams $\HD_L$, $\HD_C^i$ and $\HD_R$ for
  $Y_L$, $Y_C^i$ and $Y_R$ so that
  $\HD_C^i=(\Sigma_C,\alphas_C,\betas_C^i,z)$ (i.e., only the
  $\beta$-curves change). For $i\in\ZZ/3$, let 
  \[
  f_{OS}^{i,i+1}\co \CFa(\HD_L\cup\HD_C^i\cup\HD_R)\to \CFa(\HD_L\cup\HD_C^{i+1}\cup\HD_R)
  \]
  denote the map from the exact triangle
  in~\cite{OS04:HolDiskProperties} and let
  \[
  f_{LOT}^{i,i+1}\co \CFDAa(\HD_C^i)\to \CFDAa(\HD_C^{i+1})
  \]
  denote the map from the exact triangle in~\cite{LOT1}, as in the
  proof of Theorem~\ref{thm:Dehn-is-MC}. Then for each $i$, the following diagram
  homotopy commutes:
  \[
  \begin{tikzpicture}
    \node at (0,0) (T1) {$\CFa(\HD_L\cup\HD_C^i\cup\HD_R)$};
    \node at (8,0) (T2) {$\CFa(\HD_L\cup\HD_C^{i+1}\cup\HD_R)$};
    \node at (0,-2) (B1) {$\CFAa(\HD_L)\DT\CFDAa(\HD_C^i)\DT\CFDa(\HD_R)$};
    \node at (8,-2) (B2) {$\CFAa(\HD_L)\DT\CFDAa(\HD_C^{i+1})\DT\CFDa(\HD_R)$.};
    \draw[->] (T1) to node[above]{\lab{f_{OS}^{i,i+1}}} (T2);
    \draw[->] (T1) to (B1);
    \draw[->] (B1) to node[below]{\lab{f_{LOT}^{i,i+1}}} (B2);
    \draw[->] (T2) to (B2);
  \end{tikzpicture}
  \]
  Here, the vertical arrows are the homotopy equivalences given by the
  pairing theorem.
\end{proposition}

A version of Theorem~\ref{thm:BorderedSpectralSequence} now comes for
free.

\begin{proof}[Proof of Theorem~\ref{thm:BorderedSpectralSequence}]
  By the pairing theorem,
  \begin{equation}
    \label{eq:PairingIsomorphism}
    \CFa(\Sigma(L))
  \simeq \CFAa(B_0)\DT \CFDAa(B_1)\DT\dots \DT \CFDAa(B_n)\DT \CFDa(B_{n+1}).
  \end{equation}
  Theorem~\ref{thm:Dehn-is-MC} expresses each $\CFDAa(B_i)$ for
  $i=1,\dots,n$ as a mapping cone of 
  $$f\co
  \CFDAa(B_i^0)\to\CFDAa(B_i^1).$$
  As in Example~\ref{eg:mapping-cone}, we think of this as
  a $\{0,1\}$-filtered type \DA\ bimodule.  Thus, by
  Lemma~\ref{lem:FilteredProduct}, the right-hand-side of
  Equation~\eqref{eq:PairingIsomorphism} is naturally a
  $\{0,1\}^n$-filtered chain complex. At each vertex $v\in\{0,1\}^n$
  in the filtration, we have
  $$
  \CFAa(B_0)\DT\CFDAa(B_1^{v_1})\DT\dots\DT\CFDAa(B_n^{v_n})\DT \CFDa(B_{n+1})
  \simeq \CFa(\Sigma(v(L))),
  $$ 
  where $v(L)$ denotes the resolution
  of $L$ specified by the vector $v$;
  this follows from  another application
  of the pairing theorem. Note that $v(L)$ is an unlink, so
  $\Sigma(v(L))$ is a connected sum of copies of
  $S^2\times S^1$; and hence the $E_1$ term in the spectral sequence
  associated to the filtration
  gives $\HFa(\Sigma(v(L)))$, which in turn is identified with the reduced Khovanov
  homology of the unlink $v(L)$.

  The differentials on the $E_1$ page all have the form $\Id\DT F^+_\gamma \DT\Id$ or $\Id\DT F^-_\gamma\DT\Id$. By Proposition~\ref{prop:triangles-agree}, these differentials agree with the $E_1$ differentials in~\cite{BrDCov}. It is shown in~\cite[Theorem 6.3]{BrDCov} that these agree, in turn, with the differentials on the Khovanov chain complex.
\end{proof}

In Section~\ref{sec:edges} we will need the following property of the
maps $F^\pm_\gamma$:
\begin{definition}\label{def:triangle-like}
  Let $\HD=(\Sigma,\alphas,\betas,\gammas,z)$ be a bordered Heegaard
  triple diagram with $(\Sigma,\betas,\gammas)$ the standard Heegaard
  diagram for $\#^k(S^2\times S^1)$. Let $\theta_{\beta,\gamma}$ denote
  the top-dimensional generator of $\CFa(\Sigma,\betas,\gammas)$ and
  let $F\co \CFDDa(\Sigma,\alphas,\betas)\to
  \CFDDa(\Sigma,\alphas,\gammas)$ be a homomorphism. For each generator
  $\x\in\Gen(\Sigma,\alphas,\betas)$ write
  \[
  F(\x)=\sum_{\y\in\Gen(\Sigma,\alphas,\gammas)}f_{\x,\y}\otimes \y,
  \]
  with $f_{\x,\y}\in\Alg$.
  We say that $F$ is \emph{triangle-like} if for each $\x$, $\y$ and
  basic generator $a$ in $f_{\x,\y}$ there is a domain $B\in
  \pi_2(\x,\theta_{\beta,\gamma},\y)$ and a compatible sequence of Reeb chords
  $\vec{\rho}=(\rho_1,\dots,\rho_k)$ so that $\ind(B,\vec{\rho})=0$
  and $a=a(-\rho_1)\cdots a(-\rho_k)$. (Here, $\ind(B,\vec{\rho})$
  denotes the expected dimension of the moduli space of embedded
  triangles in the homology class~$B$ with east asymptotics
  $\vec{\rho}$.)
\end{definition}

\begin{lemma}\label{lem:triangle-like}
  The maps $F^\pm_\gamma$ are triangle-like.
\end{lemma}
\begin{proof}
  This follows from the observations that:
  \begin{itemize}
  \item The maps in~\cite{LOT1} giving the skein sequence are
    triangle-like, and
  \item If $F\co \CFDa(Y_1)\to \CFDa(Y_2)$ is triangle-like then
    $\Id\DT F\co \CFDAa(Y)\DT\CFDa(Y_1)\to \CFDAa(Y)\DT\CFDa(Y_2)$ is
    triangle-like.
  \end{itemize}
  The first point is immediate from the definitions of the maps, and
  the second follows from gluing properties of domains and of the index.
\end{proof}

\begin{remark}
  The number $\ind(B,\vec{\rho})$ is computed by a formula similar
  to~\cite[Definition~\ref*{LOT1:def:emb-ind-emb-chi}]{LOT1},
  cf.~\cite{Sarkar11:IndexTriangles}, so it is not hard to check if a
  given map is triangle-like.
\end{remark}


\section{The plat handlebody}\label{sec:diagrams}
\subsection{The linear pointed matched circle}
Recall that a pointed matched circle $\PMC$ consists of an oriented
circle $Z$, $4k$ points $\mathbf{a}\subset Z$, a matching $M$ of the
points $\mathbf{a}$ in pairs, and a basepoint $z\in Z\setminus
\mathbf{a}$, subject to the compatibility condition that 
the manifold obtained by performing
surgery on $Z$ along the $M$-paired  points in $\mathbf{a}$
is connected. The pointed matched circle $\PMC$
specifies both a surface $F(\PMC)$ and a \dg algebra $\Alg(\PMC)$.
The algebra $\Alg(\PMC)$ is
described in~\cite[Chapter 2]{LOT1}. The
surface $F(\PMC)$ is gotten as follows:
\begin{enumerate}
\item Fill $Z$ with a disk.
\item Attach a $2$-dimensional $1$-handle along each pair of points in
  $Z$. The boundary of the resulting surface is obtained by performing
  surgery along the pairs of points in $Z$ and hence is connected by hypothesis.
\item Glue a disk to this boundary.
\end{enumerate}

The pointed matched circle $\PMC$ also endows the surface $F(\PMC)$ with a preferred set of
$2k$ embedded, closed curves $\gamma_1,\dots,\gamma_{2k}$. Specifically, for each $M$-paired point in $\mathbf{a}$, we let $\gamma_i$ consist
of an arc in $Z$ which connects the two points (so as not to cross the basepoint $z$),
which is then joined to an arc which runs through the one-handle associated to the matched pair.
One can orient $\gamma_i$ so that its portion in $Z$ agrees with the orientation of $Z$.
The corresponding homology classes of the $\{\gamma_n\}_{n=1}^{2k}$  form a  basis for
$H_1(F(\PMC))$.

\begin{figure}
  \begin{center}
    \input{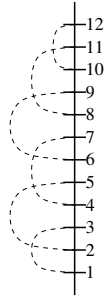}
  \end{center}
  \caption {{\bf Linear pointed matched circles.}
    \label{fig:ArcslideMatching}
    We have here an illustration of the pointed matched
    circle in genus~$3$. We have cut the circle at the basepoint $z$ to obtain
    the interval shown.}
\end{figure}

In our present considerations, we would like a pointed
matched circle for which certain Dehn twists have convenient
representations. For each $k$, define the \emph{genus $k$ linear
  pointed matched circle} as follows. In the unit circle
$Z=\{z\in\CC\mid |z|=1\}$ let $z=\{1\}$. Let $\mathbf{a}=\{e^{2\pi i
  n/(4k+1)}\mid 1\leq n\leq 4k\}$; for ease of notation, let
$a_n=e^{2\pi i n/(4k+1)}$. The matching $M$ matches up the
following pairs of points: $\{a_1,a_3\}$, $\{a_{4k-2},a_{4k}\}$, and,
for $n=1,\dots,2k-2$, $\{a_{2n}, a_{2n+3}\}$. See
Figure~\ref{fig:ArcslideMatching} for an example.

\subsection{A particular handlebody}
Consider the linear pointed matched circle $\PMC$. Let $\gamma_1$ be the
curve in $F(\PMC)$ corresponding to $\{a_1,a_3\}$, $\gamma_{n+1}$ the
curve corresponding to $\{a_{2n},a_{2n+3}\}$ for $n=1,\dots,2k-2$, and
$\gamma_{2k}$ the curve in $F(\PMC)$ corresponding to
$\{a_{4k-2},a_{4k}\}$. Note, in particular, that (a small
perturbation of) $\gamma_i$ intersects
each of $\gamma_{i-1}$ and $\gamma_{i+1}$ transversely once, and is
disjoint from the other $\gamma_j$.

\begin{figure}
\begin{center}
  \includegraphics{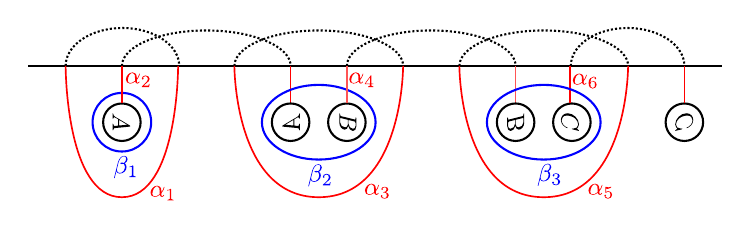}
\end{center}
\caption {{\bf Heegaard diagram for the plat handlebody.}
\label{fig:HeegaardHandlebody}
The genus $3$ case is shown.  The diagram has been rotated $90^\circ$
clockwise to fit better on the page.}
\end{figure}

We consider the handlebody in which the odd curves
$\{\gamma_{2n-1}\}_{n=1}^k$ bound disks. We call this the \emph{genus $k$ plat
  handlebody}. A Heegaard diagram for the plat handlebody is
illustrated in Figure~\ref{fig:HeegaardHandlebody}, and is gotten as
follows.  Consider a disk with boundary $\PMC$. Draw an $\alpha$-arc
$\alpha_{2n-1}$ for each $\gamma_{2n-1}$ (i.e., connecting the same
endpoints as $\gamma_{2n-1}$): thus, we have drawn $\alpha$ arcs which
connect the matched pairs of points $a_1$ and $a_3$; and also $a_{4n}$
and $a_{4n+3}$ for $n=1,\dots,k-1$.  We wish to draw $\alpha$-arcs
(disjoint from the other $\alpha$-arcs) which connect the remaining
matched pairs of points $a_{4n-2}$ and $a_{4n+1}$ with
$n=1,\dots,k-1$ and $a_{4k-2}$ to $a_{4k}$.  This, of course, cannot
be done in the disk, since $a_{2n}$ and $a_{2n+3}$ are separated by
both arcs $\alpha_{2n-1}$ and $\alpha_{2n+3}$ (also, $a_{4k-2}$ and
$a_{4k}$ are separated by $\alpha_{2k-1}$).  Instead,  add a
one-handle $A_n$ for $n=1,\dots, k$
to the disk near each of the remaining 
such pair of points $a_{4n-2}$ and $\alpha_{4n+1}$ for 
$n=1,\dots,k-1$, and when $n=k$ the pair $a_{4k-2}$ and $a_{4k}$; 
then draw arcs $\alpha_{2n}$ for $n=1,\dots, k$
running through $A_n$ connecting these pairs of points.
Finally, draw a $\beta$-circle $\beta_n$ for $n=1,\dots,k$
encircling $A_n$ and $A_{n-1}$ (so $\beta_1$ encircles only $A_1$).
This constructs the Heegaard diagram for the plat handlebody.

Write $\rho_{i,j}$ for the Reeb chord connecting $a_i$ to $a_{j}$.

There is a domain ${\mathcal D}_n$ for $n=1,\dots,k$ whose boundary
consists of $\alpha_{2n-1}$ and $\beta_n$.
We have $\bdy^\bdy{\mathcal D}_1=\rho_{1,3}$, 
$\bdy^\bdy{\mathcal D}_{n}=\rho_{4n-4,4n-1}$ for $n=1,\dots, k$.

Let $\SetS$ be the subset of $\{1,\dots,4k\}$ which are joined by the
odd $\alpha$-arcs $\alpha_{2n-1}$ for $n=1,\dots k$.  Write
\[
\xi_n = \left\{\begin{array}{ll}
a(\rho_{1,3})\cdot I(\SetS)
 & {\text{if }n=1} \\
a(\rho_{4n-4,4n-1})\cdot I(\SetS)
& {\text{if $n=2,\dots,k$.}}
\end{array}\right.
\]

\begin{proposition}\label{prop:comp-plat-handlebody}
  Let $\HD$ denote the Heegaard diagram described above for the plat handlebody. Then $\CFDa({\HD})$
  has a single generator $\x$, and
  $$\partial \x = \left(\sum_{n=1}^k \xi_n\right)\x.$$
\end{proposition}

\begin{proof}
  For $n=1,\dots,k$, the circle $\beta_n$ meets $\alpha_{2n}$ in a
  single intersection point, which we denote $x_n$. The tuple
  $\x=\{x_1,\dots,x_n\}$ is the unique generator. (For $n>1$,
  $\beta_n$ also meets $\alpha_{2n-2}$, but none of these points can
  be completed to a generator.)  By definition, $I(\SetS)\cdot \x =
  \x$.
  
  Each domain $\mathcal{D}_n\in\pi_2(\x,\x)$ has a unique holomorphic
  representative $u\co S\to \Sigma\times[0,1]\times\RR$, as
  in~\cite[Chapter~\ref*{LOT1:chap:structure-moduli}]{LOT1}. Indeed,
  each such domain can be cut (from $\x$ to $\bdy\Sigma$) to give a
  bigon from $\x$ to itself, whose contribution is the algebra element
  $\xi_n$, coming from $\partial{\mathcal D}_n$. This is the only
  representative: if the cut did not go out to the boundary, the
  resulting algebra element would not be consistent with the
  idempotents (and there also would not be a holomorphic
  representative, by a maximum modulus principal argument).

  Any element of $\pi_2(\x,\x)$ can be
  expressed as a linear combination of the ${\mathcal D}_n$. It follows that
  $\partial \x = a\x$ for some algebra element $a$ with
  $I(\SetS)\cdot a \cdot I(\SetS)=a$, and whose support is a
  combination of the intervals $[1,3]$ and $[4n,4n+3]$ for
  $n=1,\dots,k-1$.  Looking at the strands algebra, this
  forces $a$ to lie in the subalgebra generated by
  $\{\xi_1,\dots,\xi_k\}$.  (Because of the form of $\mathbf{s}$,
  $d(\xi_i)=0$ for each $i$.)
  We have already seen that each algebra element $\xi_n$
  appears with non-zero coefficients in $a$. Grading reasons prevent any
  other term from appearing in the algebra element $a$:
  since $\xi_n\x$ appears in the differential of $\x$,
  it follows that $\lambda^{-1}\gr(\x)=\gr(\xi_n)\cdot \gr(\x)$, so if $b=\xi_{i_1}\cdot\cdots\cdot\xi_{i_m}$ for some $m>1$ then 
  $\gr(b)\cdot \gr(\x)=\lambda^{-m}\cdot \gr(\x)$,
  and hence $b\otimes \x$ cannot appear in $\partial \x$.
  (We are using here the fact that multiplication by $\lambda$ gives a free $\ZZ$-action on that the grading set of the handlebody.)
\end{proof}


\newcommand\blank{-}
\section{Vertices in the hypercube}
\label{sec:vertices}

According to Section~\ref{sec:SkeinSequence}, the bimodules associated
to the branched
double-cover of the $B_i$ ($i=1,\dots,n$) are realized as mapping 
cones. The aim of the present section is to calculate the \DD-bimodules
appearing in these mapping cones.

Although most of the results we state here hold for rather general
pointed matched circles $\PMC$, we will be primarily interested in the
case of the linear pointed matched circle. Restricting attention entirely
to this case seems a little {\em ad hoc}, so we do not do it. We do, however,
make one simplifying assumption throughout, to make the proofs a little simpler.
To state it, we introduce some terminology (see~\cite{LOT4}):

 \begin{definition}\label{def:special-length-3}
  A {\em special length $3$ chord} in a pointed matched circle $\PMC$
  is a chord which connects four consecutive positions
  $[a_{i},a_{i+1},a_{i+2}, a_{i+3}]$ with the property that $a_i$ and
  $a_{i+2}$ are matched, and $a_{i+1}$ and $a_{i+3}$ are matched.
  A pointed matched circle is called {\em unexceptional}
  if it contains no special length $3$ chords.
\end{definition}

The linear pointed matched circle is unexceptional, except in the
degenerate case where it has genus one.
At various points we will assume that $\PMC$ is unexceptional. Since
we can always stabilize a braid to have at least 6 strands (3
bridges), this does not restrict our results.

\begin{remark}\label{rem:stable}
  This hypothesis can be removed by restricting to ``stable'' modules,
  in the sense of~\cite{LOT4}. Since we will not need this, we leave the details
  as an exercise for the interested reader.
\end{remark}

\subsection{Braid-like resolutions}\label{sec:braid-like}
For the braid-like resolutions, the terms correspond to the three-ball
branched at a collection of arcs. As such, its type \DD\ bimodule is
the type \DD\ bimodule of the identity cobordism, as computed
in~\cite[Theorem~\ref*{HFa:thm:DDforIdentity}]{LOT4}, as we now recall.

Fix a pointed matched circle $\PMC$, and let $\PMC'=-\PMC$ denote the
orientation reverse of $\PMC$.
Let $\SetS$ and $\SetT$ be two subsets of the set of matched pairs
for $\PMC$, so that $\SetS\cap\SetT=\emptyset$ and $\SetS\cup\SetT$
consists of all matched pairs. Let $I(\SetS)$ denote the corresponding
idempotent in $\Alg(\PMC)$ and $I'(\SetT)$ denote the corresponding
idempotent in $\Alg(\PMC')$ (via the natural identification of matched
pairs in $\PMC$ with those in $\PMC'=-\PMC$). We call the element
$I(\SetS)\otimes I'(\SetT)\in \Alg(\PMC)\otimes\Alg(\PMC')$ a pair of {\em complementary idempotents}.

\begin{definition}
  The {\em diagonal subalgebra} $\Dlg$ of $\Alg(\PMC)\otimes \Alg(\PMC')$ is the
  algebra generated by algebra elements $a\otimes a'$ where
  \begin{itemize}
  \item $i \cdot a \cdot j =a$ and $i'\cdot a'\cdot j'=a'$,
    where $i\otimes i'$ and $j\otimes j'$ are pairs of complementary
    idempotents, and
  \item the support of $a$ agrees with (minus) the support of $a'$
    under the obvious orientation-reversing identification of $\PMC$
    and $\PMC'$.
  \end{itemize}
\end{definition}

As discussed in Section~\ref{sec:Trunc}, when we write $\Alg(\PMC)$
above, we are referring to the truncated algebra. An analogous definition can be made 
in the untruncated case, but we have no need for that here.

We call an element $x$ of $\Dlg$ a {\em near-chord} if there is a pair
of elementary idempotents $i$ and $i'$ and a chord $\xi$ (not crossing the basepoint $z$, of course) with the property
that $x=(i\cdot a(\xi))\otimes (i'\cdot a'(\xi))$ where here
$a(\xi)\in\Alg(\PMC)$ is the algebra element associated to $\xi$,
and $a'(\xi)\in\Alg(\PMC')$ is the algebra element associated to the same
chord $\xi$, only now thought of as a chord in $\PMC'$. 
(In~\cite{LOT4}, these elements were called {\em chord-like}; the term near-chord was
reserved for other, similar contexts. We find it preferable to use more uniform terminology
here, however.) Let $A_\Id\in\Dlg$
be the sum of all near-chords. It is easy to see 
\cite[Theorem~\ref*{HFa:thm:DDforIdentity}]{LOT4} that 
$d A_\Id + A_\Id\cdot A_\Id=0$. It follows that if 
$\Alg=\Alg(\PMC)\otimes\Alg(\PMC')$, then the map
$$\nabla_\Id\co \Dlg \to \Dlg$$
defined by
$$\nabla_\Id (B)=dB + B \cdot A_\Id$$ 
is a differential.  This induces a differential on the
ideal $\Alg\cdot \Dlg\subset \Alg=\Alg(\PMC)\otimes\Alg(\PMC')$.
(Explicitly, if $B\in\Alg(\PMC)\otimes\Alg(\PMC')$ and $I$ is an idempotent in $\Dlg$,
$\nabla_\Id(B\cdot I) = (dB)\cdot I + B\cdot A_\Id\cdot I$.)

Indeed, we have the following:

\begin{proposition}\label{prop:RecallIdentity}
  (\cite[Theorem~\ref*{HFa:thm:DDforIdentity}]{LOT4})
  There is a Heegaard diagram for the identity cobordism whose
  associated type \DD\ bimodule $\CFDDa(\Id)$ is identified with the
  ideal $\Alg\cdot \Dlg$, endowed with the differential $d+A_\Id$.
\end{proposition}

In view of the above proposition, we call $A_\Id$ the {\em structure
constant for the \DD\ identity bimodule}.

\begin{example}\label{exam:DD-id-torus}
  Let $\PMC$ be the (unique) pointed matched circle for a genus $1$
  surface. The algebra $\Alg=\Alg(\PMC)$ is $8$-dimensional, with elements
  $\iota_0,\iota_1,\rho_1,\rho_2,\rho_3,\rho_{12},\rho_{23},\rho_{123}$,
  and non-zero products:
  \begin{align*}
    \iota_0^2&=\iota_0 & \iota_1^2&=\iota_1&
    \iota_0\rho_1\iota_1&=\rho_1& \iota_1\rho_2\iota_0&=\rho_2\\
    \iota_0\rho_3\iota_1&=\rho_3&
    \rho_1\rho_2 &= \rho_{12} & \rho_2\rho_3& = \rho_{23} & 
    \rho_1\rho_{23} &= \rho_{123} \\
    \rho_{12}\rho_{3}& = \rho_{123}.
  \end{align*}
  If we label the points in $\CircPts$ as $1,2,3,4$ according to the
  orientation on the interval $Z\setminus\{z\}$ then $\rho_1$
  corresponds to the interval $[1,2]$, $\rho_2$ corresponds to the
  interval $[2,3]$, and $\rho_3$ corresponds to the interval
  $[3,4]$. The idempotent $\iota_0$ corresponds to the matched pair
  $\{1,3\}$ and $\iota_1$ to $\{2,4\}$.  See~\cite[Section 11.1]{LOT1}
  for more details.)

  To avoid confusion, we will keep different copies of $\Alg$ separate
  by sometimes replacing the letters $\rho$ by other letters (e.g.,
  $\sigma$, $\tau$). We will write $\Alg^\sigma$ to denote $\Alg$
  where we label the elements by $\sigma$'s, and so on. Note also that
  $\Alg(\PMC)\cong \Alg(\PMC')$.

  There are two pairs of complementary idempotents in
  $\Alg(\PMC)\otimes\Alg(\PMC')$: $J=\iota_0\otimes\iota_0$ and
  $K=\iota_1\otimes\iota_1$. (The orientation-reversal in the
  isomorphism between $\Alg(\PMC)$ and $\Alg(\PMC')$ is why these
  idempotents do not look complementary.)

  The diagonal subalgebra of $\Alg^\sigma(\PMC)\otimes\Alg^\rho(\PMC')$ is spanned
  by the idempotents and $\sigma_1\otimes\rho_3$,
  $\sigma_2\otimes\rho_2$, $\sigma_3\otimes\rho_1$, and
  $\sigma_{123}\otimes\rho_{123}$. (All products of these
  non-idempotent elements vanish. Note, for instance, that
  $\sigma_{12}\otimes\rho_{12}$ does not appear because it does not
  satisfy the restriction on supports, while
  $\sigma_{12}\otimes\rho_{23}$ does not appear because it is not
  compatible with a pair of complementary idempotents.)

  All of the non-idempotent generators of the diagonal subalgebra are
  near-chords. Thus,
  \[
  A_\Id=\sigma_1\otimes\rho_3+ \sigma_2\otimes\rho_2+
  \sigma_3\otimes\rho_1+ \sigma_{123}\otimes\rho_{123}.
  \]
  So, by Proposition~\ref{prop:RecallIdentity}, the differentials of
  the generators of $\CFDDa(\Id_{T^2})$ are given by
  \begin{align*}
    \bdy(J) &=
    (\sigma_1\otimes\rho_3+\sigma_3\otimes\rho_1+\sigma_{123}\otimes\rho_{123})K\\
    \bdy(K)&= (\sigma_2\otimes\rho_2)J.
  \end{align*}
  Note that this agrees with the computation, via holomorphic curves,
  in~\cite[Section~\ref*{LOT2:subsec:AAId1}]{LOT2}.
\end{example}

In the process of proving Proposition~\ref{prop:RecallIdentity}, the
following is proved in~\cite{LOT4}:

\begin{lemma}\label{lem:factor-diag}
  Every basic generator of the diagonal subalgebra can be factored as a product of near-chords.
\end{lemma}

\begin{proof}
  This is proved
  in~\cite[Lemma~\ref*{HFa:lem:FactorDiagonalSubalgebra}]{LOT4} in the
  untruncated case. The truncated version can be thought of as a formal consequence.
  However, the proof in the truncated case is somewhat simpler, so we
  give a brief sketch.  

  Consider the topmost boundary
  point in the support of our algebra element on the $\PMC$ side.
  There is a strand $s$ which terminates in $p$, and a corresponding
  strand $t$ which starts at the corresponding point $p'$ on the
  $\PMC$ side. If $s$ is shorter than $t$, then we can factor off a
  near-chord on the left, gotten by assembling $s$ and a portion of
  $t$. If $s$ is longer than $t$, we can factor off a near-chord on the right,
  gotten by assembling a portion of $s$ and the strand $t$. If $s$ and $t$ have
  the same length, they can be combined to form a near-chord, which is factored off
  on either the left or the right. By induction on the total support, we obtain the
  desired factorization.
\end{proof}

\begin{lemma}\label{lem:diagonal-Z-graded}
  The (group-valued) gradings on $\Alg(\PMC)$ and $\Alg(\PMC')$ induce
  a $\ZZ$-grading on the diagonal subalgebra $\Dlg$.
\end{lemma}
\begin{proof}
  This follows the
  fact that, in the language of \cite[Definition
  \ref*{HFa:def:coeff-alg}]{LOT4}, $\Dlg$ is the coefficient algebra
  of the type \DD\
  bimodule associated to the standard Heegaard diagram for the
  identity map, by \cite[Lemma~\ref*{HFa:lem:grading-coeff-implies}]{LOT4}.
\end{proof}

\subsection{Anti-braid-like resolutions: Generic case}\label{sec:anti-braid-like}
The terms corresponding to the anti-braid-like resolutions can be
interpreted as $0$-framed surgery in $[0,1]\times F(\PMC)$ along a
curve supported in $F(\PMC)$, giving the manifold $Y_{0(\gamma)}$ considered in
Section~\ref{sec:SkeinSequence}. The curves we consider have a
particularly explicit form.  Fix some matched pair $\{c_1,c_2\}$. This
determines a curve in $F(\PMC)$, which we denote by $\gamma$. 

\begin{remark}
  The three-manifold $Y_{0(\gamma)}$ can be thought of as a composite of two elementary
  cobordisms: a cobordism from a surface of genus $k$ to one of genus $k-1$ 
  (attaching a $2$-handle to $[0,1]\times F$ along $\gamma\subset \{1\}\times F$),
  followed by one which is a cobordism back to $k$ again (gotten by turning around the first).
  As such, it can be described as composition of the two corresponding bimodules.
  This is not the approach we take presently; but the two constituent bimodules
  can be calculated by the same means.  This would also give
  computations of the
  invariants of cups and caps.
\end{remark}

For this subsection, we suppose we have the following condition.

\begin{hypothesis}
  \label{hyp:DistanceTwo}
  We suppose that between the ends of our matched pair $\{c_1,c_2\}$
  in $\PMC$ there are
  exactly two other distinguished points, which we label $u$ and $d$,
  so that $c_1<d<u<c_2$ with respect to the ordering given
  by then orientation of~$\PMC$.
\end{hypothesis}

\begin{figure}
\begin{center}
\input{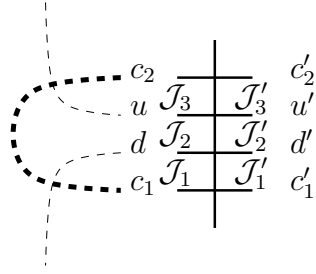}
\end{center}
\caption {{\bf Labelled diagram.}
\label{fig:LabelledDiagram}
A portion of a pointed matched circle between $c_1$ and $c_2$, including
labeling conventions.}
\end{figure}

We label the points in $\PMC'$ corresponding to $u$ and $d$
by $u'$ and $d'$. Note that in $\PMC'$, we have
$u'<d'$ with respect to the orientation on $Z'$. See
Figure~\ref{fig:LabelledDiagram} for an illustration.

\begin{remark}
  \label{def:DistanceTwo}
  The distance of exactly two between the distinguished points might
  seem arbitrary. However, this hypothesis is met this
  in the typical case for the linear pointed matched circles we
  consider in the present paper; the one other case for linear pointed
  matched circles is analyzed in Section~\ref{subsec:DegenerateVertex}.
\end{remark}

\begin{definition}
  \label{def:AntiBraidSubalgebra}
  The distinguished points $u$ and $d$ (respectively $u'$ and $d'$) divide
  the interval between $c_1$ and $c_2$ (respectively $c_1'$ and $c_2'$) into
  three consecutive intervals, labelled $\Zint_1$, $\Zint_2$, and $\Zint_3$ (respectively
  $\Zpint_1$, $\Zpint_2$, and $\Zpint_3$). 
  The {\em anti-braid subalgebra} 
  $\Blg$ is the subalgebra of $\Alg(\PMC)\otimes \Alg(\PMC')$ specified as
  follows.
  \begin{itemize}
  \item Basic idempotents correspond to pairs $\SetS$ and $\SetT$ of subsets
    of matched pairs in $\PMC$ (where $\SetT$ is thought of as a set of
    pairs for $\PMC'$) with the property that:
    \begin{itemize}
    \item The union $\SetS\cup\SetT$ contains all matched pairs
      {\em except} for one of the two containing one of $u$ or $d$.
    \item The intersection $\SetS\cap\SetT$ consists of exactly the
      matched pair $\{c_1,c_2\}$.
    \end{itemize}
    Thus, there are four
    types of basic idempotents in $I(\SetS)\otimes I(\SetT)$, 
    divided according to whether $u\in\SetS$, $d\in\SetS$, $u'\in\SetT$, or
    $d'\in\SetT$. We denote these types $\lsub{u}Y$,  $\lsub{d}Y$,
    $Y_{u}$, and $Y_{d}$ respectively.
  \item The non-idempotent elements are linear combinations of elements of the
    form $a\otimes a'$, constrained as follows:
    \begin{itemize}
    \item There are basic idempotents $i\otimes i'$ and $j\otimes j'$
      in $\Blg$ as above
      with $iaj=a$ and $i'a'j'=a'$.
    \item The local multiplicities of $a$ in $\PMC$
      at any region outside the interval
      $[c_1,c_2]$ coincide with the local multiplicities of $a'$
      at the corresponding region in~$\PMC'$.
    \item The local multiplicity of $a$ at $\Zint_1$ equals the local
      multiplicity of $a$ at $\Zint_3$, and the local multiplicity
      of $a'$ at $\Zpint_1$ equals the local multiplicity of $a'$ at~$\Zpint_3$.
    \end{itemize}
  \end{itemize}
\end{definition}

\begin{figure}
\begin{center}
\input{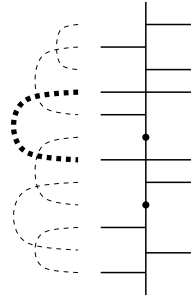}
\end{center}
\caption {{\bf Graphical representation of a typical idempotent in the anti-braid subalgebra.}
\label{fig:Idempotents}
This representation is obtained as follows. Cut
the pointed matched circle along the basepoint $z$ to obtain an
interval. At each distinguished position, draw an edge from
the left of the interval to the center
(if the position lies in the idempotent on the $\PMC$ side)
or from the right of the interval to the center
(if the position lies in the idempotent on the $\PMC'$
side), with the following exceptions: at the distinguished matched pair
$\{c_1,c_2\}$, draw edges from both the left and the right to the center, and at one
of the two positions between $c_1$ and $c_2$, and its matched pair, draw no
edges. If two positions are matched, the corresponding edges must lie on the same side of the center.}
\end{figure}

\begin{figure}
\begin{center}
\input{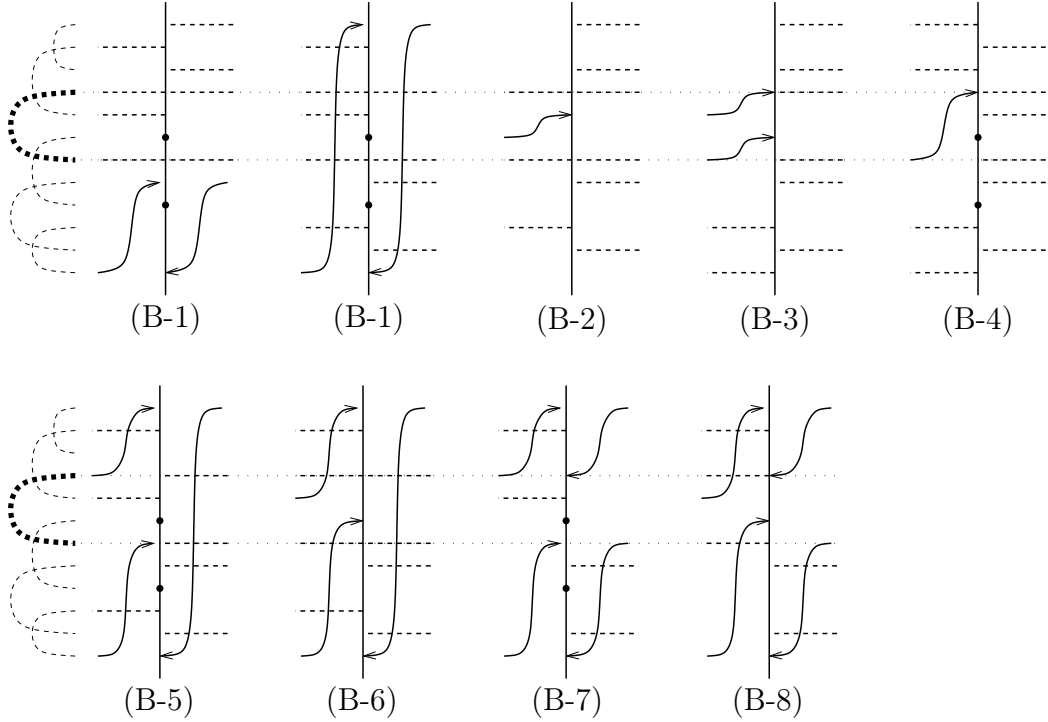}
\end{center}
\caption {{\bf Anti-braid-like near-chords.}
\label{fig:ABLNearChords}
Examples of near-chords in the sense of
Definition~\ref{def:NearChords} of all $8$ types. As the arrows
indicate, we are drawing the left action by $\Alg(\PMC')$ as a right
action by the opposite algebra $\Alg(\PMC')^\op=\Alg(-\PMC')$; this
convention persists in future figures.}
\end{figure}

See Figure~\ref{fig:Idempotents} for a picture of a typical idempotent
(for the genus $3$ linear pointed matched circle), and
Figure~\ref{fig:ABLNearChords} for some examples of elements of the
anti-braid subalgebra.

\begin{figure}
\begin{center}
\input{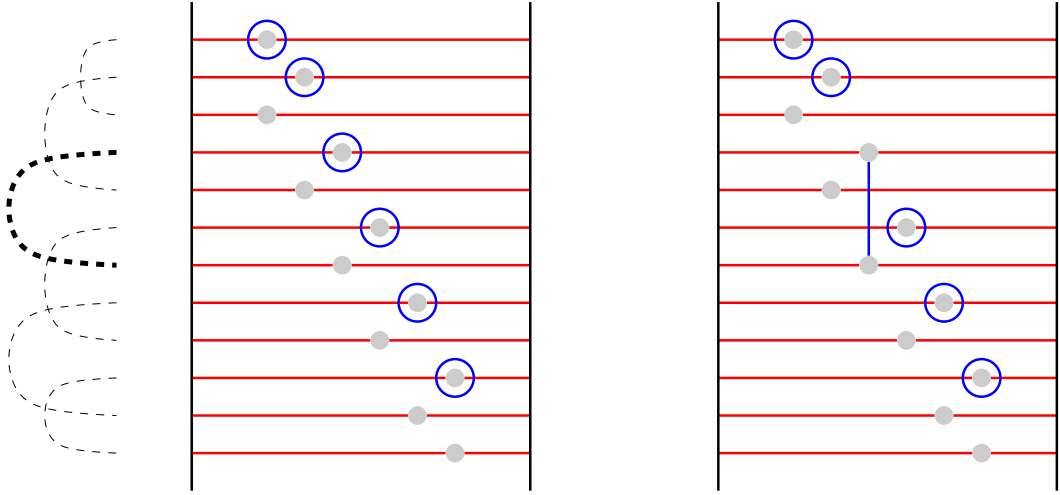}
\end{center}
\caption {{\bf Heegaard diagrams.}
\label{fig:HeegaardDiagrams}
On the left is the standard Heegaard diagram for the
identity cobordism of the genus $3$ linear pointed matched circle $\PMC$.  On the
right is a Heegaard diagram for $Y_{0(\gamma_j)}$, where
$\gamma_j\subset F(\PMC)$ is the circle corresponding to the dark
dashed line. In each diagram, gray circles correspond to handles,
which are identified in pairs according to their horizontal positions. Each
has curves running through them, so that horizontal arcs become
arcs which run through the handle, entering and exiting at the same boundary
components. The vertical arc in the right picture closes up to a closed loop.}
\end{figure}

We give a Heegaard diagram for zero-surgery along $\gamma$, as
illustrated in Figure~\ref{fig:HeegaardDiagrams} on the right.  Call
this Heegaard
diagram $\HD_0(\gamma)$.

\begin{lemma}\label{lem:anti-braid-Z-graded}
  The (group-valued) gradings on $\Alg(\PMC)$ and $\Alg(\PMC')$ induce
  a $\ZZ$-grading on the anti-braid subalgebra $\Blg$.
\end{lemma}
\begin{proof}
  The anti-braid subalgebra $\Blg$ is the coefficient algebra, in the
  sense of~\cite[Definition~\ref*{HFa:def:coeff-alg}]{LOT4}, of the bimodule
  $\CFDDa(\HD_0(\gamma))$. So this follows
  from~\cite[Lemma~\ref*{HFa:lem:grading-coeff-implies}]{LOT4}.
\end{proof}

For $\star, \bullet\in\{ \lsub{u}Y, \lsub{d}Y, Y_u, Y_d\}$, we will
say an element $x$ of $\Blg$ has \emph{type $\star \to \bullet$} (for
instance, type $Y_u\to Y_d$) if $x=I\cdot x\cdot J$ for some idempotents
$I$ of type $\star$ and $J$ of type $\bullet$.
We will say $x$ has
\emph{type $\star\to \star$} if there is some $\bullet\in\{ \lsub{u}Y,
\lsub{d}Y, Y_u, Y_d\}$ so
 $x$ has type $\bullet\to\bullet$.

An interval $\xi$ connecting two points in $\PMC$ is called
{\em generic} if none of its endpoints is $c_1$, $c_2$, $u$, or $d$.

It is convenient to extend the notation $a(\xi)$ slightly. If $S$ is a subset of $Z$ with $\bdy S\subset \mathbf{a}=Z\cap\alphas$ then we can associate a strands diagram to $S$ by turning each connected component of $S$ into a strand. We write $a(S)$ to denote the algebra element associated to this strands diagram. (This notation is used in case (\ref{typeS:I1I3}) of the following definition, for instance.)
\begin{definition}
  \label{def:NearChords}
  A non-zero element $x$ of $\Blg$ is a \emph{near-chord} if $x$
  satisfies one of the following conditions:
  \begin{enumerate}[label=(B-\arabic*),ref=B-\arabic*]
    \item
      \label{typeS:Generic}
      $x=I\cdot (a(\xi)\otimes a'(\xi))\cdot J $, where here $\xi$ is
      some generic interval (i.e., neither of its endpoints is $c_1$,
      $c_2$, $u$, or $d$). In this case, $x$ can be of any of the
      following possible types: $\star\to \star$, $\lsub{u}Y\to Y_u$,
      $Y_u\to\lsub{u}Y$, $\lsub{d}Y\to Y_d$ or $Y_d\to \lsub{d}Y.$
    \item
      \label{typeS:I2}
      $x = I \cdot (a(\Zint_2) \otimes 1) \cdot J$ or $x = I \cdot (1
      \otimes a(\Zpint_2)) \cdot J$. In these cases, $x$ has type
      $\lsub{d}Y\to \lsub{u}Y$ or $Y_u\to Y_d$, respectively.
    \item 
      \label{typeS:I1I3}
      $x=I\cdot (a(\Zint_1\cup\Zint_3)\otimes 1)\cdot J$ or $x=I\cdot
      (1\otimes a'(\Zpint_1\cup\Zpint_3))\cdot J$. In this case, $x$
      has type $\lsub{u}Y\to \lsub{d}Y$ or $Y_d\to Y_u$, respectively.
    \item
      \label{typeS:C}
      $x=I\cdot (a([c_1,c_2])\otimes 1)\cdot J$ or $x=I\cdot (1\otimes
      a'([c_2',c_1']))\cdot J$. In this case, $x$ has type $\star\to
      \star$.
    \item 
      \label{typeS:NoC}
      For some chord $\xi \supset [c_1,c_2]$, we have $x = a(\xi
      \setminus [c_1,c_2]) \otimes a'(\xi) \cdot I$ or $x = a(\xi)
      \otimes a'(\xi \setminus [c_1',c_2']) \cdot I$.  In this case,
      $x$ can be of any of the following possible types: $\star\to
      \star$, $\lsub{u}Y\to Y_u$, $Y_u\to\lsub{u}Y$, $\lsub{d}Y\to
      Y_d$ or $Y_d\to \lsub{d}Y.$
    \item
      \label{typeS:NoI2}
      For some chord $\xi \supset [d,u]$, we have $x = a(\xi \setminus
      [d,u]) \otimes a'(\xi) \cdot I$ or $x = a(\xi) \otimes a'(\xi
      \setminus [d',u']) \cdot I$.  In this case, $x$ has type
      $\lsub{u}Y\to \lsub{d}Y$ or $Y_d\to Y_u$, respectively.
    \item
      \label{typeS:NoCC}
      For some chord $\xi \supset [c_1,c_2]$, we have $x = a(\xi
      \setminus [c_1,c_2]) \otimes a'(\xi\setminus[c_1',c_2']) \cdot
      I$.  In this case, $x$ can be of any of the following possible
      types: $\star\to \star$, $\lsub{u}Y\to Y_u$, $Y_u\to\lsub{u}Y$,
      $\lsub{d}Y\to Y_d$ or $Y_d\to \lsub{d}Y.$
    \item
      \label{typeS:NoI2NoC}
      For some chord $\xi \supset [c_1,c_2]$, we have $x = a(\xi
      \setminus [c_1,c_2]) \otimes a'(\xi\setminus[d',u']) \cdot I$ or
      $x = a(\xi \setminus [d,u]) \otimes a'(\xi\setminus[c_1',c_2'])
      \cdot I$.  In this case, $x$ has type $\lsub{u}Y\to \lsub{d}Y$
      or $Y_d\to Y_u$, respectively.
  \end{enumerate}
\end{definition}

\begin{remark}
  \label{rmk:DegenerateCase}
  In the degenerate case where the basepoint $z$ is next to $c_1$ or $c_2$,
  some of the above near-chords do not exist (i.e., the support of
  these elements automatically contains the basepoint $z$, and hence
  the element does not exist in the algebra). Specifically, 
  elements of Types~(\ref{typeS:NoC}), (\ref{typeS:NoI2}),
  (\ref{typeS:NoCC}) and (\ref{typeS:NoI2NoC}) do not exist.
\end{remark}

\begin{proposition}\label{prop:factor-anti-braid} 
  Suppose $x$ is a basic generator of $\Blg$ (i.e., $x=i a(\xi)
  j\otimes i'a'(\eta)j'$ where $\xi$ and $\eta$ are strands
  diagrams). If $x$ is not an idempotent then $x$ can be factored as a
  product of near-chords.
\end{proposition}
\begin{proof}
  After reversing the roles of $\PMC$ and $\PMC'$, we can assume the initial idempotent is one of
  $\lsub{u}Y$ or $\lsub{d}Y$. This means no strand can start or end at $u'$ or $d'$.

  We call a basic generator $x$ \emph{trimmed} if it satisfies the
  following properties:
  \begin{itemize}
  \item There is no strand ending at $c_2$ or starting at~$c_1$;
    similarly, there is no strand ending at $c_1'$ or starting at~$c_2'$.
  \item There is no strand ending at $u$ or starting at~$d$.
  \end{itemize}
  If the basic generator is not trimmed, we obtain immediately a factorization:
  \begin{enumerate}
  \item If a strand ends at $u$, we can factor off a near-chord of Type~\ref{typeS:I2} (on the right).
  \item If a strand starts at $d$, we can factor off a near-chord of Type~\ref{typeS:I2} (on the left).
  \item If a strand ends at $c_2$ (and none ends at $u$ and none starts at $d$), there are subcases:
    \begin{enumerate}
    \item If the strand ending at $c_2$ starts at $u$, then
      a different strand must terminate at $d$, and we can factor off a near-chord of Type~\ref{typeS:I1I3} (on the right).
    \item If the strand ending at $c_2$ starts at $c_1$, no other strand can end at $c_1$, and none can
      start at $c_1$. Thus, we can factor off a near-chord of Type~\ref{typeS:C} (on the right or the left).
    \item If the strand ending at $c_2$ starts below $c_1$,
      we can factor off a near-chord of Type~\ref{typeS:C} (on the right).
    \end{enumerate}
  \item If a strand starts at $c_1$ (and we are in none of the earlier cases), we can factor off of 
    a near-chord of Type~\ref{typeS:C} (on the left).
  \item If a strand starts at $c_1'$ or ends at $c_2'$, we can factor
    off a near-chord of Type~\ref{typeS:C}.
  \end{enumerate}

    We next show that trimmed basic generators can be factored into
    near-chords of Types (\ref{typeS:Generic}), (\ref{typeS:NoC}),
    (\ref{typeS:NoI2}), (\ref{typeS:NoCC}), and (\ref{typeS:NoI2NoC}) by an adaptation of the argument
    used to prove Lemma~\ref{lem:factor-diag}. To adapt the argument,
    we first give a suitable generalization of
    the notion of ``strands''.
    To this end a {\em generalized strand} is any one of the following things:
    \begin{enumerate}[label=(GS-\arabic*),ref=GS-\arabic*]
    \item a strand starting and ending at points other than $u$, $d$, $u'$, $d'$, $c_1$, $c_2$, $c_1'$, or $c_2'$;
    \item a pair of strands $s_1$ and $s_2$, where $s_1$ terminates in $c_1$ and $s_2$ starts at $c_2$;
    \item a pair of strands $s_1$ and $s_2$, where $s_1$ terminates in $c_2'$ and $s_2$ starts at $c_1'$;
    \item \label{typeGS:NoI2} a pair of strands $s_1$ and $s_2$, where
      $s_1$ terminates in $d$ and $s_2$ starts at $u$.
    \end{enumerate}
    In cases where the generalize strands consists of two strands $s_1$ and $s_2$ (in the above ordering), 
    we say that the generalized strand starts at the starting point
    of~$s_1$ and terminates in the
    terminal point of~$s_2$.
    
    With this understanding, the proof of Lemma~\ref{lem:factor-diag}
    applies to give a factorization of $x$ as a product of algebra
    elements (in the anti-braid subalgebra) whose moving strands consist of matched pairs of
    generalized strands. By ``matched'', we mean here that the initial
    point $p$ of the generalized strand on the $\PMC$ side corresponds
    to a point $p'$ which is the terminal point of the generalized
    strand on the $\PMC'$ side; similarly, the terminal point $q$ of
    the generalized strand on the $\PMC$ side corresponds to a point
    $q'$ which is the initial point of the generalized strand on the
    $\PMC'$ side.

    These pairings give near-chords of Types~(\ref{typeS:Generic}), (\ref{typeS:NoC}),
    (\ref{typeS:NoI2}), (\ref{typeS:NoCC}), and
    (\ref{typeS:NoI2NoC}).
    This completes the proof of Proposition~\ref{prop:factor-anti-braid}.
\end{proof}

\begin{definition}
  A near-chord $x$ is called {\em short} if it is of one of the following types:
  \begin{itemize}
    \item $x$ is of Type~(\ref{typeS:Generic}), where $\xi$ is
      an interval of length one.
    \item $x$ is of Type~(\ref{typeS:I2}).
    \item $x$ is of Type~(\ref{typeS:I1I3}).
    \item $x$ is of Type~(\ref{typeS:NoCC}), where $\xi$ has length
      $5$ (so each connected component of the support of $x$ has
      length $1$).
    \end{itemize}
\end{definition}

\begin{figure}
\begin{center}
\input{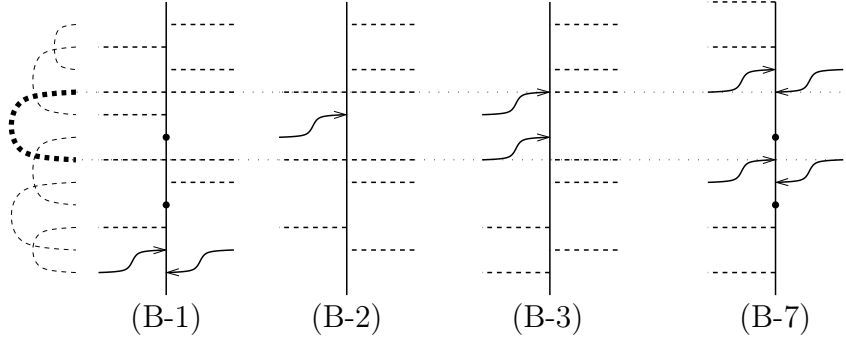}
\end{center}
\caption {{\bf Short anti-braid-like near-chords.}
  \label{fig:ShortABLNearChords}}
\end{figure}

See Figure~\ref{fig:ShortABLNearChords} for an illustration.

\begin{proposition}
  \label{prop:NonBraid} Let $\PMC$ be an unexceptional pointed matched
circle.  There is a unique solution to the equation $dA_0+A_0\cdot A_0 = 0$,
where $A_0$ is in the anti-braid-like subalgebra~$\Blg$, satisfying the
following properties:
  \begin{itemize}
  \item $A_0$ contains each short near-chord with multiplicity one.
  \item $A_0$ is in grading $-1$.
  \end{itemize} The solution $A_0$ is simply the sum of all near-chords,
  in the sense of Definition~\ref{def:NearChords}.
\end{proposition}

\begin{proof}
  The proof is similar to several arguments in~\cite{LOT4}, and we
  will give it expeditiously.

  The fact that $A_0$ contains generic chords of arbitrary length, but
  on one side of $[c_1,c_2]$, follows from a simple induction on the
  length. More precisely, if $a$ is a near-chord of
  Type (\ref{typeS:Generic}) that does not contain $[c_1,c_2]$ and is
  not short, then $da\neq 0$, and
  $da$ is a sum of
  products of elements of Type~(\ref{typeS:Generic}) whose support
  is smaller.  Each of these products occurs once in $A_0\cdot A_0$, and
  can not be canceled by any other term in $dA_0$ (see the proof
  of~\cite[Proposition~\ref*{HFa:thm:DDforIdentity}]{LOT4}).  This uses
  the hypothesis that $\PMC$ is unexceptional -- a special
  length $3$ chord gives rise to an algebra element with $da=0$. 

  This same argument can be used to show that all terms of
  Type~(\ref{typeS:NoCC}) appear in the differential, by
  induction of the length of the support: minimal ones
  exist by hypothesis, and multiplying one by an element of
  Type~(\ref{typeS:Generic}) gives a product with no alternative
  factorization, but which appears in the differential of another term
  of Type~(\ref{typeS:NoCC}) (and not in the differential of any
  other basic algebra element).

  Taking the product of elements of Type~(\ref{typeS:I1I3}) with
  of Type~(\ref{typeS:NoCC}), we get elements with unique non-trivial
  factorizations, and which appear in the differentials of elements of
  Type~(\ref{typeS:NoI2NoC}). This ensures the existence of all
  near-chords of Type~(\ref{typeS:NoI2NoC}) in the differential.

  Taking products of elements of Types~(\ref{typeS:I2})
  and (\ref{typeS:I1I3}), we get elements which cannot be factored
  in any other way, and which appear in the differentials of near-chords of
  Type~(\ref{typeS:C}) (and not in the differentials of any other
  basic algebra elements). This forces the existence of near-chords of
  Type~(\ref{typeS:C}).

  Taking products of elements of Type~(\ref{typeS:C}) with elements
  of Type~(\ref{typeS:NoCC}), we get elements which appear in
  the differentials of elements of Type~(\ref{typeS:NoC}) and which
  have no other factorizations, and appear in no other differentials. This forces
  elements of Type~(\ref{typeS:NoC}) to appear in the differential.

  Products of Type~(\ref{typeS:I1I3}) with elements of
  Type~(\ref{typeS:NoC}) we get elements which are in the
  differentials of elements of Type~(\ref{typeS:NoI2}).  
  These force the existence of elements of Type~(\ref{typeS:NoI2}).
  (Other terms in the differential of elements of
  Type~(\ref{typeS:NoI2}) arise as products of elements
  of Type~(\ref{typeS:NoI2NoC}) with elements of
  Type~(\ref{typeS:C}).)

  Finally, products of elements of
  Type~(\ref{typeS:I2}) with elements of Type~(\ref{typeS:NoI2}),
  and also products of elements of Types~(\ref{typeS:C}) with
  elements of Type~(\ref{typeS:NoC}), give terms in the differentials
  of elements of Type~(\ref{typeS:Generic}) which cross $[c_1,c_2]$.
  (Other terms in these differentials are gotten as products of pairs
  of chords of Type~(\ref{typeS:Generic}).)

  Results of this argument are summarized by the following table,
  which shows how differentials cancel against products in $dA_0 + A_0
  \cdot A_0$.  (There are
  also  cancellations of products against products, which are not
  shown.)
 \begin{center}
  \begin{tabular}{c|cccccccc}
   $\times$ &(\ref{typeS:Generic}) & (\ref{typeS:I2}) &(\ref{typeS:I1I3}) &(\ref{typeS:C}) & (\ref{typeS:NoC}) &(\ref{typeS:NoI2}) &(\ref{typeS:NoCC}) & (\ref{typeS:NoI2NoC})  \\
    \hlx{hv}
   (\ref{typeS:Generic})  & $d$(\ref{typeS:Generic}) & \blank & \blank & 
   \blank & $d$(\ref{typeS:NoC})  & $d$(\ref{typeS:NoI2}) & $d$(\ref{typeS:NoCC}) & $d$(\ref{typeS:NoI2NoC}) \\
   (\ref{typeS:I2}) & \blank & \blank & $d$(\ref{typeS:C}) & \blank & \blank & $d$(\ref{typeS:Generic}) & \blank & $d$(\ref{typeS:NoC}) \\
   (\ref{typeS:I1I3}) & \blank & $d$(\ref{typeS:C}) & \blank & \blank & 
   $d$(\ref{typeS:NoI2}) & \blank & $d$(\ref{typeS:NoI2NoC}) & \blank \\
   (\ref{typeS:C}) & \blank & \blank & \blank & \blank & $d$(\ref{typeS:Generic}) & 
   \blank & $d$(\ref{typeS:NoC}) & $d$(\ref{typeS:NoI2}) \\
   (\ref{typeS:NoC}) & $d$(\ref{typeS:NoC}) & \blank & $d$(\ref{typeS:NoI2}) & $d$(\ref{typeS:Generic}) & \blank & \blank & \blank & \blank \\
   (\ref{typeS:NoI2}) & $d$(\ref{typeS:NoI2}) & $d$(\ref{typeS:Generic}) & \blank &\blank & \blank & \blank & \blank & \blank  \\
   (\ref{typeS:NoCC}) & $d$(\ref{typeS:NoCC}) & \blank & $d$(\ref{typeS:NoI2NoC}) & $d$(\ref{typeS:NoC}) & \blank & \blank & \blank & \blank\\
   (\ref{typeS:NoI2NoC}) & $d$(\ref{typeS:NoI2NoC}) & $d$(\ref{typeS:NoC}) & \blank & $d$(\ref{typeS:NoI2}) & \blank & \blank &\blank & \blank 
  \end{tabular}
  \end{center}

  We have shown that $A_0$ contains the sum of all near-chords in the
  sense of Definition~\ref{def:NearChords}, and that the sum of all
  near chords is a solution. By Proposition~\ref{prop:factor-anti-braid},
  any non-idempotent basic element of $\Blg$
  factors as a product of such near-chords; it follows that all other
  algebra elements have grading less than~$-1$, and thus cannot appear
  in~$A_0$. This proves that the
  sum of the near-chords is the unique solution.
\end{proof}

Let $A_0$ be the sum of the near-chords (from Definition~\ref{def:NearChords}).
Let $M$ denote the $\Blg$-module $\Blg$ with differential given by
$\nabla_0(B)=dB+B\cdot A_0$.

\begin{figure}
  \begin{center}
    \input{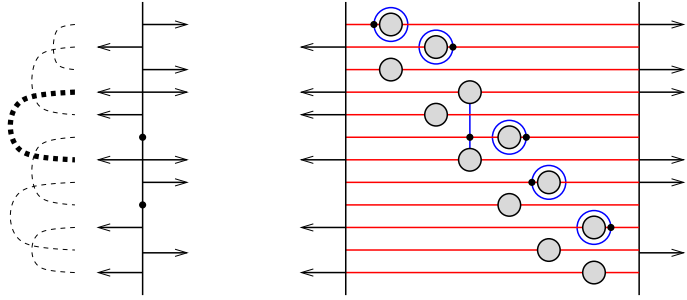}
  \end{center}
  \caption {{\bf Heegaard diagrams.}
    \label{fig:HeegaardNonBraid}
    At the left, we have illustrated an idempotent; at the right,
    a generator in the Heegaard diagram.}
\end{figure}

\begin{proposition}
  \label{prop:CalculateVertex}
  The bimodule $\CFDDa(\HD_0(\gamma))$ is isomorphic to
  $(\Alg(\PMC)\otimes\Alg(\PMC'))\otimes_\Blg M$.
\end{proposition}

\begin{proof}
  There is an obvious identification of generators of
  $\CFDDa(\HD_0(\gamma))$ with basic idempotents of $\Blg$; see
  Figure~\ref{fig:HeegaardNonBraid}. It is immediate from the topology
  of domains in the Heegaard diagram $\HD_0$ that the coefficients of
  the differential on $\CFDDa(\HD_0(\gamma))$ lie in $\Blg$; that is,
  if $\x\in\Gen(\HD)$ is a generator of $\CFDDa(\HD_0(\gamma))$ then
  $\bdy x=\sum_{\y\in\Gen(\HD)} a_{\x,\y}\y$ where each $a_{\x,\y}$
  lies in $\Blg$. It follows that
  $\CFDDa(\HD_0(\gamma))=(\Alg(\PMC)\otimes\Alg(\PMC'))\otimes_\Blg N$
  where $N$ is some free, rank $1$ \dg module over $\Blg$.

  \begin{figure}
    \centering
    \includegraphics[scale=.5]{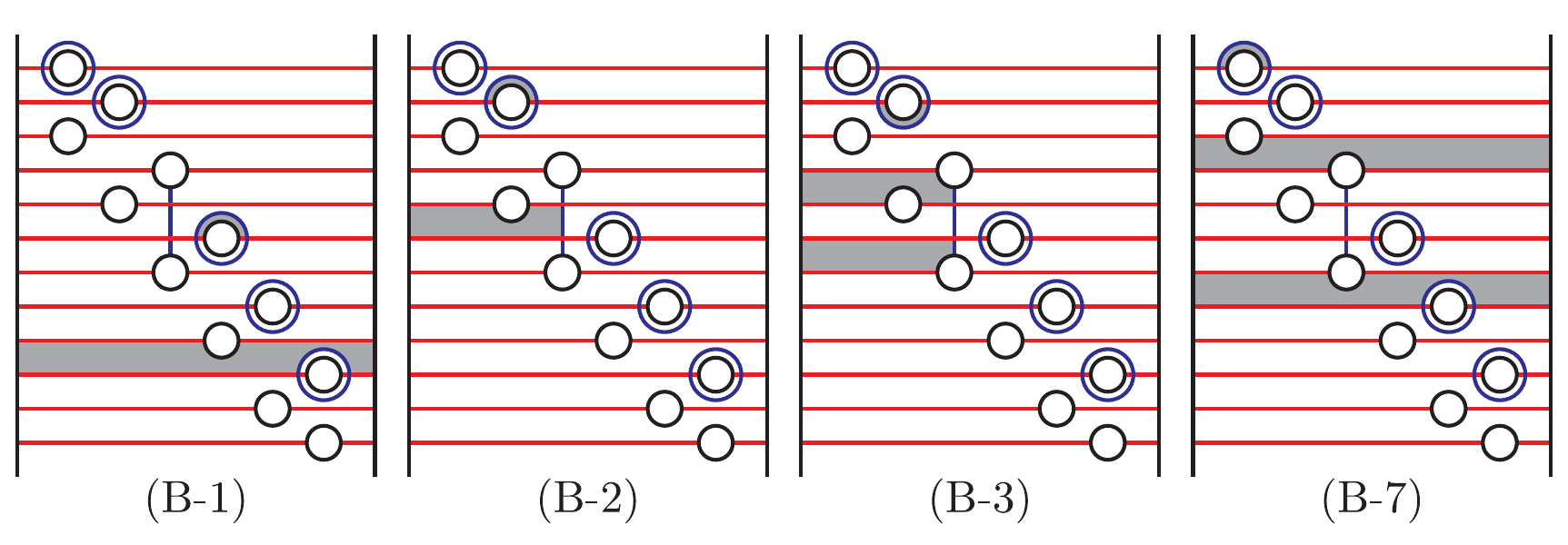}
    \caption{\textbf{Domains for short near chords.} This is the same
      Heegaard diagram as Figure~\ref{fig:HeegaardDiagrams}, but here
      we are using shading to indicate domains corresponding to
      the four types of short near-chords.}
    \label{fig:ShortNearDomains}
  \end{figure}

  Next, we check that each short near-chord contributes to the
  differential on $\CFDDa(\HD_0(\gamma))$ (and hence to the
  differential on $N$). There is a unique domain in $\HD_0(\gamma)$
  corresponding to each short near-chord; see
  Figure~\ref{fig:ShortNearDomains}. These domains are topological
  disks, so it is easy to check that they contribute to the
  differential.

  Thus, it follows from the uniqueness statement in
  Proposition~\ref{prop:NonBraid} that the differential on $N$ is the
  same as the differential $\nabla_0$ on $M$ (given by
  $\nabla_0(B)=dB+B\cdot A_0$).
\end{proof}

\subsection{Anti-braid-like resolutions: Degenerate case}
\label{subsec:DegenerateVertex}

The above discussion is adequate to describe the type \DD\ module associated to the anti-braid-like resolution of $B_i$ when $i$
is not extremal (i.e., top-most or
bottom-most). That is, for the linear pointed matched
circle, the above discussion is sufficient to describe $Y_{0(\gamma_n)}$, where
$\gamma_n$ is any of the preferred curves specified by the pointed
matched circle {\em except} when $n=1$ or $4k$. For those two cases, the
matched pair $c_1$ and $c_2$ are separated by one, rather than
two, marked points. This causes some collapse of the idempotents, and
immediately excludes several families of near-chords (those of
Types~(\ref{typeS:I2}), (\ref{typeS:I1I3}), (\ref{typeS:NoI2}),
and (\ref{typeS:NoI2NoC})). Further, in these cases, the basepoint
is also next to $c_1$ and $c_2$, which also removes several other
families of near-chords---Types~(\ref{typeS:NoC}) and~(\ref{typeS:NoCC}); see also Remark~\ref{rmk:DegenerateCase}.

For the reader's convenience, we spell out precisely what is left over
for the case of~$\gamma_1$, in which case we have the following situation.

\begin{hypothesis}
  \label{hyp:DistanceOne}
  We consider a matched pair $\{c_1,c_2\}$ in $\PMC$, where $c_1$ is just above the basepoint $z$, 
  and $c_1$ and $c_2$ are separated by a single distinguished point, which we denote $p$.
\end{hypothesis}

We call the corresponding curve $\gamma$ {\em degenerate}.

In this case, we have a much simpler analogue of Proposition~\ref{prop:CalculateVertex}. To 
state it, we describe first a corresponding version of the anti-braid subalgebra and near-chords.

\begin{definition}
  \label{def:degen:AntiBraidSubalgebra}
  The distinguished point $p$ (respectively $p'$) divides
  the interval between $c_1$ and $c_2$ (respectively $c_1'$ and $c_2'$) into
  two consecutive intervals, labelled $\Zint_1$ and $\Zint_2$ (respectively
  $\Zpint_1$ and $\Zpint_2$). 
  The {\em degenerate anti-braid subalgebra} 
  $\BlgDeg$ is the subalgebra of $\Alg(\PMC)\otimes \Alg(\PMC')$ specified as
  follows.
  \begin{itemize}
  \item Basic idempotents correspond to pairs $\SetS$ and $\SetT$ of subsets
    of matched pairs in $\PMC$ (where $\SetT$ is thought of as a set of
    pairs for $\PMC'$) with the property that:
    \begin{itemize}
    \item The union $\SetS\cup\SetT$, contains all matched pairs
      {\em except} for the one containing $p$.
    \item The intersection $\SetS\cap\SetT$ consists of exactly the
      matched pair $\{c_1,c_2\}$.
    \end{itemize}
  \item The non-idempotent elements are linear combinations of elements of the
    form $a\otimes a'$, constrained as follows:
    \begin{itemize}
    \item There are basic idempotents $i\otimes i'$ and $j\otimes j'$
      in $\BlgDeg$ as above
      with $iaj=a$ and $i'a'j'=a'$.
    \item The local multiplicities of $a$ in $\PMC$
      at any region outside the interval
      $[c_1,c_2]$ coincide with the local multiplicities of $a'$
      at the corresponding region in~$\PMC'$.
    \item The local multiplicity of $a$ at $\Zint_1$ equals the local
      multiplicity of $a$ at $\Zint_2$, and the local multiplicity
      of $a'$ at $\Zpint_1$ equals the local multiplicity of $a'$
      at~$\Zpint_2$.
    \end{itemize}
  \end{itemize}
\end{definition}

\begin{definition}
  \label{def:NearChordsDeg}
  A non-zero element $x$ of $\BlgDeg$ is a \emph{near-chord} if $x$
  satisfies one of the following conditions:
  \begin{enumerate}[label=(B$_{\degen}$-\arabic*),ref=B$_{\degen}$-\arabic*]
    \item
      \label{typeSD:Generic}
      $x=J\cdot (a(\xi)\otimes
      a'(\xi))\cdot I $, where here $\xi$ is some interval
      neither of whose endpoints is $c_1$, $c_2$, or $p$.
      (This is the analogue of Type~(\ref{typeS:Generic}).)
    \item
      \label{typeSD:C}
      $x=I\cdot (a([c_1,c_2])\otimes 1)\cdot J$ or $x=I\cdot (1\otimes
      a'([c_2',c_1']))\cdot J$. (This is the analogue of Type~(\ref{typeS:C}).)
  \end{enumerate}
\end{definition}

Let $A_{\degen}$ be the sum of the degenerate near-chords (from Definition~\ref{def:NearChordsDeg}).
Let $M_{\degen}$ denote the $\BlgDeg$-module $\BlgDeg$ with differential given by
$\nabla_\degen(\cdot) = d(\cdot)+\cdot A_\degen$, i.e., $\nabla_\degen(B)=dB+B\cdot A_{\degen}$.

\begin{proposition}
  \label{prop:CalculateVertexDeg}
  The bimodule $\CFDDa(\HD_0(\gamma))$ is isomorphic to
  $(\Alg(\PMC)\otimes\Alg(\PMC'))\otimes_\Blg M_{\degen}$.
\end{proposition}

\begin{proof}
  Follow the proof of Proposition~\ref{prop:CalculateVertex},
  and note that the region in the Heegaard diagram just below $c_1$, and hence also the region
  just above $c_2$, contains the basepoint.
\end{proof}

\begin{remark}
  \label{rmk:ChangeHypothesis}
  In the above, we have concentrated on the case of $\gamma_1$.  The
  case of $\gamma_{4k}$ can be described similarly, with a
  straightforward change of notation; in particular, in this case
  requirement that $c_1$ be just above the basepoint is replaced by
  the requirement that $c_2$ be just below the basepoint, with minor
  notational changes.  
\end{remark}

\begin{example}\label{exam:abl-torus}
  Consider the case $k=1$ and $n=4$. Continuing with the notation for
  the torus algebra from Example~\ref{exam:DD-id-torus}, the
  degenerate anti-braid subalgebra $\BlgDeg$ has a single idempotent
  $I=\iota_1\otimes \iota_0$. Note that $\Zint_1=[2,3]$ and
  $\Zint_2=[3,4]$, while $\Zpint_1=[2,3]$ and $\Zpint_2=[1,2]$
  (remember the orientation reversal). Thus, the non-idempotent
  generators are $\sigma_{23}\otimes \iota_0$, $\iota_1\otimes
  \rho_{12}$, and $\sigma_{23}\otimes\rho_{12}$. Of these, the first
  two are near-chords. So, by
  Proposition~\ref{prop:CalculateVertexDeg}, $\CFDDa(\HD_0(\gamma))$
  is generated by $I$ with differential
  \[
  \bdy(I)=(\sigma_{23}\otimes 1+1\otimes\rho_{12}) I.
  \]
  (Actually, Proposition~\ref{prop:CalculateVertexDeg} does not quite
  apply, because our pointed matched circle is not
  unexceptional. However, a slightly longer argument, using the fact
  that the bimodule $\CFDDa(Y_{0(\gamma)})$ is \emph{stable} in an
  analogous sense to~\cite[Definition~\ref*{HFa:def:StableModule}]{LOT4},
  shows that the result holds in this setting, as well. See also Remark~\ref{rem:stable}.)
\end{example}


\section{Edges}
\label{sec:edges}

In the computation of the bimodule associated to the branched
double-cover of the $B_i$ ($i=1,\dots,n$) as a mapping
cones, we computed the \DD-bimodules appearing in these mapping cones in
Section~\ref{sec:vertices}. We next
calculate the morphisms between these bimodules appearing in the
mapping cones.

More precisely, we compute the maps $F^\pm_\gamma$ from
Theorem~\ref{thm:Dehn-is-MC}, where $\gamma$ is one of the
distinguished curves in the linear pointed matched circle, subject to
Hypothesis~\ref{hyp:DistanceTwo} (or possibly
Hypothesis~\ref{hyp:DistanceOne}, in
Section~\ref{subsec:DegenerateEdge}).

\subsection{Left-handed cases}
\label{subsec:LeftHanded}

We start with the map $F^-_\gamma\co \CFDDa(\Id)\to
\CFDDa(Y_{0(\gamma)})$. Recall from Section~\ref{sec:vertices} that
$\CFDDa(\Id)$ is induced from a free, rank-one module over the diagonal
subalgebra $\Dlg$; as such, its differential is given by 
\[
\nabla_\Id(B)=d(B)+B\cdot A_\Id
\]
where $A_\Id\in B$ is a particular element, computed in
Section~\ref{sec:braid-like}. Similarly, the module
$\CFDDa(Y_{0(\gamma)})$ is induced by a free, rank-one module over the
anti-braid subalgebra $\Blg$, and so its differential is given by
\[
\nabla_0(B)=d(B)+B\cdot A_0
\]
where $A$ is a particular element, computed in
Section~\ref{sec:anti-braid-like}.

The homomorphism $F^-_\gamma$ is determined by the images of the
idempotents in $\Dlg$. The image $F^-_\gamma(I)$ must have the same
left idempotent as $I$, and its right idempotent must lie in
$\Blg$. Further, the definition of $F^-_\gamma$ gives restrictions
on the supports of the algebra elements it outputs; see
Proposition~\ref{prop:CalculateNegMorphism}. (Alternatively, these
restrictions can be obtained from grading considerations.) Taken
together, these restrictions imply that $F^-_\gamma$ is determined by
an element in a certain bimodule, which in turn turns out to be a subalgebra:

\begin{definition}\label{def:MlgMinus}
  The {\em left-handed morphism subalgebra} $\Mlg^-$ is the
  subalgebra of $\Alg(\PMC)\otimes \Alg(\PMC')$ (where $\PMC'=-\PMC$)
  which is spanned by algebra elements $a\otimes a'$ satisfying the
  following properties:
  \begin{itemize}
  \item $I\cdot a \cdot J=a$ and $I'\cdot a'\cdot J'=a'$ where
    $I\otimes I'$ is an idempotent for the diagonal subalgebra, while
    $J\otimes J'$ is an idempotent for the anti-braid subalgebra.
  \item the local multiplicities of $a$ outside the interval
    $[c_1,c_2]$ agree with the corresponding local multiplicities of
    $a'$ outside the interval $[c_1',c_2']$.
  \end{itemize}
\end{definition}

\begin{figure}
  \centering
  \includegraphics{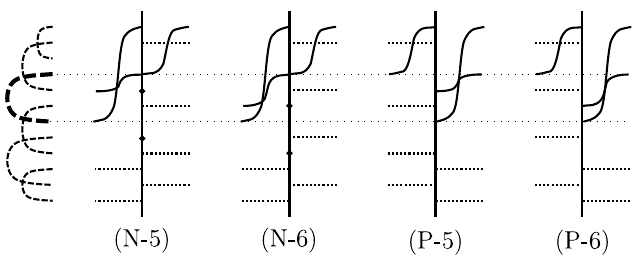}
  \caption{\textbf{Additional near-chords in for morphisms in the
      untruncated case.} The additional near-chords are shown for both
    the left-handed and right-handed morphism algebras.  We only show
    one case of each type of near-chord; the other case is gotten by
    rotating the diagram by $180^\circ$.}
  \label{fig:untrunc-near-chords}
\end{figure}

\begin{remark}
  Recall that we are using $\Alg(\PMC)$ to denote the truncated
  algebra associated to a surface, as in Section~\ref{sec:Trunc}; it
  is the algebra denoted $\Alg'(\PMC)$ in~\cite{LOT2}, and it is
  smaller than the algebra $\Alg(\PMC)$ of~\cite{LOT1}.  In the
  untruncated case, there would be two more types of near-chords
  (Definition~\ref{def:lhnc}), as shown in
  Figure~\ref{fig:untrunc-near-chords}, and the proof of
  Proposition~\ref{prop:factor-M-minus} would be correspondingly
  rather more complicated.
\end{remark}

Although, technically, $\Mlg^-$ is a subalgebra, we will not consider
multiplying two elements of $\Mlg^-$, as the product vanishes
identically.  (No idempotent of $\Blg$ is also an idempotent of
$\Dlg$; in particular $\Mlg^-$ is not unital.)
Rather, we will multiply
elements of $\Mlg^-$ by elements of $\Blg$ and $\Dlg$:

\begin{lemma}
  $\Dlg\cdot \Mlg^-\cdot \Blg =\Mlg^-$, where $\Blg$ is the
  anti-braid subalgebra of Definition~\ref{def:AntiBraidSubalgebra},
  $\Dlg$ is the diagonal subalgebra, and $\Mlg^-$
  is the left-handed morphism subalgebra.
\end{lemma}
\begin{proof}
  This is immediate from the definitions.
\end{proof}

\begin{lemma}\label{lem:M-minus-Z-graded}
  There is a $\ZZ$-grading on $\Mlg^-$ which is compatible with the
  $\ZZ$-gradings on $\Blg$ and $\Dlg$ from
  Lemmas~\ref{lem:diagonal-Z-graded} and~\ref{lem:anti-braid-Z-graded}, in the sense that
  $\gr(d\cdot m\cdot b)=\gr(d)+\gr(m)+\gr(b)$ for $b\in \Blg$, $m\in
  \Mlg^-$ and $d\in\Dlg$.
\end{lemma}
\begin{proof}
  Consider the type \DD\ bimodule $\CFDDa(\Cone(F^-_\gamma))$. The
  coefficient algebra, in the sense of~\cite[Definition~\ref*{HFa:def:coeff-alg}]{LOT4}, of
  $\CFDDa(\Cone(F^-_\gamma))$, is the set of $2\times 2$ matrices
  $
  \left(\begin{smallmatrix}
    d & m\\
    0 & b
  \end{smallmatrix}\right)
  $ where $d\in \Dlg$, $b\in\Blg$, and $m$ is
  in~$\Mlg^-$. Additionally, the grading set for the bimodule
  $\CFDDa(\Cone(F^-_\gamma))$ is $\lambda$-free; this follows from the
  fact that there are no provincial triply-periodic domains in a
  bordered Heegaard triple diagram for $F^-_\gamma$. (Here, we are
  using the fact that $F^-_\gamma$ is triangle-like,
  Lemma~\ref{lem:triangle-like}.)  Then
  \cite[Lemma~\ref*{HFa:lem:grading-coeff-implies}]{LOT4} proves the
  result.
\end{proof}
Note that Lemma~\ref{lem:M-minus-Z-graded} does not assert uniqueness
of the grading; we will prove uniqueness as a relative grading in
Corollary~\ref{cor:Mlg-minus-Z-gr-unique}.

The idempotents of $\Dlg$ are partitioned into eight types,
$\lsub{duc}X$, $\lsub{du}X_c$, $\lsub{dc}X_u$, $\lsub{uc}X_d$,
$\lsub{d}X_{uc}$, $\lsub{u}X_{dc}$, $\lsub{c}X_{dc}$, $X_{duc}$,
according to which side the matched pairs containing $d$, $u$ and
$c_1$ (and $c_2$) are occupied on. We will say that an
element $x\in\Mlg^-$ has \emph{type $\lsub{u}X_{dc}\to Y_d$}, for
example, if $x=I\cdot x\cdot J$ for idempotents $I\in\Dlg$ and
$J\in\Blg$ where $I$ has type $\lsub{u}X_{dc}$ and $J$ has type $Y_d$.

\begin{figure}
  \centering
  \includegraphics{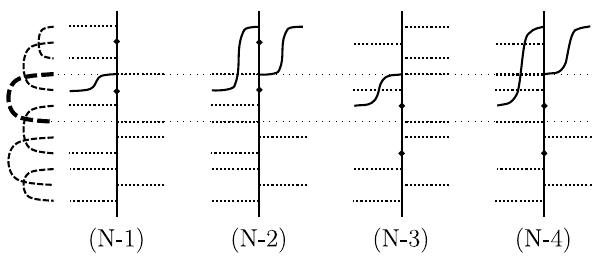}
  \caption{\textbf{Left-handed near-chords.} We have drawn the first
    case of each type of near chord from Definition~\ref{def:lhnc}.}
  \label{fig:lhnc}
\end{figure}

\begin{definition}\label{def:lhnc}
  An element $x\in \Mlg^-$ is called a \emph{(left-handed) near-chord}
  if it has one of the following forms:
  \begin{enumerate}[label=(N-\arabic*),ref=N-\arabic*]
  \item 
    \label{type:l:I1}
    $x=I\cdot (a(\Zint_3)\otimes 1)\cdot J$ or $x=I\cdot (1\otimes
    a(\Zpint_1))\cdot J$. In the first case, $x$ has type
    $\lsub{ud}X_c\to \lsub{d}Y$ or $\lsub{u}X_{dc}\to Y_d$, while in
    the second case $x$ has type $\lsub{ud}X_c\to \lsub{u}Y$ or
    $\lsub{d}X_{uc}\to Y_u$.
  \item 
    \label{type:l:I1+}
    $x=I\cdot (a(\Zint_3\cup\xi)\otimes a'(\xi))\cdot J$, where $\xi$
    is a chord starting at $c_2$; or $x=I\cdot (a(\xi)\otimes
    a'(\Zpint_1\cup\xi))\cdot J$, where $\xi$ is a chord ending at
    $c_1$. In the first case, $x$ has type $\lsub{ud}X_c\to \lsub{d}Y$
    or $\lsub{u}X_{dc}\to Y_d$, while in the second case $x$ has type
    $\lsub{ud}X_c\to \lsub{u}Y$ or $\lsub{d}X_{uc}\to Y_u$.
  \item 
    \label{type:l:I1I2}
    $x=I\cdot (a(\Zint_2\cup\Zint_3)\otimes 1)\cdot J$ or $x=I\cdot (
    1\otimes a'(\Zpint_1\cup\Zpint_2))\cdot J$. In the first case, $x$
    has type $\lsub{duc}X \to \lsub{d}Y$ or $\lsub{uc}X_d \to Y_d$
    while in the second case $x$ has type $X_{duc} \to Y_u$ or
    $\lsub{u}X_{dc} \to \lsub{u}Y$.
  \item 
    \label{type:l:I1I2+}
    $x = I \cdot (a(\Zint_2 \cup \Zint_3 \cup \xi) \otimes a'(\xi))
    \cdot J$ where $\xi$ is a chord starting at $c_2$, or $\xi = I
    \cdot (a(\xi) \otimes a'(\xi \cup \Zpint_1 \cup \Zpint_2)) \cdot
    J$ where $\xi$ is a
    chord ending at
    $c_1$. The first of these corresponds to $\lsub{duc}X \to
    \lsub{u}Y$ or $\lsub{uc}X_{d} \to Y_d$ while the second to
    $\lsub{u}X_{dc} \to \lsub{u}Y$ or $X_{duc} \to Y_u$.
  \end{enumerate}
  See Figure~\ref{fig:lhnc}.

  The \emph{structure constant for the left-handed morphism}, which we
  write $F^-$, is the sum of all (left-handed) near-chords in~$\Mlg^-$.
\end{definition}

\begin{proposition}\label{prop:factor-M-minus} 
  Every element $m$ of $\Mlg^-$ can be factored as a product $m=d\cdot
  m'\cdot b$ where $d\in\Dlg$, $b\in\Blg$, and $m'$ is a near-chord for $\Mlg^-$.
\end{proposition}
\begin{proof}
  We divide the proof into cases, according to the initial and terminal idempotents of our 
  algebra element. After rotating the diagram, we can assume that the terminal
  idempotent is of type $\lsub{u}Y$ or $\lsub{d}Y$. 

  Fixing the initial and terminal idempotent types determines many of
  the local multiplicities
  in a neighborhood of $[c_1,c_2]\times [c_2',c_1']$, with only a few
  choices. We encode these local multiplicities
  in a matrix form
  \[\left[\begin{smallmatrix}
      m_1 & m_1'\\
      m_2 & m_2'\\
      m_3 & m_3'\\
      m_4 & m_4'\\
      m_5 & m_5'
    \end{smallmatrix}\right],\] where
  $m_1$ is the local multiplicity of our element just above $c_2$,
  $m_2$ is the local multiplicity just above $u$, $m_3$ is the local
  multiplicity just above $d$, $m_4$ is the local multiplicity just
  above $c_1$, and $m_5$ is the local multiplicity just below $c_1$;
  and the $m_i'$ are the corresponding local multiplicities on the
  $\PMC'$ side. We annotate these blocks of multiplicities with a $V$
  or an $E$, as needed, where here ``$V$'' means that there is a
  strand on the left terminating at the point $v$ matched with $u$,
  while the notation ``$E$'' means that there is a strand on the left
  terminating at the point $e$ matched with $d$. The matrix and
  marking is sufficient to reconstruct the idempotent types of the
  start and finish; moreover, each matrix marked by $E$ or $V$ occurs
  also as a matrix without any marking.
 
  These marked matrices are shown in
  Table~\ref{tab:Mlg-minus-mults}. Some of the matrix entries are
  $\epsilon$ or $\delta$, which may each be $0$ or $1$, but all the
  $\epsilon$'s (respectively $\delta$'s) in each matrix must be the
  same. So, for instance, there are $5$ different matrices for the
  case $\lsub{dc}X_u\to \lsub{u}Y$.

  (If we were working in the untruncated algebra from,
  say, \cite{LOT1}, the corresponding case analysis would have to
  include multiplicities greater than $1$, as well.)

    \begin{table}
    \centering
   \begin{tabular}{l@{\quad}|@{\quad}l@{\qquad}l}
     & $\to \lsub{u}Y$ & $\to \lsub{d}Y$\\
     \hlx{vhvv}
     $\lsub{duc}X\to$ & 
     $\left[
       \begin{smallmatrix}
         1 & 1\\
         1 & \epsilon\\
         1 & \epsilon\\
         0 & \epsilon\\
         0 & 0
       \end{smallmatrix}\right]$
   & 
     $\left[
       \begin{smallmatrix}
         1 & 1\\
         1 & \epsilon\\
         0 & \epsilon\\
         0 & \epsilon\\
         0 & 0
       \end{smallmatrix}\right]$
     \\
     \hlx{vv}
     $\lsub{du}X_c\to $ &
     $\left[
       \begin{smallmatrix}
         \delta & \delta\\
         1 & \epsilon\\
         1 & \epsilon\\
         0 & \epsilon\\
         \delta & \delta
       \end{smallmatrix}\right]$
     &
     $\left[
       \begin{smallmatrix}
         \delta & \delta\\
         1 & \epsilon\\
         0 & \epsilon\\
         0 & \epsilon\\
         \delta & \delta
       \end{smallmatrix}\right]$
     \\
     \hlx{vv}
     $\lsub{dc}X_u\to $ & 
     $\left[
       \begin{smallmatrix}
         \delta & \delta\\
         0 & 0\\
         1 & 1\\
         0 & 1\\
         \delta & \delta
       \end{smallmatrix}\right]$,
     $\left[
       \begin{smallmatrix}
         1 & 1\\
         1 & \epsilon\\
         1 & \epsilon\\
         0 & \epsilon\\
         0 & 0
       \end{smallmatrix}\right]\!\!V$
     & 
     $\left[
       \begin{smallmatrix}
         \delta & \delta\\
         \epsilon & 0\\
         \epsilon & 1\\
         \epsilon & 1\\
         \delta & \delta
       \end{smallmatrix}\right]$
     \\
     \hlx{vv}
     $\lsub{d}X_{uc}\to$ &
          $\left[
       \begin{smallmatrix}
         0 & 0\\
         0 & 0\\
         1 & 1\\
         0 & 1\\
         1 & 1
       \end{smallmatrix}\right]$,
     $\left[
       \begin{smallmatrix}
         \delta & \delta\\
         1 & \epsilon\\
         1 & \epsilon\\
         0 & \epsilon\\
         \delta & \delta
       \end{smallmatrix}\right]\!V$
     &
     $\left[
       \begin{smallmatrix}
         0 & 0\\
         \epsilon & 0\\
         \epsilon & 1\\
         \epsilon & 1\\
         1 & 1
       \end{smallmatrix}\right]$
     \\
     \hlx{vv}
     $\lsub{uc}X_d\to $ &
     $\left[
       \begin{smallmatrix}
         \delta & \delta\\
         \epsilon & 0\\
         \epsilon & 0\\
         \epsilon & 1\\
         \delta & \delta
       \end{smallmatrix}\right]$
     & 
     $\left[
       \begin{smallmatrix}
         \delta & \delta\\
         1 & 0\\
         0 & 0\\
         1 & 1\\
         \delta & \delta
       \end{smallmatrix}\right]$,
     $\left[
       \begin{smallmatrix}
         1 & 1\\
         1 & \epsilon\\
         0 & \epsilon\\
         0 & \epsilon\\
         0 & 0
       \end{smallmatrix}\right]\!E$
     \\
     \hlx{vv}
     $\lsub{u}X_{dc}\to$ &
     $\left[
       \begin{smallmatrix}
         0 & 0\\
         \epsilon & 0\\
         \epsilon & 0\\
         \epsilon & 1\\
         1 & 1
       \end{smallmatrix}\right]$
     & 
     $\left[
       \begin{smallmatrix}
         0 & 0\\
         1 & 0\\
         0 & 0\\
         1 & 1\\
         1 & 1
       \end{smallmatrix}\right]$,
     $\left[
       \begin{smallmatrix}
         \delta & \delta\\
         1 & \epsilon\\
         0 & \epsilon\\
         0 & \epsilon\\
         \delta & \delta
       \end{smallmatrix}\right]\!E$
     \\
     \hlx{vv}
     $\lsub{c}X_{du}\to$ &
     $\left[
       \begin{smallmatrix}
         \delta & \delta\\
         \epsilon & 0\\
         \epsilon & 0\\
         \epsilon & 1\\
         \delta & \delta
       \end{smallmatrix}\right]\!V$
     & 
     $\left[
       \begin{smallmatrix}
         \delta & \delta\\
         \epsilon & 0\\
         \epsilon & 1\\
         \epsilon & 1\\
         \delta & \delta
       \end{smallmatrix}\right]\!E$
     \\
     \hlx{vv}
     $X_{duc}\to$ & 
     $\left[
       \begin{smallmatrix}
         0 & 0\\
         \epsilon & 0\\
         \epsilon & 0\\
         \epsilon & 1\\
         1 & 1
       \end{smallmatrix}\right]\!V$
     &
     $\left[
       \begin{smallmatrix}
         0 & 0\\
         \epsilon & 0\\
         \epsilon & 1\\
         \epsilon & 1\\
         1 & 1
       \end{smallmatrix}\right]\!E$
     \\
     \hlx{v}
   \end{tabular}
    \caption{\textbf{Local multiplicities of elements of $\Mlg^-$.}}
    \label{tab:Mlg-minus-mults}
  \end{table}

  We will examine most of the cases separately.  But first, observe
  that, by the argument from Lemma~\ref{lem:factor-diag}, if some
  strand in $m$ on the $\PMC$ side terminates below $c_1$ or some
  strand on the $\PMC'$ side terminates above $c_2'$\footnote{More
    precisely: a strand terminates at some position $p'$ in $\PMC'$
    whose corresponding point $p\in\PMC$ is above~$c_2$.},
  we can factor off a
  near-chord of Type~(\ref{typeS:Generic}) on the right; similarly, if
  some strand on the $\PMC$ side starts above $c_2$ or some strand on the
  $\PMC'$ side starts below $c_1'$, we can factor off a near-chord in
  $\Dlg$ on the left. So, we can
  assume that no strand terminates below $c_1$ or above $c_2'$ or starts
  above $c_2$ or below $c_1'$.

  Next, for each case marked $V$ or $E$, there is a corresponding
  unmarked case; and the factorization for the $V$ or $E$ case is the
  same as the factorization in the unmarked case, except that a few
  cases cannot occur.
  So, it suffices to consider only the unmarked cases.

  We now analyze the remaining $12$ cases in turn.

  \textbf{Case $\lsub{duc}X\to \lsub{u}Y$.} 
  We have the following subcases:
  \begin{enumerate}
  \item \emph{If a strand starts at $u$,}  factor off a near-chord of
    Type~\eqref{typeS:I2} on the right.
  \item \emph{If a strand starts at $c_2$, but no strand starts
      at~$u$,} factor off an element of $\Dlg$ above~$c_2$ on the
    left.  (That is, factor off $a([c_2,p]) \otimes a'([c_2,p])$ for
    some $p$ above~$c_2$.)
  \item \emph{If no strand starts at $u$ or $c_2$,}  factor off a near-chord
    of Type~\eqref{type:l:I1I2+} on the left. (Note that, when
    $\epsilon=0$, there is, in fact, nothing left after we factor off
    the chord of Type~\eqref{type:l:I1I2+}. This happens again at
    various points later in the argument.)
  \end{enumerate}
  
  \textbf{Case $\lsub{duc}X\to \lsub{d}Y$.}  We have the following subcases:
  \begin{enumerate}
  \item \emph{If a strand starts at $c_2$,}
    factor off an element of $\Dlg$ above~$c_2$ on the left.
  \item \emph{If no strand starts at $c_2$,}
    factor off a near-chord of Type~\eqref{type:l:I1+} on the left.
  \end{enumerate}

  \textbf{Case $\lsub{du}X_c\to \lsub{u}Y$.}  We have the following
  subcases:
  \begin{enumerate}
  \item  \emph{If a strand starts at~$u$,}
    factor off a near-chord of Type~\eqref{typeS:I2} on the right.
  \item \emph{If no strand starts at~$u$,} factor off a near-chord of
    Type~\eqref{type:l:I1I2} on the left.
  \end{enumerate}

  \textbf{Case $\lsub{du}X_c\to \lsub{d}Y$.} Factor off a near-chord of Type~\eqref{type:l:I1} on the left.

  \textbf{Case $\lsub{dc}X_u\to \lsub{u}Y$.} Factor off a near-chord of
  Type~\eqref{type:l:I1I2} on the left.

  \textbf{Case $\lsub{dc}X_u\to \lsub{d}Y$.}
  Factor off a near-chord of Type~\eqref{type:l:I1I2} (in $\PMC'$),
  on the left.
  
  \textbf{Case $\lsub{d}X_{uc}\to \lsub{u}Y$.} Factor off a near-chord of Type~\eqref{typeS:I2} on the right.

  \textbf{Case $\lsub{d}X_{uc}\to \lsub{d}Y$.} We have the following subcases:
  \begin{enumerate}
  \item  \emph{If a strand terminates
      at $c_1'$,} factor off an element of $\Dlg$ below~$c_1$ on the
    left.
  \item \emph{If no strand terminates at $c_1'$,} factor off a near-chord of
    Type~\eqref{type:l:I1I2+} on the left.
  \end{enumerate}

  \textbf{Case $\lsub{uc}X_d \to \lsub{u}Y$.}
  Factor off a near-chord of Type~\eqref{type:l:I1} (in $\PMC'$), on the left.

  \textbf{Case $\lsub{uc}X_d \to \lsub{d}Y$.}
  Factor off a near-chord of Type~\eqref{type:l:I1} (in $\PMC'$), on the left.

  \textbf{Case $\lsub{u}X_{dc} \to \lsub{u}Y$.} We have the following subcases:
  \begin{enumerate}
  \item \emph{If a strand starts at $c_1'$,} factor off an element
    of~$\Dlg$ below~$c_1$ on the left.
  \item \emph{If no strand starts at $c_1'$,} factor off a
    near-chord of Type~\eqref{type:l:I1+} on the left.
  \end{enumerate}

  \textbf{Case $\lsub{u}X_{dc} \to \lsub{d}Y$.}
  Factor off a near-chord of Type~\eqref{type:l:I1} on the right.
\end{proof}

\begin{corollary}\label{cor:factor-M-minus-chords} 
  Every element $m$ of $\Mlg^-$ can be factored as a product
  $m=d_1\cdot\cdots d_m\cdot m'\cdot b_1\cdot\cdots\cdot b_n$ where
  the $d_i$ are near-chords in $\Dlg$, the $b_i$ are near-chords for
  $\Blg$, and $m'$ is a near-chord for $\Mlg^-$.
\end{corollary}
\begin{proof}
  This follows Proposition~\ref{prop:factor-M-minus} and the facts
  that elements of $\Blg$ and $\Dlg$ can be factored into
  near-chords---Lemma~\ref{lem:factor-diag}
and Proposition~\ref{prop:factor-anti-braid}, respectively. 
\end{proof}

\begin{corollary}\label{cor:Mlg-minus-Z-gr-unique}
  There is a unique relative $\ZZ$-grading on $\Mlg^-$ so that basic
  generators are homogeneous and for $d\in\Dlg$, $m\in\Mlg^-$ and
  $b\in\Blg$, we have $\gr(dmb)-\gr(m)=\gr(d)+\gr(b)$.
\end{corollary}
\begin{proof}
  This is immediate from Lemma~\ref{lem:M-minus-Z-graded} (existence)
  and Corollary~\ref{cor:factor-M-minus-chords} (uniqueness).
\end{proof}

\begin{proposition}
  \label{prop:CharacterizeF-}
  Let $\PMC$ be an unexceptional pointed matched circle (in the sense
  of Definition~\ref{def:special-length-3}).
  Let $F\in \Mlg^-$ be an element such that:
  \begin{itemize}
    \item $F$ contains all near chords of Type~(\ref{type:l:I1}).
    \item $F$ lies in a single grading with respect to the relative
      $\ZZ$-grading on $\Mlg^-$.
    \item $F$ solves the equation
      \begin{equation}
        \label{eq:MorphismEquation}
        dF+A_\Id\cdot F + F\cdot A_0=0.
      \end{equation}
    \end{itemize}
    Then $F$ is $F^-$, the sum of all near-chords in $\Mlg^-$.
\end{proposition}

\begin{proof}
  By a \emph{$\Dlg\Mlg^-\Blg$-factorization} of an element $x\in\Mlg^-$ we
  mean a nontrivial factorization $x=d\cdot m\cdot b$ where
  $d\in\Dlg$, $m\in\Mlg^-$ and $b\in\Blg$.

  Consider an element $x$ of Type~(\ref{type:l:I1}) and a chord $\xi$
  with initial endpoint $c_2$. The product $(a(\xi)\otimes a'(\xi))\cdot x$ occurs in
  $A_\Id\cdot F$, and has no other $\Dlg\Mlg^-\Blg$-factorization. Thus,
  by Equation~\eqref{eq:MorphismEquation}, $(a(\xi)\otimes
  a'(\xi))\cdot x$ must be the differential of some element $y$ in
  $F$. The only $y$ with $d(y)=(a(\xi)\otimes a'(\xi))\cdot x$ is a
  near-chord of Type~(\ref{type:l:I1+}); varying~$\xi$, this guarantees that all
  near-chords of Type~(\ref{type:l:I1+}) occur in $F$.

  Next, consider a near-chord $x$ of Type~(\ref{type:l:I1}),
  supported in $\PMC$, say. Let $\xi=[d,u]$; note that $a(\xi)\otimes
  1$ is a sum of near-chords of Type~(\ref{typeS:I2}). The product
  $x\cdot (a(\xi)\otimes 1)$ occurs in $F\cdot A$, and has no other
  $\Dlg\Mlg^-\Blg$-factorization. Thus, it must be canceled by a term
  of Type~(\ref{type:l:I1I2}) in $F$. This, and the corresponding
  argument on the $\PMC'$-side, force all terms of
  Type~(\ref{type:l:I1I2}) which are not cycles to occur in $F$.  The
  near-chords of Type~(\ref{type:l:I1I2}) which are cycles 
  are then forced to appear in $F$, as well, as we can see
  by commuting with $a(\eta) \otimes a'(\eta)$, where $\eta$ is a
  chord with one end point which
  is matched with~$u$, if the near-chord was supported in~$\PMC$,
  or~$d$, if the near-chord was supported in~$\PMC'$. Thus, all
  near-chords of Type~(\ref{type:l:I1I2}) occur in $F$.

  The existence of all chords of Type~(\ref{type:l:I1I2}) implies
  the existence of all chords of Type~(\ref{type:l:I1I2+}) by the same
  mechanism which established the existence of chords of
  Type~(\ref{type:l:I1+}) from the existence of chords of 
  Type~(\ref{type:l:I1}).

  Finally, it follows from Corollary~\ref{cor:factor-M-minus-chords}
  and the homogeneity of $F$ that $F$ contains no other
  terms: all other elements of $\Mlg^-$ must have strictly smaller grading.
\end{proof}

\begin{proposition}  
  \label{prop:CalculateNegMorphism}
  Let $\PMC$ be an unexceptional pointed matched circle and $\gamma$
  be a curve in
  $F(\PMC)$ corresponding to a pair of matched points in $\PMC$ whose
  feet are $3$ units apart.  Recall that there is an identification
  between generators of $\CFDDa(\Id)$ and basic idempotents in the
  diagonal subalgebra, and between generators of
  $\CFDDa(Y_{0(\gamma)})$ and basic idempotents in the anti-braid
  subalgebra. With respect to these identifications, the morphism
  $F^-_\gamma\co \CFDDa(\Id)\to \CFDDa(Y_{0(\gamma)})$ is given by
  \[
  F^-_\gamma(a\cdot I)=a\cdot I\cdot F^-.
  \]
  for any pair of complementary idempotents $I$ and element
  $a\in\Alg(\PMC)\otimes\Alg(\PMC')$.
\end{proposition}

\begin{proof}
  From the forms of the modules, it is immediate that $F^-_\gamma$ is
  given by 
  \[
  F^-_\gamma(a\cdot I)=a\cdot I\cdot F
  \]
  for some $F\in\Alg(\PMC)\otimes\Alg(\PMC')$. It follows from the
  fact that $F^-_\gamma$ is a homomorphism of type \DD\ bimodules that
  $F$ satisfies Equation~\eqref{eq:MorphismEquation}. The map
  $F_\gamma^-$ is triangle-like with respect to the Heegaard triple
  diagram shown in Figure~\ref{fig:TripleDiagram}. So, it follows from
  inspecting the diagram that $F\in \Mlg^-$, and from the definition
  of $\ZZ$-grading (from Lemma~\ref{lem:M-minus-Z-graded}) that $F$
  lies in a single grading.

  \begin{figure}
    \centering
    \input{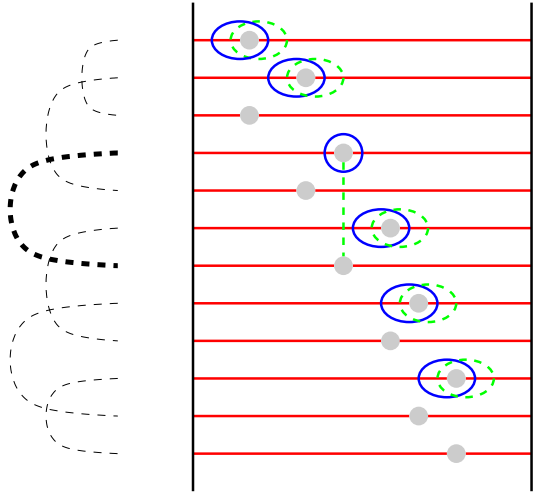}    
    \caption{\textbf{A Heegaard triple diagram for $F^-_\gamma$.}
    This is gotten by overlaying the two diagrams from Figure~\ref{fig:HeegaardDiagrams}.
    The ordering of the circles in the triple determines whether we are considering
    $F^-_\gamma$ or $F^+_\gamma$.
    \label{fig:TripleDiagram}}
  \end{figure}

  To see that $F$ contains all near chords of
  Type~(\ref{type:l:I1}), consider the effect of setting to $0$ any
  element of $\Alg(\PMC)\otimes\Alg(\PMC')$ whose support is not
  entirely contained in $\Zint_3\subset\PMC$. The differential on
  $\CFDDa(\Dehn_\gamma)$ contains the terms $a(\Zint_3)\otimes 1$, and
  if $F$ does not contain these terms then the mapping cone of
  $F_\gamma^-$ is not homotopy equivalent to
  $\CFDDa(\Dehn_\gamma)$. Applying the same argument to
  $\Zpint_1\subset\PMC'$ implies that $F$ also contains the terms
  $1\otimes a'(\Zpint_1)$.

  Thus, Proposition~\ref{prop:CharacterizeF-} implies that $F=F^-$.
\end{proof}

\begin{remark}
  Proposition~\ref{prop:CalculateNegMorphism} also implies that $F^-$
  does solve Equation~\eqref{eq:MorphismEquation}. (Of
  course, this can be checked directly, but it is tedious to do so.)
\end{remark}

\subsection{Right-handed cases}
\label{subsec:RightHanded}

Next, we turn to the map $F^+\co \CFDDa(Y_{0(\gamma)})\to \CFDDa(\Id)$.

\begin{definition}\label{def:MlgPlus}
  The {\em right-handed morphism sub-algebra} $\Mlg^+$ is a subalgebra
  of $\Alg(\PMC)\otimes \Alg(\PMC')$ which is spanned by algebra
  elements $a\otimes a'$ satisfying the following properties:
  \begin{itemize}
  \item $I\cdot a \cdot J=a$ and $I'\cdot a'\cdot J'=a'$ where
    $I\otimes I'$ is an idempotent for the anti-braid subalgebra,
    while $J\otimes J'$ is an idempotent for the diagonal subalgebra.
  \item the local multiplicities of $a$ outside the interval
    $[c_1,c_2]$ agree with the corresponding local multiplicities of
    $a'$ outside the interval $[c_1',c_2']$.
  \end{itemize}
\end{definition}
\begin{lemma}
  $\Blg\cdot \Mlg^+ \cdot \Dlg =\Mlg^+$, where here $\Blg$ is the
  anti-braid subalgebra of Definition~\ref{def:AntiBraidSubalgebra},
  $\Dlg$ is the diagonal subalgebra and $\Mlg^+$
  is the right-handed morphism subalgebra.
\end{lemma}
\begin{proof}
  This is immediate from the definitions.
\end{proof}

\begin{lemma}\label{lem:plus-is-op-of-minus}
  Fix a pointed matched circle $\PMC$ and a matched pair $\{c_1,c_2\}$
  in $\PMC$ so that the distance between $c_1$ and $c_2$ is $3$. Then
  the right-handed morphism subalgebra $\Mlg^+(\PMC,\{c_1,c_2\})$ of
  $\PMC$ with respect to $\{c_1,c_2\}$ is the opposite algebra of the
  left-handed morphism subalgebra of $-\PMC$ with respect to
  $\{c_1,c_2\}$. This identification intertwines the actions of
  $\Blg(\PMC)$ and $\Dlg(\PMC)$ with the actions of $\Blg(-\PMC)$ and
  $\Dlg(-\PMC)$.
\end{lemma}
\begin{proof}
  This is clear.
\end{proof}

\begin{lemma}\label{lem:M-plus-Z-graded}
  There is a $\ZZ$-grading on $\Mlg^+$ which is compatible with the
  $\ZZ$-gradings on $\Blg$ and $\Dlg$, in the sense that
  $\gr(d\cdot m\cdot b)=\gr(d)+\gr(m)+\gr(b)$ for $d\in \Dlg$, $m\in
  \Mlg^+$ and $b\in\Blg$.
\end{lemma}
\begin{proof}
  The proof is the same as the proof of
  Lemma~\ref{lem:M-minus-Z-graded}. Alternatively, this follows from
  Lemmas~\ref{lem:M-minus-Z-graded} and~\ref{lem:plus-is-op-of-minus}.
\end{proof}

\begin{figure}
  \centering
  \includegraphics{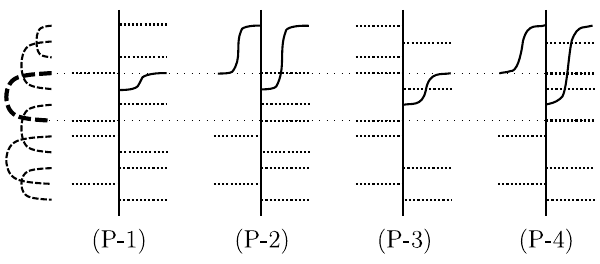}
  \caption{\textbf{Right-handed near-chords.} We have drawn the second
    case of each type of near chord from Definition~\ref{def:rhnc}.}
  \label{fig:rhnc}
\end{figure}

\begin{definition}\label{def:rhnc}
  An element $x\in \Mlg^+$ is called a \emph{(right-handed)
    near-chord} if it has one
  of the following forms:
  \begin{enumerate}[label=(P-\arabic*),ref=P-\arabic*]
  \item 
    \label{type:r:I1}
    $x=I\cdot (a(\Zint_1)\otimes 1)\cdot J$ or $x=I\cdot ( 1\otimes
    a'(\Zpint_3))\cdot J$. 
  \item 
    \label{type:r:I1+}
    $x=I\cdot (a(\Zint_1\cup\xi)\otimes a'(\xi))\cdot J$, where $\xi$
    is a chord ending at $c_1$; or $x=I\cdot (a(\xi)\otimes
    a'(\Zpint_3\cup\xi))\cdot J$, where $\xi$ is a chord starting at
    $c_2$.
  \item 
    \label{type:r:I1I2}
    $x=I\cdot (a(\Zint_1\cup\Zint_2)\otimes 1)\cdot J$ or $x=I\cdot (
    1\otimes a'(\Zpint_2\cup\Zpint_3))\cdot J$. 
  \item 
    \label{type:r:I1I2+}
    $x=I\cdot (a(\Zint_1\cup\Zint_2\cup \xi)\otimes a'(\xi))\cdot J$
    where $\xi$ is a chord ending at $c_1$, or $x=I\cdot
    (a(\xi)\otimes a'(\Zpint_2\cup\Zpint_3\cup\xi))\cdot J$ where
    $\xi$ is a chord starting at $c_2$. 
  \end{enumerate}
  See Figure~\ref{fig:rhnc}.

  The {\em structure constant for the right-handed morphism} 
  $F^+$ is the sum of all (right-handed) near-chords in $\Mlg^+$.
\end{definition}

\begin{proposition}\label{prop:factor-M-plus} 
  Every element $m$ of $\Mlg^+$ can be factored as a product $m=b\cdot
  m'\cdot d$ where $b\in\Blg$, $d\in\Dlg$, and $m'$ is a near-chord for $\Mlg^+$.
\end{proposition}
\begin{proof}
  This follows from Proposition~\ref{prop:factor-M-minus} and
  Lemma~\ref{lem:plus-is-op-of-minus}. (Alternatively, the proof is a
  straightforward adaptation of the proof of
  Proposition~\ref{prop:factor-M-minus}.)
\end{proof}

\begin{corollary}\label{cor:factor-M-plus-chords} 
  Every element $m$ of $\Mlg^+$ can be factored as a product
  $m=b_1\cdot\cdots\cdot b_n\cdot m'\cdot d_1\cdot\cdots d_m$ where
  the $b_i$ are near-chords for $\Blg$, the $d_i$ are near-chords in
  $\Dlg$, and $m'$ is a near-chord for $\Mlg^-$.
\end{corollary}
\begin{proof}
  This is immediate from Proposition~\ref{prop:factor-M-plus},
  Lemma~\ref{lem:factor-diag} and
  Proposition~\ref{prop:factor-anti-braid}.
\end{proof}

\begin{corollary}\label{cor:Mlg-plus-Z-gr-unique}
  There is a unique relative $\ZZ$-grading on $\Mlg^+$ so that basic
  generators are homogeneous and for $d\in\Dlg$, $m\in\Mlg^+$ and
  $b\in\Blg$, $\gr(dmb)-\gr(m)=\gr(d)+\gr(b)$.
\end{corollary}
\begin{proof}
  This is immediate from Lemma~\ref{lem:M-plus-Z-graded} (existence)
  and Corollary~\ref{cor:factor-M-plus-chords} (uniqueness).
\end{proof}

\begin{proposition}
  \label{prop:CharacterizeF+}
  Let $\PMC$ be an unexceptional pointed matched circle (in the sense
  of Definition~\ref{def:special-length-3}).
  Let $F\in \Mlg^+$ be an element such that:
  \begin{itemize}
    \item $F$ contains all near chords of Type~(\ref{type:r:I1}).
    \item $F$ lies in a single grading with respect to the relative
      $\ZZ$-grading on $\Mlg^+$.
    \item $F$ solves the equation
      \begin{equation}
        \label{eq:MorphismEquationPlus}
        dF+A_0\cdot F + F\cdot A_\Id=0.
      \end{equation}
    \end{itemize}
    Then $F=F^+$.
\end{proposition}
\begin{proof}
  By Lemma~\ref{lem:plus-is-op-of-minus}, $F$ corresponds to an
  element $F'$ of $\Mlg^-(-\PMC)$ which satisfies the conditions of
  Proposition~\ref{prop:CharacterizeF-}. Thus, $F'=F^-(-\PMC)$, and so
  $F=F^+(\PMC)$, as desired. (Alternatively, one could imitate the
  proof of Proposition~\ref{prop:CharacterizeF-}, switching the
  orientations and the orders of all the products.)
\end{proof}

\begin{proposition}  
  \label{prop:CalculatePosMorphism}
  Let $\PMC$ be a non-exceptional pointed matched circle and $\gamma$
  be a curve in
  $F(\PMC)$ corresponding to a pair of matched points in $\PMC$ whose
  feet are $3$ units apart.  Recall that there is an identification
  between generators of $\CFDDa(\Id)$ and basic idempotents in the
  diagonal subalgebra, and between generators of
  $\CFDDa(Y_{0(\gamma)})$ and basic idempotents in the anti-braid
  subalgebra. With respect to these identifications, the morphism
  $F^+_\gamma\co \CFDDa(Y_{0(\gamma)})\to \CFDDa(\Id)$ is given by
  \[
  F^+_\gamma(a\cdot I)=a\cdot I\cdot F^+
  \]
  for any pair of complementary idempotents $I$ and element
  $a\in\Alg(\PMC)\otimes\Alg(\PMC')$.
\end{proposition}
\begin{proof}
  The proof is essentially the same as the proof of
  Proposition~\ref{prop:CalculateNegMorphism}.
\end{proof}

\subsection{Degenerate left- and right-handed morphisms}
\label{subsec:DegenerateEdge}

For computations it is useful to have one more, degenerate case: when
the points $c_1$ and $c_2$ are separated by a single point $p$, but
lie at an end of the pointed matched circle
(Hypothesis~\ref{hyp:DistanceOne}). 

The definitions of the negative
and positive morphism subalgebras for this case are analogous to those
in the non-degenerate case (Definitions~\ref{def:MlgMinus}
and~\ref{def:MlgPlus} respectively), except that now we use the
degenerate anti-braid subalgebra from
Definition~\ref{def:AntiBraidSubalgebra} in place of the anti-braid
subalgebra. We denote these two morphism subalgebras
$\Mlg^\pm_{\degen}$.

As in Subsection~\ref{subsec:DegenerateVertex} (see especially
Hypothesis~\ref{hyp:DistanceOne} and
Remark~\ref{rmk:ChangeHypothesis}), we focus here on the case where
the basepoint $z$ lies just below $c_1$.  The other extreme, where $z$
lies just above $c_2$, follows similarly, with routine notational
changes.

Now, we have the following types of near-chords in the left-handed
case:
\begin{definition}   
  In the degenerate case (Hypothesis~\ref{hyp:DistanceOne}),
  an element $x\in \Mlg^-_{\degen}$ is called a
  \emph{left-handed degenerate near-chord} if it has the following form:
  \begin{enumerate}[label=(N$_\degen$-\arabic*),ref=N$_\degen$-\arabic*]
  \item 
    \label{type:degen:l:I1}
    $x=I\cdot (a(\Zint_2)\otimes 1)\cdot J$ or $x=I\cdot (1\otimes
    a(\Zpint_1))\cdot J$. 
  \item     
    \label{type:degen:l:I1+}
    $x=I\cdot (a(\Zint_2\cup\xi)\otimes a'(\xi))\cdot J$, where $\xi=[c_2,e]$
    for some $e$ above $c_2$.
  \end{enumerate}
  The {\em structure constant for the degenerate left-handed morphism} 
  $F^-_{\degen}\in \Mlg^-_{\degen}$ is the sum of all left-handed
  degenerate near-chords.
\end{definition}

Similarly we have the following definition in the right-handed case:

\begin{definition}\label{def:degen:rhnc}
  An element $x\in \Mlg^+_{\degen}$ is called a \emph{right-handed degenerate
    near-chord} if it has the following form:
  \begin{enumerate}[label=(P$_\degen$-\arabic*),ref=P$_\degen$-\arabic*]
  \item 
    \label{type:degen:r:I1}
    $x=I\cdot (a(\Zint_1)\otimes 1)\cdot J$ or $x=I\cdot ( 1\otimes
    a'(\Zpint_2))\cdot J$. 
  \item 
    \label{type:degen:r:I1+}
    $x=I\cdot (a(\xi)\otimes
    a'(\Zpint_2\cup\xi))\cdot J$, where $\xi=[c_2,e]$ for some $e$ above $c_2$.
  \end{enumerate}
  The {\em structure constant for the degenerate right-handed morphism} 
  $F^+_{\degen}\in \Mlg^+_{\degen}$ is the sum of all right-handed
  degenerate near-chords.
\end{definition}

With these definitions in place, we have the following analogues of
Proposition~\ref{prop:CalculateNegMorphism} and~\ref{prop:CalculatePosMorphism}:

\begin{proposition}  
  \label{prop:CalculateNegMorphismDeg}
  Let $\PMC$ be a non-exceptional pointed matched circle and $\gamma$
  be a curve in
  $F(\PMC)$ satisfying Hypothesis~\ref{hyp:DistanceTwo}.  Recall that there is an identification
  between generators of $\CFDDa(\Id)$ and basic idempotents in the
  diagonal subalgebra, and between generators of
  $\CFDDa(Y_{0(\gamma)})$ and basic idempotents in the anti-braid
  subalgebra. With respect to these identifications, the morphism
  $F^-_\gamma\co \CFDDa(\Id)\to \CFDDa(Y_{0(\gamma)})$ is given by
  \[
  F^-_\gamma(a\cdot I)=a\cdot I\cdot F^-_{\degen}
  \]
  for any pair of complementary idempotents $I$ and element
  $a\in\Alg(\PMC)\otimes\Alg(\PMC')$.
\end{proposition}

\begin{proof}
  Adapting the proof of Proposition~\ref{prop:CalculateNegMorphism},
  this follows from a degenerate analogue of
  Proposition~\ref{prop:CharacterizeF-}. 
\end{proof}

\begin{proposition}  
  \label{prop:CalculatePosMorphismDeg}
  Let $\PMC$ be a non-exceptional pointed matched circle and $\gamma$
  be a curve in
 $F(\PMC)$ satisfying Hypothesis~\ref{hyp:DistanceTwo}.
 Recall that there is an identification
  between generators of $\CFDDa(\Id)$ and basic idempotents in the
  diagonal subalgebra, and between generators of
  $\CFDDa(Y_{0(\gamma)})$ and basic idempotents in the anti-braid
  subalgebra. With respect to these identifications, the morphism
  $F^+_\gamma\co \CFDDa(Y_{0(\gamma)})\to \CFDDa(\Id)$ is given by
  \[
  F^+_\gamma(a\cdot I)=a\cdot I\cdot F^+_{\degen}
  \]
  for any pair of complementary idempotents $I$ and element
  $a\in\Alg(\PMC)\otimes\Alg(\PMC')$.
\end{proposition}

\begin{proof}
  This proceeds exactly as Proposition~\ref{prop:CalculateNegMorphismDeg}.
\end{proof}

\begin{example}\label{exam:pos-morph}
  Consider the positive morphism in the genus $1$ case from
  Example~\ref{exam:abl-torus} (i.e., $k=1$, $n=4$). In this case,
  there are two right-handed degenerate near-chords of
  type~(\ref{type:degen:r:I1}): $\sigma_2\otimes \iota_0$ and
  $\iota_1\otimes \rho_1$ . There are no right-handed degenerate
  near-chords of type~(\ref{type:degen:r:I1+}). Thus,
  $F^+_{\degen}=\sigma_2\otimes \iota_0+\iota_1\otimes \rho_1$, and
  the map $F^+_\gamma$ is given by
  \[
  F^+_\gamma(I)=(\sigma_2\otimes 1)J+(1\otimes\rho_1)K.
  \]
  (Similarly to Example~\ref{exam:abl-torus}, since $\PMC$ is
  exceptional a little more work is needed to verify this claim.)
\end{example}

\subsection{Consequences}
\label{subsec:Consequences}

We can now combine the above results to deduce Theorem~\ref{thm:MappingCone}.

\begin{proof}[Proof of Theorem~\ref{thm:MappingCone}]
  Start from Theorem~\ref{thm:Dehn-is-MC}. In the negative
  case, we apply   Proposition~\ref{prop:CalculateNegMorphism}
  or~\ref{prop:CalculateNegMorphismDeg}  (the latter in the degenerate case);
  in the positive case, we apply 
  Proposition~\ref{prop:CalculatePosMorphism}
  or~\ref{prop:CalculatePosMorphismDeg}  (the latter in the degenerate case).
\end{proof}


\section{Putting it all together: computing the spectral sequence}\label{sec:cube}
\begin{figure}
  \centering
  \includegraphics{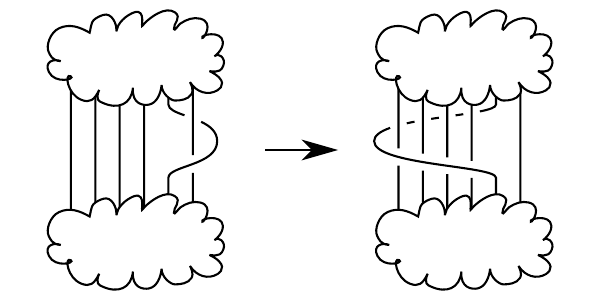}
  \caption{\textbf{Finding a braid representing $L$ and not using
      $s_{n-1}$.} The clouds represent the rest of the braid and the
    plats closing it.}
  \label{fig:clouds}
\end{figure}

Let $L$ be a link. Present $L$ as the plat closure of a braid $B$. We
may arrange that, in~$B$, the last strand is not involved in any
crossings; that is, $s_{n-1}$ does not occur in $B$.
(To do this, permute the bottom so that the last strand gets mapped to
itself, then pull the last strand tight so we only see $(s_{n-1})^2$,
and finally replace $(s_{n-1})^2$ by an equivalent sequence involving
the other generators; see Figure~\ref{fig:clouds}.) The branched
double cover of $B$ is the mapping cylinder of a word in Dehn twists along curves
$\gamma_1,\dots,\gamma_{n-2}$ in a surface $S$ of genus $n/2-1$. These
curves form a linear chain in $S$, in the sense that $\gamma_i$
intersects $\gamma_{i+1}$ in a single point for each $i$; thus, it is
convenient to view these as the standard curves in the surface
associated to the linear pointed matched circle from
Section~\ref{sec:diagrams}.

Now, fix a link $L$. To compute the spectral sequence from Khovanov
homology to $\HFa(\Sigma(L))$, one proceeds in six steps:
\begin{enumerate}
\item Present $L$ as the plat closure of a braid $B$, and slice $B$
  into a sequence of braid generators $B_1,\dots, B_n$. Let $B_0$ and
  $B_{n+1}$ denote the crossingless matchings used to plat-close the
  braid, so $L=B_0\cup B_1\cup\cdots\cup B_n\cup B_{n+1}$.
\item For each $B_i$, let $\Sigma(B_i)$ denote the branched double
  cover of $B_i$ bordered by $F(\PMC)$ where $\PMC$ is the linear
  pointed matched circle. By Theorem~\ref{thm:Dehn-is-MC},
  $\CFDDa(\Sigma(B_i))=\Cone(F_\gamma^+)$ if $B_i$ is a positive braid
  generator (with respect to our conventions, from
  Figure~\ref{fig:resolve-braid-gen}) and 
  $\CFDDa(\Sigma(B_i))=\Cone(F_\gamma^-)$ if $B_i$ is a negative braid.
  Let $\widecheck{B}_i$ denote the anti-braid resolution of $B_i$.

  Compute $\CFDDa(\Id)$ via Proposition~\ref{prop:RecallIdentity} and
  $\CFDDa(\Sigma(\widecheck{B}_i))$ via
  Proposition~\ref{prop:CalculateVertex}. When $\gamma$ is not the
  first or last
  braid generator, compute the map $F_\gamma^\pm$ via
  Proposition~\ref{prop:CalculateNegMorphism}
  or~\ref{prop:CalculatePosMorphism}; when $\gamma$ is 
  the first or last braid generator, compute $F_\gamma^\pm$
  via 
  Proposition~\ref{prop:CalculateNegMorphismDeg}
  or~\ref{prop:CalculatePosMorphismDeg}.
\item Compute $\CFAAa(\Id)$,
  using~\cite[Theorem~\ref*{Bimodules:prop:DDAA-duality}]{LOT2}, which
  shows that
  $\CFAAa(\Id)=\Mor_{\Alg(\PMC)}(\CFDDa(\Id),\Alg(\PMC))$. As
  explained in \cite{LOT2},
  one should view $\CFDDa(\Id)$ as a bimodule, not a type \DD\
  structure, when taking morphisms. (The difference is whether a copy
  of the algebra is dualized; the other way, one picks up some
  boundary Dehn twists, though they will disappear in the tensor
  product~\eqref{eq:final-cone}.)

  In practice, one might want to simplify this complex before
  proceeding. See \cite[Section~\ref*{HFa:sec:hpt}]{LOT4} for some relevant
  discussion.
\item\label{step:tensorAA} Turn the type \DD\ modules and maps into type \DA\ modules and
  maps by tensoring with $\CFAAa(\Id)$. 
  That is, if $B_i$ is a negative braid generator, 
  \[
  \CFDAa(B_i)=
  \Cone(\Id\DT F_\gamma^-\co \CFAAa(\Id)\DT\CFDDa(\Id)\to
  \CFAAa(\Id)\DT \CFDDa(Y_{0(\gamma)})); \]
  while if $B_i$ is a positive braid generator,
  \[
  \CFDAa(B_i)=
  \Cone(\Id\DT F_\gamma^+\co \CFAAa(\Id)\DT
  \CFDDa(Y_{0(\gamma)}) \to \CFAAa(\Id)\DT\CFDDa(\Id)).
  \]
\item Compute the modules $\CFDa(B_0)$ and $\CFDa(B_{n+1})$ using
  Proposition~\ref{prop:comp-plat-handlebody}. Tensor $\CFDa(B_0)$
  with $\CFAAa(\Id)$ to obtain $\CFAa(B_0)$.
\item Let $F_{B_i}$ be the map $\Id\DT F_\gamma^\pm$ associated to
  piece $B_i$ in step~(\ref{step:tensorAA}).
  Compute the tensor product
  \begin{equation}
    \CFAa(B_0)\DT \Cone(F_{B_1})\DT\cdots\DT\Cone(F_{B_n})\DT\CFDa(B_{n+1}).\label{eq:final-cone}  
  \end{equation}
  By Lemma~\ref{lem:FilteredProduct}, this is a $\{0,1\}^n$-filtered
  chain complex. By Theorem~\ref{thm:BorderedSpectralSequence}, this
  spectral sequence has $E_2$-page equal to $\rKh(r(L))$ and
  $E_\infty$-page equal to $\HFa(\Sigma(L))$.
\end{enumerate}

\subsection{An example}
\begin{figure}
  \centering
  \includegraphics{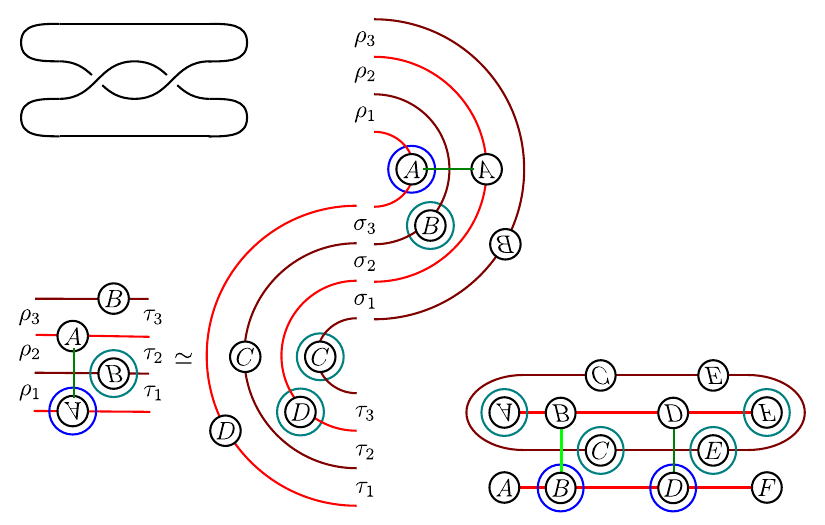}
  \caption{\textbf{Diagrams for the Hopf link.} Top left: the Hopf link as
    braid closure. Center: The first step in the computation is
    turning the filtered type \DD\ bimodule corresponding to the
    displayed multi-diagram into a filtered type \DA\ bimodule by
    tensoring with the type \AAm\ identity module. Lower right: A Heegaard
    multi-diagram for the cube of resolutions of the Hopf link. (Note
    that repeated curves in the multi-diagrams are only drawn once.)}
  \label{fig:hopf}
\end{figure}

We will use the techniques of this paper to compute the spectral
sequence when $L$ is the Hopf link, which is the plat closure of a
$4$-strand braid with braid word
$s_2^2$; see Figure~\ref{fig:hopf}. Since the double
cover of a sphere branched over $4$ points is a torus, this
computation takes place in $\Alg=\Alg(T^2)$. The algebra $\Alg(T^2)$ is
$8$-dimensional, and was described in Example~\ref{exam:DD-id-torus}.

Let $\gamma$ be the curve in the torus corresponding to the braid
generator $s_2$, and let \[F^+_\gamma\co \CFDDa(Y_{0(\gamma)})\to
\CFDDa(\Id)\] be the map from Theorem~\ref{thm:MappingCone}.  We start by
computing the filtered type \DD\ bimodule $\Cone(F^+_\gamma)$.

Recall from Example~\ref{exam:DD-id-torus} that
$\CFDDa(\Id)$ has two generators $J$ and $K$ over $\Alg^\sigma$ and
$\Alg^\rho$, and differentials
\begin{align*}
  \nabla_0 J &=
  (\sigma_1\otimes\rho_3+\sigma_3\otimes\rho_1+\sigma_{123}\otimes\rho_{123})K\\
  \nabla_0 K&= (\sigma_2\otimes\rho_2)J.
\end{align*}
Similarly, recall from Example~\ref{exam:abl-torus} that
$\CFDDa(Y_{0(\gamma)})$ has one generator $I$ over $\Alg^\sigma$ and
$\Alg^\rho$, and differential
\[
\nabla_0 I=(\sigma_{23}\otimes 1+1\otimes\rho_{12}) I.
\]
Finally, recall from Example~\ref{exam:pos-morph} that the map
$F^+_\gamma$ is given by
\[
F^+_\gamma(I)=(\sigma_2\otimes 1)J+(1\otimes\rho_1)K.
\]

So, $\Cone(F^+_\gamma)$ is given by
\[
\begin{tikzpicture}
  \node at (0,0) (I) {$I^0$};
  \node at (4,1) (J) {$J^1$};
  \node at (4,-1) (K) {$K^1$};
  \draw[->] (I) [loop left] to
  node[left]{\lab{\sigma_{23}\otimes 1 + 1\otimes \rho_{12}}}
  (I);
  \draw[->] (I) to node[above, sloped]{\lab{\sigma_2\otimes 1}} (J);
  \draw[->] (I) to node[below, sloped]{\lab{1\otimes \rho_1}} (K);
  \draw[->, bend left=20] (K) to node[left]{\lab{\sigma_2\otimes \rho_2}} (J);
  \draw[->, bend left=20] (J) to node[right]{\lab{\sigma_1\otimes \rho_3+\sigma_3\otimes \rho_1+\sigma_{123}\otimes \rho_{123}.}} (K);
\end{tikzpicture}
\]
where the filtration is indicated by the exponent.

Recall from~\cite[Section~\ref*{Bimodules:subsec:AAId1}]{LOT2} that, in the torus case, $\CFAAa(\Id)$ is
given by:
\[
  \begin{tikzpicture}[y=54pt,x=1in]
    \node at (0,3) (w1) {$W_1$} ;
    \node at (2,3) (z1) {$Z_1$} ;
    \node at (1,2) (y)  {$Y$} ;
    \node at (1,1) (x)  {$X$} ;
    \node at (0,0) (w2) {$W_2$} ;
    \node at (2,0) (z2) {$Z_2$} ;
    \draw[->] (w1) to node[above,sloped] {\lab{\sigma_1}} (y) ;
    \draw[->] (z1) to node[above] {\lab{\tau_1}} (y) ;
    \draw[->] (y)  to node[above,sloped] {\lab{(\sigma_2,\tau_2)}} (x) ;
    \draw[->] (x)  to node[below,sloped] {\lab{\tau_3}} (w2) ;
    \draw[->] (x)  to node[below] {\lab{\sigma_3}} (z2) ;
    \draw[->] (w1) to node[below,sloped] {\lab{1+(\sigma_{12},\tau_{23})}} (w2) ;
    \draw[->] (z1) to node[above,sloped] {\lab{1+(\sigma_{23},\tau_{12})}} (z2) ;
    \draw[->] (w1) to[pos=0.4] node[below,sloped] {\lab{(\sigma_{12},\tau_2)}} (x) ;
    \draw[->] (y)  to[pos=0.6] node[above,sloped] {\lab{(\sigma_2,\tau_{23})}} (w2) ;
    \draw[->] (z1) to[pos=0.4] node[below,sloped] {\lab{(\sigma_2,\tau_{12})}} (x) ;
    \draw[->] (y)  to[pos=0.6] node[above,sloped] {\lab{(\sigma_{23},\tau_2)}} (z2) ;
    \draw[->] (w1) to[out=-125,in=145] (-0.25,-0.25) to[out=-35,in=-150] node[pos=0,below,sloped]
       {\lab{(\sigma_{123},\tau_2)+(\sigma_3,\sigma_2,\sigma_1,\tau_2)}} (z2) ;
    \draw[->] (z1) to[out=-55,in=35] (2.25,-0.25) to[out=-145,in=-30] node[pos=0,below,sloped]
       {\lab{(\sigma_{2},\tau_{123})}} (w2) ;
  \end{tikzpicture}
\]

Tensoring $\Cone(F^+)$ with $\CFAAa(\Id)$ over $\Alg^\sigma$ gives
\[
\begin{tikzpicture}
\node at (-1,0) (Z2K) {$Z_2K^1$};
\node at (3,0) (Z1K) {$Z_1K^1$};
\node at (8,0) (YK) {$YK^1$};
\node at (7,-2) (W1J) {$W_1J^1$};
\node at (4,-2) (W2J) {$W_2J^1$};
\node at (1,-2) (XJ) {$XJ^1$};
\node at (-1,-6) (Z2I) {$Z_2I^0$};
\node at (3,-5) (Z1I) {$Z_1I^0$};
\node at (8,-6) (YI) {$YI^0$};
\draw[->] (Z1K) to node[above, sloped]{\lab{1}} (Z2K);
\draw[->] (Z1K) to node[above, sloped]{\lab{\tau_1}} (YK);
\draw[->] (Z1I) to node[above, sloped]{\lab{1+\tau_{12}}} (Z2I);
\draw[->] (Z1I) to node[above, sloped]{\lab{\tau_{1}}} (YI);
\draw[->, bend left=10] (YI) to node[below, sloped]{\lab{\tau_{2}}} (Z2I);
\draw[->] (W1J) to node[above, sloped]{\lab{1}} (W2J);
\draw[->] (XJ) to node[above, sloped]{\lab{\tau_3}} (W2J);
\draw[->] (W1J) to node[above, sloped]{\lab{\rho_3}} (YK);
\draw[->] (XJ) to node[below, sloped]{\lab{\rho_1}} (Z2K);
\draw[->] (YI) to node[above, sloped]{\lab{\tau_{23}}} (W2J);
\draw[->] (Z1I) to node[above, sloped]{\lab{\tau_{12}}} (XJ);
\draw[->] (YK) to node[above, sloped]{\lab{\rho_2\tau_{2}}} (XJ);
\draw[->] (YK) to node[above, sloped]{\lab{\rho_2\tau_{23}}} (W2J);
\draw[->] (Z1K) to node[below, sloped]{\lab{\rho_2\tau_{123}}} (W2J);
\draw[->] (Z1K) to node[above, sloped]{\lab{\rho_2\tau_{12}}} (XJ);
\draw[->] (Z2I) to node[left]{\lab{\rho_1}} (Z2K);
\draw[->] (Z1I) to node[left]{\lab{\rho_1}} (Z1K);
\draw[->] (YI) to node[right]{\lab{\rho_1}} (YK);
\draw[->] (Z1I) to node[below, sloped]{\lab{\tau_{123}}} (W2J);
\draw[->] (YI) to node[above, sloped]{\lab{\tau_{2}}} (XJ);
\draw[->] (Z2I) [loop below] to node[below]{\lab{\rho_{12}}} (Z2I);
\draw[->] (Z1I) [loop below] to node[below]{\lab{\rho_{12}}} (Z1I);
\draw[->] (YI) to [loop below] node[below]{\lab{\rho_{12}}} (YI);
\end{tikzpicture}
\]
where again the filtration level is indicated by the exponent. (The
fact that we are tensoring over the idempotents means certain pairs of
generators do not appear.)
Here, for instance, the arrow from $YK$ to $W_2J$ labeled by
$\rho_2\tau_{23}$ means that $\delta^1_2(YK,\tau_{23})=\rho_{2}W_2J$.

Cancelling the arrows from $Z_1K$ to $Z_2K$ and from $W_1J$ to $W_2J$
gives us the somewhat simpler module:
\[
\begin{tikzpicture}
\node at (8,0) (YK) {$YK^1$};
\node at (0,0) (XJ) {$XJ^1$};
\node at (0,-4) (Z2I) {$Z_2I^0$};
\node at (4,-3) (Z1I) {$Z_1I^0$};
\node at (8,-4) (YI) {$YI^0$};
\draw[->, bend left=8] (Z2I) to node[above, sloped]{\lab{\rho_{1}\tau_{1}+\rho_{123}\tau_{123}}} (YK);%
\draw[->] (Z2I) to node[left]{\lab{\rho_{12}\tau_{12}}} (XJ);%
\draw[->] (Z1I) to node[above, sloped]{\lab{\rho_3\tau_{123}}} (YK);%
\draw[->] (Z1I) to node[above, sloped]{\lab{1+\tau_{12}}} (Z2I);%
\draw[->] (Z1I) to node[above, sloped]{\lab{\tau_{1}}} (YI);%
\draw[->, bend left=10] (YI) to node[below, sloped]{\lab{\tau_{2}}} (Z2I);%
\draw[->] (Z1I) to node[above, sloped]{\lab{\tau_{12}}} (XJ);%
\draw[->, bend left=15] (YK) to node[below, sloped]{\lab{\rho_2\tau_{2}}} (XJ);%
\draw[->, bend left=15] (XJ) to node[above, sloped]{\lab{\rho_1\tau_{1}+\rho_3\tau_3+\rho_{123}\tau_{123}}} (YK);%
\draw[->] (YI) to node[right]{\lab{\rho_1+\rho_{3}\tau_{23}}} (YK);%
\draw[->] (YI) to node[above, sloped]{\lab{\tau_{2}}} (XJ);%
\draw[->] (Z2I) [loop below] to node[below]{\lab{\rho_{12}}} (Z2I);%
\draw[->] (Z1I) [loop below] to node[below]{\lab{\rho_{12}}} (Z1I);%
\draw[->] (YI) [loop below] to node[below]{\lab{\rho_{12}}} (YI);%
\draw[->] (YK) [loop above] to node[above]{\lab{\rho_{23}\tau_{23}}} (YK);%
\draw[->] (XJ) [loop above] to node[above]{\lab{\rho_{12}\tau_{12}}} (XJ);%
\end{tikzpicture}
\]
We will take this as our
type \DA\ model for $\Cone(F^+_\gamma)$.

The type $D$ module for the plat solid torus has one generator $p$ and
\[
\bdy (p) =\tau_{23}p.
\]

Tensoring this with our type \DA\ model for  $\Cone(F^+_\gamma)$ gives 
\[
\begin{tikzpicture}
  \node at (0,0) (YIp) {$YIp^0$};
  \node at (4,0) (YKp) {$YKp^0$};
  \draw[->] (YIp) to node[above]{\lab{\rho_1+\rho_3}} (YKp);
  \draw[->] (YIp) [loop left] to node[left]{\lab{\rho_{12}}} (YIp);
  \draw[->] (YKp) [loop right] to node[right]{\lab{\rho_{23}}} (YKp);
\end{tikzpicture}
\]
Call
this module $M^\rho$. Let $M^\tau$ be the corresponding type $D$
module over $\Alg^\tau$:
\[
\begin{tikzpicture}
  \node at (0,0) (YIp) {$a^0$};
  \node at (4,0) (YKp) {$b^1$};
  \draw[->] (YIp) to node[above]{\lab{\tau_1+\tau_3}} (YKp);
  \draw[->] (YIp) [loop left] to node[left]{\lab{\tau_{12}}} (YIp);
  \draw[->] (YKp) [loop right] to node[right]{\lab{\tau_{23}}} (YKp);
\end{tikzpicture}
\]

Tensoring $M^\tau$ with our type \DA\ model for $\Cone(F^+_\gamma)$
gives the $\{0,1\}^2$-filtered module
\[
\begin{tikzpicture}
  \node at (0,0) (XJa) {$XJa^{01}$};
  \node at (6,0) (YKb) {$YKb^{11}$};
  \node at (1,-3.5) (Z2Ia) {$Z_2Ia^{00}$};
  \node at (0,-5) (Z1Ia) {$Z_1Ia^{00}$};
  \node at (6, -5) (YIb) {$YIb^{10}$};
  \draw[->] (Z1Ia) to node[below]{\lab{1}} (YIb);
  \draw[->] (Z1Ia) to node[left]{\lab{1}} (XJa);
  \draw[->] (Z2Ia) to node[right]{\lab{\rho_{12}}} (XJa);
  \draw[->] (Z2Ia) to node[above, sloped]{\lab{\rho_1}} (YKb);
  \draw[->] (YIb) to node[right]{\lab{\rho_1+\rho_3}} (YKb);
  \draw[->] (XJa) to node[above]{\lab{\rho_1+\rho_3}} (YKb);
  \draw[->] (Z2Ia) [loop below] to node[below]{\lab{\rho_{12}}} (Z2Ia);%
  \draw[->] (Z1Ia) [loop below] to node[below]{\lab{\rho_{12}}} (Z1Ia);%
  \draw[->] (XJa) [loop above] to node[above]{\lab{\rho_{12}}} (XJa);%
  \draw[->] (YIb) [loop below] to node[below]{\lab{\rho_{12}}} (YIb);%
  \draw[->] (YKb) [loop above] to node[above]{\lab{\rho_{23}}} (YKb);%
\end{tikzpicture}
\]
Call this
filtered module $N$.

Finally, a bounded type $A$ module for the plat solid torus is given
by
\[
\begin{tikzpicture}
  \node at (0,0) (r) {$r$};
  \node at (0,-2) (s) {$s$};
  \node at (-2,-1) (m) {$m$};
  \draw[->] (r) to node[right]{\lab{1+\rho_{23}.}} (s);
  \draw[->] (r) to node[above]{\lab{\rho_{2}}} (m);
  \draw[->] (m) to node[above]{\lab{\rho_{3}}} (s);
\end{tikzpicture}
\]
Taking the $\DT$-product of $N$ with this type $A$ module gives the
$\{0,1\}^2$-filtered complex:
\[
\begin{tikzpicture}
  \node at (0,0) (mXJa) {$mXJa^{01}$};
  \node at (5,1) (rYKb) {$rYKb^{11}$};
  \node at (6,0) (sYKb) {$sYKb^{11}$};
  \node at (1,-4) (mZ2Ia) {$mZ_2Ia^{00}$};
  \node at (0,-5) (mZ1Ia) {$mZ_1Ia^{00}$};
  \node at (6, -5) (mYIb) {$mYIb^{10}.$};
  \draw[->] (mZ1Ia) to (mXJa);
  \draw[->] (mZ1Ia) to (mYIb);
  \draw[->] (mYIb) to (sYKb);
  \draw[->] (mXJa) to (sYKb);
\end{tikzpicture}
\]

This gives a spectral sequence with $E_0$- and $E_1$-page
\begin{center}
  \begin{tabular}{cc}
    $\Field$ & $(\mathbb{F}_2)^2$\\[4pt]
    $(\mathbb{F}_2)^2$ & $\Field$,
  \end{tabular}
\end{center}
and $E_2$- through $E_\infty$-page
\begin{center}
  \begin{tabular}{cc}
    0 & $\Field$\\[4pt]
    $\Field$ & 0.
  \end{tabular}
\end{center}

In particular, the spectral sequence from the reduced Khovanov
homology to $\HFa(\Sigma(L))$ collapses at the $E_2$-page.
(This, of course, was already known, since $L$ is alternating,
so the rank of its reduced Khovanov homology coincides with its
determinant according to~\cite{Lee05:Khovanov}; see also~\cite{BrDCov}.)


\bibliographystyle{hamsalpha}\bibliography{heegaardfloer}
\end{document}